\documentclass[12pt]{article}
\usepackage{theorem}
\usepackage{amsmath}
\usepackage{amssymb}
\usepackage[mathscr]{eucal}
\usepackage{bm}
\usepackage{graphicx}
\usepackage{caption}
\usepackage{subcaption}
\usepackage{color}

\newtheorem{prop}{}[section]

{\theorembodyfont{\upshape} \newtheorem{rema}[prop]{}}
\begin{document}
\newcommand{\ben}{\begin{eqnarray}}
\newcommand{\een}{\end{eqnarray}}
\newcommand{\beq}{\begin{equation}}
\newcommand{\eeq}{\end{equation}}
\newcommand{\bsp}{\begin{split}}
\newcommand{\esp}{\end{split}}
\newcommand{\boma}[1]{{\mbox{\boldmath $#1$} }}
\newcommand{\pa}{\partial}
\newcommand{\HD}[1]{\mathbf{H}^{{#1}}_{\Sigma 0}}
\newcommand{\HM}[1]{\mathbb{H}^{{#1}}_{\Sigma 0}}
\newcommand{\CM}[1]{\mathbb{C}^{{#1}}_{\Sigma 0}}
\newcommand{\hnpps}{\mathbb{H}_{\Sigma}^{n+2}}
\newcommand{\hns}{\mathbb{H}_{\Sigma}^{n}}
\newcommand{\hnps}{\mathbb{H}_{\Sigma}^{n+1}}
\newcommand{\hnpp}{\mathbb{H}_{\Sigma 0}^{n+2}}
\newcommand{\hnz}{\mathbb{H}_{0}^{n}}
\newcommand{\hnzp}{\mathbb{H}_{0}^{n+1}}
\newcommand{\hpzp}{\mathbb{H}_{0}^{p+1}}
\newcommand{\hn}{\mathbb{H}_{\Sigma 0}^{n}}
\newcommand{\Dd}{\mathbb{D}}
\newcommand{\bDd}{\mathbf{D}}
\newcommand{\Ll}{\mathbb{L}}
\newcommand{\bL}{\mathbf{L}}
\newcommand{\hnp}{\mathbb{H}_{\Sigma 0}^{n+1}}
\newcommand{\hpp}{\mathbb{H}_{\Sigma 0}^{p+1}}
\newcommand{\hpz}{\mathbb{H}_{0}^{p}}
\newcommand{\hp}{\mathbb{H}_{\Sigma 0}^p}
\newcommand{\ler}{\mathcal{L}}
\newcommand{\tor}{\mathbf{T}}
\newcommand{\real}{\mathbb{R}}
\newcommand{\comp}{\mathbb{C}}
\newcommand{\zt}{\mathbf{Z}^3}
\newcommand{\zto}{\mathbf{Z}_0^3}
\newcommand{\bP}{\boma{\mathscr{P}}}
\newcommand{\bA}{\boma{\mathscr{A}}}
\newcommand{\Hn}{\mathcal{H}_{\Sigma 0}^n}
\newcommand{\Hnp}{\mathcal{H}_{\Sigma 0}^{n+1}}
\newcommand{\Hnpp}{\mathcal{H}_{\Sigma 0}^{n+2}}
\newcommand{\ralap}{\sqrt{-\Delta}^n}
\newcommand{\ralapP}{\sqrt{-\Delta}^p}
\newcommand{\lap}{\Delta}
\newcommand{\cald}{\mathcal{D}}
\newcommand{\calw}{\mathcal{W}}
\newcommand{\calr}{\mathcal{R}}
\newcommand{\rref}[1]{(\ref{#1})}
\newcommand{\Hinf}{\mathcal{H}_{\Sigma 0}^\infty}
\newcommand{\hinf}{\mathbb{H}_{\Sigma 0}^\infty}
\newcommand{\Hp}{\mathcal{H}_{\Sigma 0}^p}
\newcommand{\Hpp}{\mathcal{H}_{\Sigma 0}^{p+1}}
\newcommand{\Hppp}{\mathcal{H}_{\Sigma 0}^{p+2}}
\newcommand{\bb}[1]{\mathbb{{#1}}}

\providecommand{\normn}[1]{\lVert#1\rVert_n}
\providecommand{\normp}[1]{\lVert#1\rVert_{n+1}}
\providecommand{\normpp}[1]{\lVert#1\rVert_{n+2}}
\providecommand{\normnb}[1]{\lVert#1\rVert_{{n}}}
\providecommand{\normpnb}[1]{\lVert#1\rVert_{{{n}+1}}}
\providecommand{\normppnb}[1]{\lVert#1\rVert_{{n+2}}}
\providecommand{\normld}[1]{\lVert#1\rVert_{L^2}}
\providecommand{\normm}[1]{\lVert#1\rVert_{{m}}}
\providecommand{\normP}[1]{\lVert#1\rVert_{{p}}}
\providecommand{\normPp}[1]{\lVert#1\rVert_{{p+1}}}

\def\arsub{\rightsquigarrow}
\def\ua{u_a}
\def\LP{\ler}
\def\Tg{T_{G}}
\def\so{\Sigma 0}
\def\ug{\bu_G}
\def\DD{D}
\def\Q{Q}
\def\Dz{\bb{\DD}'_{0}}
\def\Dsz{\bb{\DD}'_{\Sigma 0}}
\def\HGG{\bb{H}^\G}
\def\HG{\bb{H}^\G_{\Sigma 0}}
\def\EG{{\mathscr{E}}^{\G}}
\def\bEG{{\boma{\mathscr{E}}^{\G}}}
\def\G{G}
\def\Kp{\hat{K}}
\def\Gp{\hat{G}}
\def\la{\langle}
\def\ra{\rangle}
\def\leqs{\leqslant}
\def\geqs{\geqslant}
\def\To{ {\mathbb{T} }}
\def\Td{ {\To}^d }
\def\Tt{ {\To}^3 }
\def\Zd{ \interi^d }
\def\Zt{ \interi^3 }
\def\reali{{\mathbb{R}}}
\def\complessi{{\mathbb{C}}}
\def\interi{{\mathbb{Z}}}
\def\Z{{\mathbb{Z}}}
\def\naturali{{\mathbb{N}}}
\def\vain{\rightarrow}
\def\sc{ {\scriptstyle{\bullet} }}
\def\beq{\begin{equation}}
\def\feq{\end{equation}}
\def\beqq{\begin{eqnarray}}
\def\feqq{\end{eqnarray}}
\def\barray{\begin{array}}
\def\farray{\end{array}}
\def\Tc{T_{\tt{c}}}
\def\dd{\displaystyle}
\def\Rr{{\mathcal{R}}}
\newcommand{\bphi}{\bm{u}}
\newcommand{\bu}{\mathbf{u}}
\newcommand{\bee}{\mathbf{e}}
\newcommand{\bv}{\mathbf{v}}
\newcommand{\bw}{\mathbf{w}}
\newcommand{\bz}{\mathbf{z}}
\def\vopdue{$\VP_2(f_0)\,$}
\def\leqs{\leqslant}
\def\geqs{\geqslant}
\def\mat{{\frak g}}
\def\tG{t_{\scriptscriptstyle{G}}}
\def\tN{t_{\scriptscriptstyle{N}}}
\def\TK{t_{\scriptscriptstyle{K}}}
\def\CK{C_{\scriptscriptstyle{K}}}
\def\CN{C_{\scriptscriptstyle{N}}}
\def\CG{C_{\scriptscriptstyle{G}}}
\def\CCG{{\mathscr{C}}_{\scriptscriptstyle{G}}}
\def\tf{{\tt f}}
\def\ti{{\tt t}}
\def\ta{{\tt a}}
\def\tc{{\tt c}}
\def\tF{{\tt R}}
\def\C{{\mathscr C}}
\def\P{{\mathscr P}}
\def\V{{\mathscr V}}
\def\TI{\tilde{I}}
\def\TJ{\tilde{J}}
\def\Lin{\mbox{Lin}}
\def\Hinfc{ H^{\infty}(\reali^d, \complessi) }
\def\Hnc{ H^{n}(\reali^d, \complessi) }
\def\Hmc{ H^{m}(\reali^d, \complessi) }
\def\Hac{ H^{a}(\reali^d, \complessi) }
\def\Dc{\DD(\reali^d, \complessi)}
\def\Dpc{\DD'(\reali^d, \complessi)}
\def\Sc{\SS(\reali^d, \complessi)}
\def\Spc{\SS'(\reali^d, \complessi)}
\def\Ldc{L^{2}(\reali^d, \complessi)}
\def\Lpc{L^{p}(\reali^d, \complessi)}
\def\Lqc{L^{q}(\reali^d, \complessi)}
\def\Lrc{L^{r}(\reali^d, \complessi)}
\def\Hinfr{ H^{\infty}(\reali^d, \reali) }
\def\Hnr{ H^{n}(\reali^d, \reali) }
\def\Hmr{ H^{m}(\reali^d, \reali) }
\def\Har{ H^{a}(\reali^d, \reali) }
\def\Dr{\DD(\reali^d, \reali)}
\def\Dpr{\DD'(\reali^d, \reali)}
\def\Sr{\SS(\reali^d, \reali)}
\def\Spr{\SS'(\reali^d, \reali)}
\def\Ldr{L^{2}(\reali^d, \reali)}
\def\Hinfk{ H^{\infty}(\reali^d, \KKK) }
\def\Hnk{ H^{n}(\reali^d, \KKK) }
\def\Hmk{ H^{m}(\reali^d, \KKK) }
\def\Hak{ H^{a}(\reali^d, \KKK) }
\def\Dk{\DD(\reali^d, \KKK)}
\def\Dpk{\DD'(\reali^d, \KKK)}
\def\Sk{\SS(\reali^d, \KKK)}
\def\Spk{\SS'(\reali^d, \KKK)}
\def\Ldk{L^{2}(\reali^d, \KKK)}
\def\Knb{K^{best}_n}
\def\sc{\cdot}
\def\k{\mbox{{\tt k}}}
\def\x{\mbox{{\tt x}}}
\def\QQQ{ {\textbf Q} }
\def\AAA{ {\textbf A} }
\def\gr{\mbox{gr}}
\def\sgr{\mbox{sgr}}
\def\loc{\mbox{loc}}
\def\PZ{{\Lambda}}
\def\PZAL{\mbox{P}^{0}_\alpha}
\def\epsilona{\epsilon^{\scriptscriptstyle{<}}}
\def\epsilonb{\epsilon^{\scriptscriptstyle{>}}}
\def\lgraffa{ \mbox{\Large $\{$ } \hskip -0.2cm}
\def\rgraffa{ \mbox{\Large $\}$ } }
\def\restriction{\upharpoonright}
\def\M{{\scriptscriptstyle{M}}}
\def\Fre{Fr\'echet~}
\def\I{{\mathcal N}}
\def\ap{{\scriptscriptstyle{ap}}}
\def\fiap{\varphi_{\ap}}
\def\dfiap{{\dot \varphi}_{\ap}}
\def\DDD{ {\mathfrak D} }
\def\BBB{ {\textbf B} }
\def\EEE{ {\textbf E} }
\def\GGG{ {\textbf G} }
\def\TTT{ {\textbf T} }
\def\KKK{ {\textbf K} }
\def\HHH{ {\textbf K} }
\def\FFi{ {\bf \Phi} }
\def\GGam{ {\bf \Gamma} }
\def\sc{ {\scriptstyle{\bullet} }}
\def\a{a}
\def\ep{\epsilon}
\def\c{\kappa}
\def\parn{\par \noindent}
\def\teta{M}
\def\elle{L}
\def\ro{\rho}
\def\al{\alpha}
\def\si{\sigma}
\def\be{\beta}
\def\ga{\gamma}
\def\te{\vartheta}
\def\ch{\chi}
\def\et{\eta}
\def\len{{\bf L}}
\def\Sfe{ {\bf S} }
\def\Zd{ \interi^d }
\def\Zt{ \interi^3 }
\def\Zet{{\mathscr{Z}}}
\def\Ze{\Zet^d}
\def\T1{{\textbf To}^{1}}
\def\es{s}
\def\ee{{E}}
\def\FF{\mathcal F}
\def\FFu{ {\textbf F_{1}} }
\def\FFd{ {\textbf F_{2}} }
\def\GG{{\mathcal G} }
\def\EE{{\mathcal E}}
\def\KK{{\mathcal K}}
\def\PP{{\mathcal P}}
\def\PPP{{\mathscr P}}
\def\PN{{\mathcal P}}
\def\PPN{{\mathscr P}}
\def\QQ{{\mathcal Q}}
\def\J{J}
\def\Np{{\hat{N}}}
\def\Lp{{\hat{L}}}
\def\Jp{{\hat{J}}}
\def\Pp{{\hat{P}}}
\def\Pip{{\hat{\Pi}}}
\def\Vp{{\hat{V}}}
\def\Ep{{\hat{E}}}
\def\Gp{{\hat{G}}}
\def\Kp{{\hat{K}}}
\def\Ip{{\hat{I}}}
\def\Tp{{\hat{T}}}
\def\Mp{{\hat{M}}}
\def\La{\Lambda}
\def\Ga{\Gamma}
\def\Si{\Sigma}
\def\Upsi{\Upsilon}
\def\Gam{\Gamma}
\def\Gag{{\check{\Gamma}}}
\def\Lap{{\hat{\Lambda}}}
\def\Upsig{{\check{\Upsilon}}}
\def\Kg{{\check{K}}}
\def\ellp{{\hat{\ell}}}
\def\j{j}
\def\jp{{\hat{j}}}
\def\BB{{\mathcal B}}
\def\LL{{\mathcal L}}
\def\MM{{\mathcal U}}
\def\SS{{\mathcal S}}
\def\DD{D}
\def\VV{{\mathcal V}}
\def\WW{{\mathcal W}}
\def\OO{{\mathcal O}}
\def\RR{{\mathcal R}}
\def\TT{{\mathcal T}}
\def\AA{{\mathcal A}}
\def\CC{{\mathcal C}}
\def\JJ{{\mathcal J}}
\def\NN{{\mathcal N}}
\def\HH{{\mathcal H}}
\def\XX{{\mathcal X}}
\def\XXX{{\mathscr X}}
\def\YY{{\mathcal Y}}
\def\ZZ{{\mathcal Z}}
\def\CC{{\mathcal C}}
\def\cir{{\scriptscriptstyle \circ}}
\def\circa{\thickapprox}
\def\vain{\rightarrow}
\def\salto{\vskip 0.2truecm \noindent}
\def\spazio{\vskip 0.5truecm \noindent}
\def\vs1{\vskip 1cm \noindent}
\def\fine{\hfill $\square$ \vskip 0.2cm \noindent}
\def\ffine{\hfill $\lozenge$ \vskip 0.2cm \noindent}

\makeatletter \@addtoreset{equation}{section}
\renewcommand{\theequation}{\thesection.\arabic{equation}}
\makeatother

\begin{titlepage}
{~}
\vspace{-2cm}
\begin{center}
{\huge On approximate solutions of the equations of incompressible magnetohydrodynamics}
\end{center}
\vspace{0.5truecm}
\begin{center}
{\large
Livio Pizzocchero$\,{}^a$ ({\footnote{Corresponding author}}), Emanuele Tassi  $\,{}^{b,c}$}\\
\vspace{0.5truecm}
${}^a$ Dipartimento di Matematica, Universit\`a di Milano\\
Via C. Saldini 50, I-20133 Milano, Italy\\
and Istituto Nazionale di Fisica Nucleare, Sezione di Milano, Italy \\
e--mail: livio.pizzocchero@unimi.it \\
${}^b$ Aix Marseille Univ, Univ Toulon, CNRS, CPT,  Marseille, France\\
${}^c$ Universit\'e C\^ote d'Azur,  Observatoire de la C\^ote d'Azur, CNRS,\\
Laboratoire Lagrange, France
\\
e--mail: etassi@oca.eu \\
\end{center}
\begin{abstract}
Inspired by an approach proposed previously for the incompressible
Navier-Stokes (NS) equations, we present a general framework
for the a posteriori analysis of the equations of incompressible
magnetohydrodynamics (MHD) on a torus of arbitrary dimension $d$;
this setting involves a Sobolev space of infinite order,
made of $C^\infty$ vector fields (with vanishing divergence and
mean) on the torus. Given any approximate solution of the MHD
Cauchy problem, its a posteriori analysis with the method
of the present work allows to infer
a lower bound on the time of existence of the exact solution,
and to bound from above the Sobolev distance of any order
between the exact and the approximate solution. In certain
cases the above mentioned lower bound on the time of existence
is found to be infinite, so one infers the global existence of
the exact MHD solution. We present some applications of this
general scheme; the most sophisticated one
lives in dimension $d=3$, with the ABC flow (perturbed magnetically)
as an initial datum,
and uses for the Cauchy problem a Galerkin approximate solution in $124$
Fourier modes. We illustrate
the conclusions arising in this case from
the a posteriori analysis of the Galerkin
approximant; these include the derivation
of global existence of the exact MHD solution with the
ABC datum,
when the dimensionless viscosity and resistivity are equal and
stay above an explicitly given threshold value.
\end{abstract}
\vspace{1cm} \noindent
\hfill \parn
\par \vspace{0.05truecm} \noindent \textbf{AMS 2000 Subject classifications:} 35Q35, 76W05.
\end{titlepage}
\section{Introduction}
\textbf{Magnetohydrodynamics (MHD) and the Navier-Stokes (NS)
equations.} The incompressible MHD equations are usually written
as follows (in dimensionless form):
\ben
\dot{u}= \nu \Delta u - u \sc \pa u + b \sc \pa b - \pa (p + {1 \over 2} |b|^2)~,   \label{uueq}\\
\dot{b}= \eta \Delta b - u \sc \pa b + b \sc \pa u~,  \label{bbeq}
\een
\beq \mbox{div} u = 0~, \quad \mbox{div} b = 0~. \eeq
Here: $u= u(x, t)$ and $b = b(x,t)$ are, respectively, the velocity and
magnetic field, depending
on the space variables $x = (x_1,...,x_d)$ and on time $t$, whereas
$p = p(x,t)$ is the pressure; the constants $\nu, \eta \geqs 0$ are the viscosity and resistivity.
Throughout the paper we consider
periodic boundary conditions, or, more precisely, we assume $x$ to range in the
$d$-dimensional torus $\Td := (\reali/2 \pi \interi)^d$.
Thus $u,b : \Td \times [0,T) \vain \reali^d$ (and $p: \Td \times [0,T)
\vain \reali$). In Eqs. (\ref{uueq}) (\ref{bbeq}) and in the rest
of the paper, $\pa$ stands for the gradient
and, for all (sufficiently regular) vector fields $v,w$, we indicate with $v \sc \pa w$
the vector field with components $(v \sc \pa w)_r :=
\sum_{s=1}^d v_s \partial_s w_r$; of course $\Delta := \sum_{r=1}^d \pa_{r r}$
is the Laplacian.
The space dimension $d$ is arbitrary, but we are typically interested in the case
$d=3$. For $d=3$, one can write the above equations in a more
familiar form using the identities
$- u \sc \pa b + b \sc \pa u = \mbox{rot}(u \wedge b)$
and $b \sc \partial b - \partial(|b|^2/2) = (\mbox{rot} \, b) \wedge b$
(the first one holding for all divergence free vector fields $u,b$ and
the second one valid for any vector field $b$). \parn
One can reexpress Eqs. \rref{uueq} \rref{bbeq} applying
to both sides the Leray projection $\ler$, which transforms any vector field
(on the torus) into its divergence free part; this operator annihilates gradients,
so that Eqs. \rref{uueq} \rref{bbeq} become
\ben
\dot{u}= \nu \Delta u -\ler ( u \sc \pa u) + \ler ( b \sc \pa b)~,   \label{ueq}\\
\dot{b}= \eta \Delta b -\ler ( u \sc \pa b) + \ler ( b \sc \pa u)~,  \label{beq}
\een
(and no longer contain the pressure $p$). It should be noted that
the vector field $- u \sc \pa b + b \sc \pa u$
is divergence free like $u$ and $b$, so that
$-\ler ( u \sc \pa b) + \ler ( b \sc \pa u)=
- u \sc \pa b + b \sc \pa u$. In spite of this,
for our purposes it is convenient to indicate
explicitly $\ler$ in these terms of Eq. \rref{beq};
one advantage is that, in this formulation, all
bilinear terms in Eqs. \rref{ueq} \rref{beq} involve
a unique bilinear map
\beq \P : (v,w) \mapsto \P(v,w) := - \ler(v \sc \pa w) \label{fundam} \eeq
(where $v,w: \Td \vain \reali^d$ are any two sufficiently
smooth vector fields). This ``fundamental'' bilinear map
is the same governing the NS equations of incompressible
fluids, which read
\beq \dot{u} = \nu \Delta u + \P(u,u) \label{ns} \feq
(with $u$ representing
again the velocity field, and
$\nu \geqs 0$ the viscosity; for $\nu=0$, these become the Euler equations).
\salto
\textbf{A posteriori analysis of NS approximate solutions:
a review of known results.}
The structural analogies between the MHD equations \rref{ueq}\rref{beq}
and the NS equations \rref{ns} suggest the possibility to extend
to the MHD case an approach developed in the last years for
the NS equations, allowing to infer rigorous results on
their exact solutions from the a posteriori analysis of
approximate solutions.
This a posteriori approach to the NS equations was started
in \cite{Che} \cite{DR} \cite{Rob}
and continued in a series of papers co-authored by one of us
\cite{appeul} \cite{reylarge} \cite{hyp2012} \cite{Mor15};
there are close relations between this scheme and a strategy
proposed for other nonlinear PDEs
(especially, the equations of surface growth), which has even been extended
to stochastic PDEs \cite{BlomNol} \cite{BlomRom} \cite{Nol}
\cite{BlomKam}. \parn
Let us give some more information about \cite{appeul} \cite{reylarge} \cite{hyp2012} \cite{Mor15};
here one works in a rigorous functional setting, where the exact or
approximate solutions of the NS Cauchy problem take values
in suitable Sobolev spaces of (divergence free, mean zero) vector fields on the torus $\Td$.
These Sobolev spaces are based on $L^2$, and their order is either finite \cite{appeul}
or infinite \cite{Mor15}; the case of infinite order amounts to work
in a space of $C^\infty$ vector fields. In this framework
one considers an approximate solution of the NS Cauchy problem,
i.e., a function fulfilling the NS equations with a given initial
datum up to errors affecting both the evolution equations
and the initial condition. Setting up an a posteriori analysis
centered about the Sobolev norms of the above errors,
one obtains a lower bound on the time of existence
of the exact solution of the NS Cauchy problem, and
also derives upper bounds on the Sobolev distances
between the exact and the approximate solution at any instant.
The previously mentioned lower bound on the time
of existence of the NS solution can be $+\infty$;
in this case, one concludes that the solution
of the NS Cauchy problem is defined on the
whole interval $[0,+\infty)$, i.e., it is global.
\parn
The key ingredient in the above constructions are certain
differential inequalities supplemented with suitable ``initial
value inequalities'', built up from the norms of the errors
mentioned previously; these are referred
to as the \textsl{control inequalities}. The unknowns
in these inequalities are real valued functions
of time; there is a pair of control inequalities
(a differential and an initial value inequality)
for any Sobolev order, and a solution is an upper bound on the
Sobolev distance of that order between the exact
and the approximate NS solution. The time of
existence of a solution of the control inequalities
of some basic Sobolev order also gives a lower bound on
the time of existence for the exact solution
of the NS Cauchy problem. The simplest way to
solve a pair of control inequalities is to
fulfill them as equalities: in this case we have
a differential equation and an initial condition
for an unknown real function of time,
forming what we call a \textsl{control Cauchy problem}
and possessing a unique solution. \parn
In the applications already proposed for the above scheme,
the approximate NS solutions are obtained using
the Galerkin method \cite{appeul} or
a truncated expansion with respect to a suitable quantity,
which can be the reciprocal of the
viscosity \cite{reylarge} or the time variable
\cite{hyp2012}. The fully quantitative
implementation of the a posteriori analysis requires accurate
estimates on the constants in certain
inequalities about the fundamental bilinear map \rref{fundam},
involving the Sobolev norms; rather
accurate upper bounds on these constants
have been given in \cite{Mor12b} \cite{Mor13}
\cite{Mor17}.
\salto
\textbf{Contents of the present work.} The aim of
this paper is to transfer some results on
the NS equations to the MHD case; these
results concern mainly the a posteriori
analysis of approximate solutions,
as developed in \cite{Mor15} (and in
the previous work \cite{appeul}) for NS equations. \parn
The key point in our constructions is
a formulation of the MHD equations, emphasizing strong
analogies with the setting of \cite{Mor15}
for NS equations. In few words,
the pair of functions $\bu := (u,b)$
appearing in the MHD equations
\rref{ueq} \rref{beq} is viewed as taking values
in the product of two copies
of an infinite order Sobolev space
(made of divergence free and mean zero
vector fields on $\Td$), and the cited
equations are written as
\beq \dot{\bu} = \bA \bu + \bP(\bu,\bu)~; \label{wras} \feq
here $\bA$ is the operator
$(u,b)$ $\mapsto$ $(\nu \Delta u, \eta \Delta b)$
and $\bP$ is a ``two component'' bilinear map
whose definition is suggested by the structure
of the bilinear terms in Eqs. \rref{ueq} \rref{beq}
(see Eq. \rref{debp} for the necessary details).
One can notice that the NS equations \rref{ns} are
formally converted into the MHD equations \rref{wras}
with the substitutions
\beq u \arsub \bu, \quad \nu \Delta \arsub \bA, \quad \P \arsub \bP~. \label{substi} \feq
These structural similarities
are very deep. In fact, as shown in the present paper,
some important Sobolev norm inequalities fulfilled by the
NS bilinear map $\P$ have essentially identical
counterparts for $\bP$; this is a not-so-trivial
fact, whose proof requires a minimum of effort. In addition,
some Sobolev norm inequalities for $\nu \Delta$
have counterparts for $\bA$, based on the
parameter $\mu := \min(\nu,\eta)$.
\parn
After pointing out the above structural analogies,
in the present work we consider any approximate
solution of the MHD Cauchy problem
and we analyze it a posteriori,
using ideas developed in \cite{Mor15}
for the NS equations and adapting
them to the MHD case. In this way
we derive lower bounds on the time
of existence of the exact MHD solution
and upper bounds on the Sobolev distances
(of any order)
between the exact and the approximate solution;
suitable control inequalities (conceptually
similar to those mentioned in the previous paragraph)
are developed for this purpose. In some cases,
this construction ensures the global
existence in time (i.e., a domain
$[0,+\infty)$) for
the exact solution of the MHD Cauchy problem. \parn
Some basic estimates on the exact
solution of the MHD Cauchy problem
can be obtained applying the previous scheme
to a very simple approximate solution,
namely, the zero function. The a posteriori
analysis of the zero function shows,
amongst else, that the solution of the MHD Cauchy
problem with any smooth initial
datum is global if $\mu := \min(\nu,\eta)$
is above a computable threshold value, depending
on the datum.
As examples we present these basic estimates
in space dimension $d=3$,
choosing as initial data the
Orszag-Tang vortex  and an Arnold-Beltrami-Childress (ABC) flow
with a perturbing magnetic field \cite{Min06} \cite{Car09}. \parn
A second, more sophisticated application is developed subsequently;
this uses an approximate
solution provided by the Galerkin
method (i.e., by the truncation of
the MHD equations to a finite set
of Fourier modes). For this construction
we take inspiration from
the Galerkin method for NS equations,
in the approach described by \cite{appeul}. \parn
The explicit (numerical)
construction
of the Galerkin approximants for
the MHD Cauchy problem is exemplified
in space dimension $d=3$,
assuming a common value $\mu$
for the viscosity and
the resistivity ($\nu=\eta \equiv \mu$)
and choosing the (magnetically perturbed) ABC initial datum;
the Galerkin approximant is supported by a set
of $124$ Fourier modes. The a posteriori
analysis of this approximate solution
shows, for example, that the MHD
equations \rref{ueq} \rref{beq} with
the ABC initial datum have a global
solution if $\mu$ is above a known threshold
value, determined by the Galerkin
approximant and by the control inequalities
(this estimate is sensibly
better than the one previously
mentioned for the ABC flow, based
on the zero approximate solution).
For $\mu$ below the threshold value,
our approach grants existence
of the MHD exact solution
up to a finite, explicitly computable
time.
\salto
\textbf{Organization of the paper.}
Section \ref{prelim} describes
some general facts on Sobolev
spaces on $\Td$; it also
reviews some results on
the bilinear map $\P$ of Eq. \rref{fundam},
including the inequalities mentioned before. \parn
Section \ref{mhd} discusses the MHD equations
Cauchy problem, in a setting based
on the infinite order Sobolev space mentioned
before (made of $C^\infty$ vector fields on $\Td$);
local in time existence of the exact
solution is reviewed, making reference
to the available literature (see Proposition
\ref{procau} and the discussion that accompanies it).
In the same section, we emphasize the
analogies between the NS equations \rref{ns}
and the MHD equations in the formulation
\rref{wras}, with the definition
\rref{debp} for $\bP$.
We have already mentioned that certain
Sobolev norm inequalities for the
NS bilinear map $\P$ (see Eq. \rref{fundam})
have counterparts for $\bP$: this fact
is presented in Section \ref{mhd} and
proved in Appendix \ref{appehel}, where we
also estimate certain related constants (making
reference to results of \cite{Mor12b} \cite{Mor13}
\cite{Mor17} about $\P$). Again in section \ref{mhd},
we write down some natural inequalities
for the operator $\bA$ of Eq. \rref{wras}. \parn
Section \ref{secapp} presents our general setting
for approximate solutions of the MHD Cauchy problem
and their a posteriori analysis, based on
the previously mentioned control inequalities.
This framework is applied in Section
\ref{secappan} to the zero approximate solution,
and in Section \ref{galegen} to general Galerkin
approximants. Finally, in Section
\ref{galespec} we construct (in dimension $3$) the previously
mentioned Galerkin approximate solution
with $124$ Fourier modes for the perturbed ABC
initial datum, and describe the
results on the exact solution
arising from its a posteriori analysis.
\salto
\textbf{Notice.}
After the acquisition of the structural analogies
between the NS equations \rref{ns} and
the MHD equations \rref{wras}, the main propositions about
the MHD approximate solutions can be derived by a simple
translation of similar propositions proved in \cite{Mor15} for the
NS approximate solutions; essentially, one applies
the ``correspondence principle" \rref{substi}. To some
extent, a similar remark also applies to the analysis of the Galerkin
approximation for the MHD equations; many results on
this subject are obtained translating via \rref{substi}
the analysis of the Galerkin method performed in
\cite{appeul} for the NS equations. \parn
In spite of this, in writing the present paper we have decided
to give explicitly the above mentioned ``translations''
to the MHD framework, even at the price of textual similarities
with the corresponding statements of \cite{appeul} \cite{Mor15} on NS equations.
This choice makes the present paper self-contained, a feature that we think
could be useful since the a posteriori analysis
of the MHD approximate solutions is (to the best of our
knowledge) an essentially new subject. In any case,
the connections of the present results
with \cite{appeul} \cite{Mor15} are indicated explicitly whenever they occur.
\section{General functional setting}  \label{prelim}
\noindent
We work in any space dimension
\beq d \in \{2,3,...\}~; \eeq
for $a, b \in \complessi^d$ we write $a \sc b := \sum_{r=1}^d a_r b_r$.
We often use the lattice $\Zd$, where $\interi$ is the set of integers; denoting
with $0$ its zero element, we write $\Zd_{0} := \Zd \setminus \{0 \}$.
\par \noindent
\textbf{Sobolev spaces on the torus.} Throughout
the paper we stick rather closely to the functional setting
adopted in \cite{appeul} \cite{Mor15} for the NS equations on $\Td$;
for convenience of the reader, let us re-propose here the basic
function spaces involved in this setting. \parn
First of all, we write
$\Dd'$ for the space of $\reali^d$-valued distributions on $\Td$
(distributional vector fields). Each $v \in \Dd'$ has weakly
convergent Fourier expansion $ v=\sum_{k \in \Zd} v_k e_k$, where
$e_k(x) := (2 \pi)^{-d/2} \mathrm{e}^{i k \sc x}$ and $v_k =
\overline{v_{-k}} \in \complessi^d$ are the Fourier coefficients.
The spaces of divergence free or zero mean distributional vector
fields and their intersection are \beq \Dd'_{\Sigma} :=\{ v \in
\Dd'~|~\mbox{div} v = 0~\} = \{ v \in \Dd'~|~k \sc v_k = 0\,
\mbox{for $k \in \Zd$} \}~; \label{dsig} \eeq \beq \Dd'_{0} :=\{ v
\in \Dd'~|~\int_{\Td} v\, d x = 0~\} = \{ v \in \Dd'~|~v_0 = 0 \}~;
\quad \Dd'_{\Sigma 0} := \Dd'_{\Sigma} \cap \Dd'_{0}~. \eeq
Let us consider the space $\Ll^2$ of square integrable vector fields on $\Td$,
and its standard inner product $\la~|~\ra_{L^2}$. For
any $p \in \reali$, we define the Sobolev space $\HM{p}$ of divergence free, zero mean
vector fields on $\Td$ as \beq \label{sobol}
\begin{split}
\HM{p} & :=\{ v \in \Dd'~|~\mbox{div} v=0,~\int_{\Td} v \, d x = 0,~
\sqrt{-\Delta}^{\,p} v \in \Ll^2 \}  \\
& = \{ v \in \Dd'~|~
k \sc v_k = 0,~v_0 = 0,~ \sum_{k \in \Zd_0}  \vert k \vert^{2 p} \vert v_k \vert^2 < + \infty \}~.
\end{split}
\eeq
(Here $\Delta$ is the Laplacian; the fractional power $\sqrt{-\Delta}^p$ is
defined by $(\sqrt{-\Delta}^p v)_k =|k|^p v_k$, as suggested
by the obvious Fourier representation $(-\Delta v)_k = |k|^2 v_k$).
$\HM{p}$ is a real Hilbert space with the inner product
\beq \label{inpn}
\la v \vert w\ra _p :=\la{\sqrt{- \Delta}}^{\,p} v \vert {\sqrt{- \Delta }}^{\,p} w \ra _{L^2}
=\sum_{k \in \Zd_0} \vert k \vert^{2 p} \bar{v}_k \sc w_k
\eeq
and the induced norm
\beq  \label{normv}
\| v \|_p :=\sqrt{\la v \vert v \ra _p}~; \eeq
of course, for $p=0$ we have
\beq \HM{0} = \Dd'_{\Sigma 0} \cap \Ll^2~,~~\la~|~\ra_{0} = \la~|~\ra_{L^2}~. \label{accazero} \feq
From now on we indicate with $\hookrightarrow$ a continuous imbedding.
With this notation, for $p \geqs q$ we have $\HM{p} \hookrightarrow \HM{q}$
(and, more quantitatively, $\|~\|_p \geqs \|~\|_q$).
Now we introduce, analogously to Ref. \cite{Mor15},
the infinite order Sobolev space
\beq \HM{\infty} := \cap_{p \in \reali} \HM{p}~. \eeq
This carries the complete topology induced by the infinitely many norms
$\|~\|_{p}$ ($p \in \reali$); indeed, this family of norms is equivalent
to the countable family $\|~\|_{p}$ ($p=0,1,2,...$), so $\HM{\infty}$ is
a \Fre space.
For $k \in \naturali \cup \{ \infty \}$ we consider the space
\beq \CM{k}(\Td) := \{ v \in C^k(\Td,\reali^d)~|~~
\mbox{div} v = 0,~ \int_{\Td} v\, d x=0~\}~, \feq
which is a Banach space for $k < \infty$ and a \Fre space for $k=\infty$, with
the usual sup norms of the derivatives of all involved orders.
Let $h, k \in \naturali$, $p \in \reali$;
then $\CM{h} \hookrightarrow \HM{p}$ if $h \geqs p$ and, by the Sobolev lemma,
$\HM{p} \hookrightarrow \CM{k}$ if $p > k + d/2$. These imbeddings imply
\beq \HM{\infty} = \CM{\infty} \feq
(equality as topological vector spaces).
\salto
\textbf{Laplacian.} Let us consider the Laplacian $\Delta : \Dd' \vain \Dd'$;
from the Fourier representation $(\Delta v)_k = -|k|^2 v_k$
we readily infer the following: for each real $p$ and $v \in \HM{p+2}$, one
has $\Delta v \in \HM{p}$ and
\beq \| \Delta v \|_p = \| v \|_{p+2}~, \label{lap1} \feq
\beq \la \Delta v | v \ra_p = - \| v \|^2_{p+1} \leqs - \| v \|^2_p~. \label{lap2} \feq
Using Eq.\rref{lap1}, one infers that $\Delta$ is continuous from $\HM{p+2}$ to $\HM{p}$
for each real $p$,
and from $\HM{\infty}$ to $\HM{\infty}$.
\salto
\textbf{Leray projection.} This is the map
\beq \ler : \Dd' \vain \Dd'_{\Sigma}~,\quad
v \mapsto \ler v~\mbox{such that
$(\ler v)_k = \ler_k v_k$ for $k \in \Zd$}~; \label{defler} \eeq
here $\ler_k$ is the orthogonal projection of $\complessi^d$ onto $k^\perp = \{ a \in
\complessi^d~|~k \sc a = 0 \}$ (so that $\ler_k c = c
- (k \sc c) k/|k|^2$ and $\ler_0 c = c$, for $k \in \Zd_0$ and
$c \in \complessi^d$).
One proves that $\ler \Dd' = \Dd'_{\Sigma}$, $\ler \Dd'_{0} = \Dd'_{\Sigma 0}$,
$\ler \Ll^2 = \Dd'_{\Sigma} \cap \Ll^2$.
\salto
\textbf{Fundamental bilinear map.}
Let us consider two vector fields $v,w \in \Dd'$ such that
\beq v \in \Ll^2, \quad \partial_s w \in \Ll^2~ \mbox{for $s=1,...,d$}; \label{vw} \feq
then $v \sc \pa w$ belongs to the space $\Ll^1$ of
integrable vector fields on $\Td$. The bilinear map
sending $v, w$ as in \rref{vw} into
\beq \P(v, w) := -\ler(v \sc \pa w) \in \ler \Ll^1 \label{defpi} \feq
will be referred to as the ``fundamental bilinear map''.
In terms of Fourier components, we have
\beq (v \sc \pa w)_k=\frac{i}{(2 \pi)^{d/2}} \sum_{h \in \Zd} [v_h
\sc (k-h)] w_{k-h}. \label{eq1} \eeq
for all $k \in \Zd$; this implies that the $k$-th Fourier component of $\P(v,w)$ is
\beq \P_k(v,w) =- \frac{i}{(2 \pi)^{d/2}} \sum_{h \in \Zd} [v_h
\sc (k-h)] \LL_k w_{k-h}. \label{eq2} \eeq
(where, as in the previous paragraph, $\LL_k$ indicates the orthogonal projection of $\complessi^d$
onto $k^{\perp}$). Of course, in Eqs. \rref{eq1} \rref{eq2} the sum
over $\Zd$ can be replaced with a sum over $\Zd \setminus \{k\}$; moreover, if
$v$ has mean zero we can sum over the set $\Zd \setminus \{0,k\}$, hereafter denoted
with $\Zd_{0 k}$. \parn
To go on let us remark that, for $v, w$ as in \rref{vw},
\beq \la v \sc \partial w | w \ra_{L^2} =  \la \P(v, w)  | w \ra_{L^2} = 0 \quad
\mbox{if $v \sc \partial w \in \Ll^2$ and $\mbox{div} v = \mbox{div} w = 0$}
\label{recall0} \eeq
(this follows, e.g., from Eq. (1.8) and Lemma 2.3 of \cite{Mor12b}; note that
$v \sc \partial w \in \Ll^2$ implies $\P(v,w) \in \Ll^2$).
We now add much more regularity.
Let $n, p$ denote two real numbers; it is known that
\beq p > d/2,~ v \in \HM{p},~ w \in \HM{p+1} \quad \Rightarrow \quad \P(v,w) \in \HM{p} \label{known} \feq
and that there are constants $K_{p n}$, $G_{p n}$ $\in (0,+\infty)$ such that the following holds:
\beq \| \P(v, w) \|_p \leqs {1 \over 2} K_{p n} ( \| v \|_p \| w \|_{n+1} + \| v \|_n \| w \|_{p+1})
\label{basineqa} \feq
$$ \qquad \mbox{for $p \geqs n > d/2$, $v \in \HM{p}$, $w \in \HM{p+1}$}~, $$
\beq | \la \P(v, w) | w \ra_p | \leqs
{1 \over 2} G_{p n} (\| v \|_p \| w \|_n + \| v \|_n \| w \|_p)\| w \|_p \label{katineqa} \feq
$$ \mbox{for $p \geqs n > d/2 + 1$, $v \in \HM{p}$, $w \in \HM{p+1}$}~. $$
Of course, with $p=n$ and $K_{p} := K_{p p}$, $G_p := G_{p p}$ the above
inequalities become
\beq \| \P(v, w) \|_p \leqs K_p \| v \|_p \| w \|_{p+1} \qquad \mbox{for $p > d/2$,
$v \in \HM{p}$, $w \in \HM{p+1}$}~, \label{basineq} \feq
\beq | \la \P(v, w) | w \ra_p | \leqs G_p \| v \|_p \| w \|^2_p
\qquad \mbox{for $p > d/2 + 1$, $v \in \HM{p}$,
$w \in \HM{p+1}$}~, \label{katineq} \feq
Statements \rref{known} \rref{basineq} indicate that $\P$ maps continuously
$\HM{p} \times \HM{p+1}$ to $\HM{p}$. The same statements can be used
in an obvious way to prove that $\P$ maps continuously $\HM{\infty}
\times \HM{\infty}$ to $\HM{\infty}$. \parn
Eq.\,\rref{basineq} will be referred to as the ``basic''
inequality for $\P$, since it is closely related to the standard norm inequalities about
multiplication in Sobolev spaces; Eq. \rref{katineq}
was established by Kato \cite{Kat72} for integer $p$,
and generalized to noninteger cases in \cite{CoFo}; it will be
referred to as the ``Kato inequality''.
Eqs. \rref{basineqa} \rref{katineqa} are ``tame'' generalizations
(in the Nash-Moser sense) of the basic and Kato inequalities
for $\P$ (for some inequalities very similar to \rref{katineqa},
see \cite{Tem} \cite{BKM} \cite{RSS}).
\parn
The inequalities \rref{basineqa} \rref{katineqa} and
the related constants were discussed in \cite{Mor17},
generalizing previous results of \cite{Mor12b} \cite{Mor13}
on the special case $p=n$. The analysis of \cite{Mor17}
shows that the cited relations are fulfilled with
\beq
K_{pn} =\frac{1}{(2 \pi)^{d/2}}\sqrt{\sup_{k \in \Zd_0} \mathcal{K}_{pn} (k)}~,
\label{dkpn} \eeq
\beq
G_{pn} =\frac{1}{(2 \pi)^{d/2}}\sqrt{\sup_{k \in \Zd_0} \mathcal{G}_{pn} (k)}~,
\label{dgpn} \eeq
where $\mathcal{K}_{pn}$, $\mathcal{G}_{pn}$ $: \interi^d_0 \to (0,+\infty)$ are the functions defined by
\beq \KK_{p n}(k) := 4 |k|^{2 p}
\sum_{h \in \Zd_{0 k}} {\Q^2_{h, k - h} \over (|h|^{p} |k-h|^{n} + |h|^n |k-h|^{p})^2}~;
\label{kknd} \feq
\beq
\mathcal{G}_{pn}(k):=
4 \sum_{h \in \Zd_{0 k}}\frac{(\vert k \vert^p - \vert k-h\vert^p)^2
Q^2_{h, k-h}}{(\vert h \vert^p\vert k-h\vert^{n-1}+\vert h\vert^n \vert k-h\vert^{p-1})^2}~.
\label{ggnd} \eeq
Here $\Zd_{0 k} := \Zd \setminus \{0 , k\}$ (as already defined),
and for all $q, r, h, \ell \in \reali^d \setminus \{0 \}$ we stipulate the
following:
\beq \te_{q r} := \mbox{convex angle between $q,r$} \quad (\te_{q r} \in [0,\pi])~, \feq
\beq \Q_{h \ell} := \left\{ \barray{ll} \sin \te_{h \ell}  & \mbox{if $d \geqs 3$}\,, \\
\sin \te_{h \ell} \cos \te_{h + \ell, \ell} & \mbox{if $d = 2$}
\farray \right.
\label{eqnorm}
\feq
({\footnote{Of course, $\cos \te_{q r} = \dd{q \sc r \over |q| |r|}$ and $\sin \te_{q r} =
\sqrt{1 - \dd{(q \sc r)^2 \over |q|^2 |r|^2}}$.
In the definition of $\Q_{h \ell}$ for $d \geqs 3$,
$\te_{h + \ell, \ell}$ is meant to indicate any angle in $[0,\pi]$
if $h + \ell=0$; the chosen value is immaterial, since
in this case $\te_{h \ell}=\pi$ and $\sin \te_{h \ell}=0$.
The coefficient
$\Q_{h \ell}$ arises in \cite{Mor17} as the norm of a certain bilinear
map acting on vectors of $\reali^d$, a fact not relevant for
our present purposes.}}).
As examples for later use, let us give explicit values for the constants
$K_{p p} \equiv K_p, G_{p p} \equiv G_p$ and $G_{p n}$ in space dimension $3$,
for some values of $p,n$ of interest for the sequel.
From \cite{Mor12b} \cite{Mor13} \cite{Mor17}  we know that we can take
({\footnote{The value of $K_3$ employed here is taken from \cite{Mor17}; this value slightly improves
the estimate given previously in \cite{Mor13}.}})
\beq K_3 = 0.320, ~ G_3 = 0.438, ~ K_5 = 0.657, ~ G_5 = 0.749, ~ G_{5 3} = 1.26
\quad (d=3). \label{costns} \feq
\salto
\section{The MHD Cauchy problem}
\label{mhd}
\textbf{Formulation of the problem.}
Let us choose two parameters
\beq \nu, \eta \in [0,+\infty)~, \feq
that we call the viscosity and the resistivity following the Introduction.
Moreover, we fix a couple of initial data
\beq   \label{indata}
u_0 , b_0 \in \HM{\infty}~.
\eeq
The MHD Cauchy\ problem with viscosity, resistivity and initial data as above
reads:
\beq \mbox{Find} \, u, b \in C^\infty([0,T), \HM{\infty})~\mbox{(with $T \in (0,+\infty]$) such that}~ \label{causm} \eeq
$$
\dot{u}=  \nu \Delta u + \P(u,u) - \P( b, b)~, $$
$$ \dot{b}= \eta \Delta b + \P(u, b) - \P(b,u) ~. $$
$$ u(0) = u_0, ~~b(0) = b_0 $$
(with $\P$ as in Eq. \rref{defpi}).
In the above, one recognizes Eqs. \rref{ueq} \rref{beq} of the Introduction;
the length $T$ of the time interval considered in \rref{causm} is unspecified, and depends on $(u,b)$.
\salto
\textbf{A reformulation of the previous setting for MHD.}
For the sake of brevity, let
\beq \bDd' := \Dd'\times \Dd'~, \qquad
\bDd'_{\Sigma 0} := \Dd'_{\Sigma 0} \times \Dd'_{\Sigma 0} \feq
\beq \bL^2 := \Ll^2 \times \Ll^2~; \feq
the last space is a Hilbert space with the inner product
\beq \la \bv | \bw \ra_{L^2} := \la v | w \ra_{L^2} + \la b | c \ra_{L^2}
\quad \mbox{for $\bv = (v,b)$, $\bw = (w,c) \in \bL^2$}~. \feq
By comparison with the Cauchy problem \rref{causm}, we see that this
involves the linear operator
\beq \bA : \bDd' \rightarrow \bDd'~, \qquad
\bv := (v,b) \mapsto \bA \bv := (\nu \Delta v , \eta \Delta b )
\label{deba} \feq
and the bilinear map
\beq \bv=(v,b), \bw=(w,c) \mapsto
\bP (\bv,\bw):=(\P(v,w)-\P(b,c) , \P(v,c)-\P(b,w))~. \label{debp} \eeq
The largest domain on which $\bP$ is well defined is formed by the
pairs $(\bv, \bw)$ as above with $v, b \in \Ll^2$ and $\partial_s w, \partial_s c \in \Ll^2$
for $s \in \{1,...,d\}$; $\bP$ maps this domain to $\LL \Ll^1 \times \LL \Ll^1$. \par
To go on, for any real $p$ we
introduce the Hilbert space
\beq
\HD{p} := \HM{p} \times \HM{p},
\eeq
equipped with the inner product
\beq   \label{inpbn} \la\bv \vert \bw\ra _{p}:=\la v \vert w\ra _p + \la b
\vert c \ra _p, \eeq
for $\bv=(v,b), \bw=(w,c)$. From this inner product
we also derive the norm
\beq \lVert \bv \rVert_{{p}}:=\sqrt{\la \bv \vert \bv\ra _{{p}}}=\sqrt{
\| v \|_p^2 + \| b \|^2_p}. \eeq
We also set
\beq
\HD{\infty} := \HM{\infty} \times \HM{\infty}~;
\eeq
this is a \Fre space with the infinitely many norms
$\|~\|_p$ ($p \in \reali$ or, equivalently, $p = 0,1,2,...$).
\parn
Keeping in mind Eqs. \rref{lap1} \rref{lap2}, one readily obtains
the following: for each real $p$ and $\bv = (v,b) \in \HD{p+2}$, one has
$\bA \bv \in \HD{p}$, and
\beq \| \bA \bv \|_p = \sqrt{ \nu^2 \| v \|^2_{p+2} + \eta^2 \| b \|^2_{p+2}}
\leqs \mu \| \bv \|_{p+2}~, \label{ba1} \feq
\beq {~} \hspace{-0.5cm} \la \bA \bv | \bv \ra_p = - \nu \| v \|^2_{p+1} - \eta \| b \|^2_{p+1}
\leqs - \mu \| \bv \|^2_{p+1}
\leqs - \mu \| \bv \|^2_p~, \label{ba2} \feq
where we have set
\beq \mu := \min(\nu,\eta)~. \label{demu} \feq
Eq. \rref{ba1} implies that $\bA$ is continuous
from $\HD{p+2}$ to $\HD{p}$ for each real $p$, and
from $\HD{\infty}$ to $\HD{\infty}$. In addition, due to the properties
of $\P$ reviewed in the previous section, $\bP$
maps continuously $\HD{p} \times \HD{p+1}$ to $\HD{p}$ for
each $p > d/2$, and $\HD{\infty} \times \HD{\infty}$ to $\HD{\infty}$. \parn
Let $\bv = (v,b)$, $\bw = (w,c)$,
with $v,b,w,c,\partial_s w, \partial_s c \in \Ll^2$, $\mbox{div} v = \mbox{div} b = \mbox{div} w = \mbox{div} c =0$
and $v \sc \partial w, b \sc \partial c, v \sc \partial c, b \sc \partial w \in \Ll^2$; then, using
Eq. \rref{recall0} one proves that
({\footnote{in fact
$$ \la \bP(\bv,\bw)|\bw \ra_{L^2} = \la \P(v,w) | w \ra_{L^2} - \la \P(b,c) | w \ra_{L^2} +
\la \P(v,c) | c \ra_{L^2} - \la \P(b,w) | c \ra_{L^2}~. $$
But $\la \P(v,w) | w \ra_{L^2}=0$ and $\la \P(v,c) | c \ra_{L^2}=0$ due to \rref{recall0}; moreover
$\la \P(b,c) | w \ra_{L^2} + \la \P(b,w) | c \ra_{L^2} = (1/2) \la \P(b,c+w)|c + w \ra_{L^2}
- (1/2) \la \P(b,c - w)|c - w \ra_{L^2} = 0$, where the first equality follows from
the bilinearity of $\P$, $\la~|~\ra_{L^2}$ and the second equality relies again on \rref{recall0}.}})
\beq \la \bP(\bv,\bw)|\bw \ra_{L^2} =  0~. \label{brecall0} \feq
For the sequel of this paper, it is essential to point out that the map $\bP$
fulfills the following inequalities, containing suitable constants $\Kp_{p n}$ and $\Gp_{p n}$
and for the rest structurally identical to
the inequalities \rref{basineqa} \rref{katineqa} for $\P$:
\beq
\begin{split}
& \| \bP (\bv, \bw) \|_p \leqs \frac{1}{2} \Kp_{pn} (\| \bv \|_p \| \bw \|_{n + 1}+\| \bv \|_n \| \bw \|_{p+1} ) \\
& \mbox{for $p\geqs n > d/2$, \,
$\bv \in \HD{p}$, \,  $\bw \in \HD{p+1}$}~; \label{basic2np}
\end{split}
 \eeq
\beq  \label{katop2np}
\begin{split}
&\vert \la \bP ( \bv,\bw) \vert \bw \ra _{{p}}\vert
\leqs  \frac{1}{2} \Gp_{pn}( \normP{\bv} \normnb{\bw}+
\normnb{\bv}\normP{\bw})\normP{\bw}\\
& \mbox{for  $p \geqs n > d/2 +1$, \,
$\bv \in \HD{p}$, \, $\bw \in \HD{p+1}$}~.
\end{split}
\eeq
We refer to Appendix \ref{appehel} for the derivation of Eqs. \rref{basic2np} and
\rref{katop2np} from Eqs. \rref{basineqa} and \rref{katineqa}, respectively
(such a derivation is not so obvious, especially in the case
of \rref{katop2np}). This appendix shows that the constants in
\rref{basic2np} \rref{katop2np} can be taken as follows:
\beq \Kp_{p n} := \sqrt{2} \, K_{p n}~,  \label{weknowk} \feq
\beq \Gp_{p n} := \sqrt{2} \, G_{p n}~, \label{weknowg} \feq
where $K_{p n}, G_{p n}$ are constants fulfilling \rref{basineqa} \rref{katineqa}
(these could be taken as in Eqs. \rref{dkpn} \rref{dgpn}, respectively).
Of course, with $p=n$ and $\Kp_{p} := \Kp_{p p}$, $\Gp_{p} := \Gp_{pp}$ we get
\beq
\begin{split}
& \| \bP (\bv, \bw) \|_p \leqs \Kp_p \| \bv \|_p \| \bw \|_{p + 1} \\
 &\mbox{for $p > d/2$, \, $\bv \in \HD{p}$, \, $\bw \in \HD{p+1}$}~; \label{basic2}
 \end{split}
 \eeq
\beq  \label{katop2}
\begin{split}
&\vert \la \bP ( \bv,\bw) \vert \bw \ra _{{p}}\vert \leqs \Gp_p \| \bv \|_p \| \bw\|^2_p \\
& \mbox{for $p > d/2 + 1$, \,
$\bv \in \HD{p}$, \, $\bw \in \HD{p+1}$}~;
\end{split}
\eeq
As an example, let us consider the case of space dimension $d=3$
and give for later use the explicit values for the constants
$\Kp_{p p} \equiv \Kp_p, \Gp_{p p} \equiv \Gp_p$ and $\Kp_{p n}, \Gp_{p n}$
for some values of $p,n$. Taking the values of $K_3$, etc. reported
in Eq. \rref{costns}, multiplying each one of these values by $\sqrt{2}$
and rounding up the results to three digits, we conclude that we can take
\beq \Kp_3 = 0.453~, ~ \Gp_3 = 0.620, ~ \Kp_5 = 0.930, ~ \Gp_5 = 1.06, ~ \Gp_{5 3} = 1.79
\quad (d=3)~.
\label{k3g3} \feq
Let us return to the case of any space dimension $d \geqs 2$.
With the previous notations, for $\bu_0 = (u_0, b_0) \in \HD{\infty}$, we
can rephrase as follows the Cauchy problem \rref{causm}:
$$ \mbox{Find $\bu = (u,b) \in C^\infty([0,T), \HD{\infty}  )$ ($T \in (0 , + \infty]$) such that} $$
\beq
\dot{\bu} = \bA \bu + \bP ( \bu, \bu), \qquad
\bu (0)  =  \bu_0~.  \label{cmhd} \eeq
This formulation makes evident the analogies with the
NS Cauchy problem $\dot{u} = \nu \Delta u + \P(u,u)$, $u(0)=u_0$,
on the grounds of a ``correspondence principle"
already mentioned in the Introduction (see Eq. \rref{substi}).
\salto
\textbf{Local existence and uniqueness results for the Cauchy problem.}
The incompressible MHD Cauchy problem has been extensively discussed in the literature
in appropriate functional settings, not necessarily coinciding with
ours. To our knowledge, the case $\nu,\eta > 0$
was first studied in \cite{Lions}; a subsequent, influential
work on the same case is \cite{Ser83}. Reference \cite{Sch88}
first treated the case $\nu = \eta= 0$ which is technically
harder, proving local existence and uniqueness and
deriving a blow-up criterion by means of techniques
which in fact work for arbitrary $\nu, \eta \geqs 0$.
\parn
Paper \cite{Sch88} considers MHD on $\reali^d$,
\cite{Lions} works on a domain in $\reali^d$, \cite{Ser83}
also considers the case of $\Td$
({\footnote{To be precise, in \cite{Lions} \cite{Ser83} $d$ is
$2$ or $3$.}});
the reformulation of the main results from
\cite{Lions} \cite{Sch88} in the case of $\Td$ is straightforward.
Another feature of the three cited works is that they consider
(strong) solutions of the Cauchy problem taking values
in Sobolev spaces of finite order. However, the
adaptation of their results to a $C^\infty$
framework is obtained by standard
arguments, as indicated explicitly in
\cite{Sch88};
the infinite
order Sobolev spaces $\HD{\infty}$ considered
here allow a precise definition
of the $C^\infty$ framework.
({\footnote{
A similar situation occurs for the incompressible
NS equations; the arguments to extend existence
theorems and blowup criteria
from a finite order to an infinite order
Sobolev setting are reviewed,
e.g., in Appendix B
of \cite{Mor15}.}}) \parn
Summing up we can refer to the existing literature, and
especially to \cite{Sch88}, to account for the following statement:
\vfill \eject \noindent
\begin{prop}
\textbf{Proposition.}
\label{procau}
For all $\nu, \eta \geqs 0$, $\bu_0 = (u_0, b_0) \in \HD{\infty}$, the following holds. \parn
i) Problem \rref{causm} (or \rref{cmhd}) has a unique maximal (i.e., unextendable)
solution $\bu = (u,b)$, for suitable $T \in (0,+\infty]$; any other solution
is a restriction of the maximal one. \parn
ii) (Blow-up criterion). If $T < + \infty$, for any $n > d/2 + 1$ one has
\beq \limsup_{t \vain T^{-}} \| \bu(t) \|_n = + \infty~. \label{blowup} \eeq
\end{prop}
\noindent
For completeness, let us add two remarks: \parn
(i) There exist blow-up conditions finer than \rref{blowup},
and similar to the Beale-Kato-Majda criterion for incompressible NS equations.
For simplicity, let us consider the case of space dimension $d=3$. In \cite{Cafl},
the following criterion was derived: if $T < + \infty$, then
\beq \int_{0}^{T} d t (\| \, \mbox{rot} u(t) \|_{L^{\infty}} + \| \mbox{rot} b(t) \|_{L^{\infty}})
= + \infty~. \feq
Let us note that, for $n > d/2 + 1 = 5/2$, the Sobolev imbedding gives
$\| \mbox{rot} \, u(t) \|_{L^{\infty}} + \| \mbox{rot} \, b(t) \|_{L^{\infty}}
\leqs \mbox{const.} \| \bu(t) \|_n$; thus if $T < + \infty$ we also have
$\int_{0}^{T} d t \| \bu(t) \|_n = + \infty$,
which of course implies \rref{blowup}. \parn
(ii) The local existence (of strong solutions) for
the MHD Cauchy problem in Sobolev spaces of \textsl{minimal}
order is obviously outside the scope of this paper; however
we would mention that the problem is especially hard
when $\nu$ or $\eta$ vanish, and that such cases
have been treated only in recent times \cite{Fan} \cite{Fef14} \cite{Fef17}.
\salto
\textbf{Energy balance law.} This is a well known fact, that we review just for
convenience. Let us consider the (maximal) solution $\bu = (u,b)$ of
the Cauchy problem \rref{causm}\rref{cmhd}; for all $t$ in its domain $[0,T)$, the squared norm
\beq \| \bu(t) \|^2_{L^2} = \| u(t) \|^2_{L^2} + \| b(t) \|^2_{L^2} \feq
represents (twice) the total energy of the system at time $t$.
\begin{prop}
\textbf{Proposition.} One has
\beq {d \over d t} \| \bu \|^2_{L^2} = 2 \la \bA \bu | \bu \ra_{L^2} \leqs - 2 \mu \| \bu \|^2_{L^2} \label{thelast} \feq
(with $\mu$ as in Eq. \rref{demu}).
This implies the following, for $t \in [0,T)$:
\beq \| \bu(t) \|_{L^2} \left\{ \barray{ll} = \| \bu_0 \|_{L^2} & \mbox{if $\nu=\eta=0$,} \\
\leqs \| \bu_0 \|_{L^2} e^{- \mu t} & \mbox{for all $\nu,\eta \geqs 0$}. \farray \right.
\label{imp} \feq
\end{prop}
\textbf{Proof.}
Writing $\| \bu \|^2_{L^2} = \la \bu | \bu \ra_{L^2}$, taking the $t$ derivative and
using Eqs. \rref{cmhd} \rref{brecall0} \rref{ba2} we get
$$ {d \over d t} \| \bu \|^2_{L^2} = 2 \la {d \bu \over d t} | \bu \ra_{L^2}
=  2 \la \bA \bu | \bu \ra_{L^2} + 2 \la \bP(\bu,\bu) | \bu \ra_{L^2} =
2 \la \bA \bu | \bu \ra_{L^2} \leqs - 2 \mu \| \bu \|^2_{L^2}~. $$
This proves Eq. \rref{thelast}.
If $\nu=\eta=0$ we have $\bA=0$, and Eq. \rref{thelast} yields the first statement in \rref{imp}.
For all $\nu,\eta \geqs 0$, Eq. \rref{thelast} gives as well the second statement in \rref{imp}.  \fine
\section{Approximate solutions of the Cauchy problem for incompressible MHD}  \label{approxmhd}
\label{secapp}
\noindent
The purpose of this section (and of the subsequent
Section \ref{secappan}) is to convert
to the MHD case the framework developed in
\cite{Mor15} for the approximate solutions
of the NS Cauchy problem; concerning
this construction and its similarities with the cited work, we recall the
notice at the end of the Introduction.
\parn
From here to the end of the present section we fix a viscosity, a resistivity and a MHD initial datum, namely
\beq \nu, \eta \in [0,+\infty)~, \qquad \bu_0 = (u_0, b_0) \in \HD{\infty}~. \eeq
The definition, the lemma and the proposition which follow correspond to
Definition 4.1, Lemma 4.2 and Proposition 4.3 of \cite{Mor15} on NS equations; this remark
applies as well to the related proofs.
\begin{prop}
\textbf{Definition.}
\label{deferror}
An approximate solution of the Cauchy problem (\ref{causm}) (or \rref{cmhd})
is any map $\bu_a = (u_a , b_a) \in  C^1([0, T_a), \HD{\infty} ) $,
with $T_a \in (0,+\infty]$. Given a map of this kind, we use the following terminology: \parn
(i) The differential error of $\bu_a$ is
\beq
\bee(\bu_a) :=\frac{d \bu_a}{dt} - \bA \bu_a - \bP (\bu_a , \bu_a)  \, \in C([0, T_a) , \HD{\infty}).
\label{deferrsm}
\eeq
A differential error estimator of order $p \in \reali$ for $\bu_a$ is a
function $\epsilon_p \in C( [0, T_a) , [0, + \infty))$ such that
\beq
\lVert \bee(\bu_a)(t) \rVert_{{p}} \leqs \epsilon_p (t) \qquad \text{for} \quad t \in [0, T_a).
\label{estiepsm}
\eeq
(ii) The datum error of $\bu_a$ is
\beq \bu_a (0) - \mathbf{u}_0  \, \in \HD{\infty}. \eeq
A datum error estimator of order $p \in \reali$ for $\bu_a$ is a real number $\delta_p \geqs 0$ such that
\beq \lVert \bu_a (0) - \mathbf{u}_0 \rVert_{{p}} \leqs
\delta_p. \label{estdatumsm} \eeq
(iii)
A growth estimator of order order $p \in \reali$ for $\bu_a$ is a function
$\cald_p \in C([0, T_a) , [0, + \infty))$ such that
\beq
\lVert \bu_a (t) \rVert_{p} \leqs \cald_p (t), \qquad \mbox{for} \, t \in [0, T_a).
\label{dednsm} \eeq
\end{prop}
From here to the end of the section, we assume the following: \parn
i) $\bu =(u,b)$ is the maximal (exact) solution of the Cauchy problem
(\ref{causm})(\ref{cmhd}), of domain $[0,T)$; \parn
ii) $\bu_a =(u_a,b_a)$ is an approximate solution of the same Cauchy problem, of domain $[0,T_a)$.
\parn
We also introduce the following notations: \parn
iii) for any real $p$, $\epsilon_p$ and $\delta_p$ are estimators of order $p$
for the differential and datum error; $\cald_p$ is a growth
estimator of the same order (see Definition \ref{deferror}); \parn
iv) As in Eq. \rref{demu}, we put
$$ \mu := \min ( \nu , \eta)~. $$
v) $d^{+}/ d t$ stands for the \textsl{right, upper Dini derivative}; so,
for each function $f: [0, \tau) \to \reali$ (with $0 < \tau \leqs + \infty$) we have
\beq {d^{+} f \over d t} : [0,\tau) \to [-\infty,+\infty],  \quad t \mapsto {d^{+} f \over d t}(t)
:= \limsup_{h \to 0^{+}} {f(t+h) - f(t) \over h}\,.\label{dpiuf} \feq
\vskip 0.1cm \noindent
\begin{prop}  \label{lemmadissmo}
\textbf{Lemma.} For any real $p$ consider the function
$t \in [0, \min(T,T_a))$ $\mapsto$ \hbox{$\| \bu(t) - \bu_a(t) \|_{p}$}
(which is continuous, possibly non derivable at times
$t$ such that $\bu(t) = \bu_a(t)$). If $n, p \in \reali$ are such that $d/2 +1 < n \leqs p < +\infty$, this function fulfills
the inequality
\begin{align}
&{d^+ \over dt} \, \| \bu- \bu_a\|_{p} \leqs \label{diseqwsmo} \\
&\leqs -\mu \|\bu - \bu_a \|_p + (\Gp_p \cald_p
+ \Kp_p \cald_{p+1}))\| \bu - \bu_a \|_p
+ \Gp_{pn}\|\bu -\bu_a \|_n \| \bu -\bu_a \|_p +\epsilon_p \nonumber
\end{align}
$$ \mbox{everywhere on \, $[0,\min(T,T_a))$}~. $$
\end{prop}
\textbf{Proof.}
We work systematically on the time interval
$[0,\min(T,T_a))$, using the abbreviations
\beq
\bw := \bu - \bu_a, \qquad \bee \equiv \bee(\bu_a)~.
\eeq
The definition \rref{deferrsm} of the differential error
amounts to
\beq\nonumber
\frac{d \bu_a}{dt}= \bA \bu_a + \bP (\bu_a , \bu_a )  +\bee~;
\eeq
making use of (\ref{cmhd}) we obtain
\beq\nonumber
\begin{split}
& \frac{d \bw}{dt}= \frac{d \bu}{dt} - \frac{d \bu_a}{dt}=
\bA \bu + \bP (\bu, \bu)  - \bA \bu_a - \bP (\bu_a , \bu_a )  -\bee\\
& = \bA \bu_a + \bA \bw + \bP (\bu_a + \bw, \bu_a + \bw)
- \bA \bu_a - \bP (\bu_a , \bu_a ) -\bee~, \\
\end{split}
\eeq
i.e.,
\beq
\frac{d \bw}{dt} = \bA \bw + \bP (\bu_a , \bw) + \bP (\bw, \bu_a) + \bP (\bw, \bw)  -\bee~. \label{eqw}
\eeq
Let us consider an instant $t_0$ such that $\bw(t_0) \neq 0$.
In a neighborhood $I$ of this instant, the function $\| \bw \|_p$ is derivable and
\parn
\vbox{
\beq
\frac{d^+ \| \bw \|_p}{dt}=\frac{d \| \bw \|_p}{dt}= \frac{1}{2 \| \bw \|_p} \frac{d \| \bw \|^2_p}{dt}
= \frac{1}{\| \bw \|_p} \la \frac{d \bw}{d t} \vert \bw \ra_{p}
\label{lastlismo} \eeq
$$ =\frac{1}{\| \bw \|_p}\left(\la  \bA \bw \vert \bw \ra _{{p}} + \la  \bP (\bu_a, \bw) \vert \bw\ra _{{p}} + \la  \bP (\bw, \bu_a) \vert \bw \ra _{{p}}
+ \la  \bP (\bw , \bw) \vert \bw \ra _{{p}} \right. \\
\left. - \la  \bee \vert \bw \ra _{{p}}\right). $$
}
In the sequel we estimate the summands in the right hand side  of Eq. \rref{lastlismo}.
To this purpose we note that
the inequalities \rref{katop2np} \rref{basic2} \rref{katop2} for $\bP$,
the Schwarz inequality,
the inequalities  \rref{estiepsm} \rref{dednsm} defining the estimators $\epsilon_p$, $\cald_p$
and the inequality in \rref{ba2} for $\bA$  give
\begin{align}
& \la  \bA \bw  \vert \bw \ra _{{p}} \leqs - \mu \| \bw \|^2_p~, \label{nocalasmo} \\
& \la \bP (\bu_a , \bw) \vert \bw \ra _{{p}} \, \,  \leqs   \Gp_p \normP{\bu_a} \| \bw \|_p^2
\leqs \Gp_p \cald_p \| \bw \|_p^2~, \label{diseqpwwsmo} \\
& \la \bP(\bw , \bu_a) \vert \bw \ra _{{p}} \, \, \leqs \| \bP(\bw , \bu_a) \|_{p} \| \bw \|_{{p}} \leqs
  \Kp_p \normPp{\bu_a} \| \bw \|_p^2
\leqs  \Kp_p \cald_{p+1} \| \bw \|_p^2~,  \label{secondsmo}\\
& \la \bP (\bw , \bw ) \vert \bw \ra_p \leqs \Gp_{pn}\normn{\bw} \normP{\bw}^2~,
\label{estwwwsmo} \\
&-\la \bee \vert \bw \ra _{{p}}\, \,
\leqs \normP{\bee} \normP{\bw} \leqs \epsilon_p \| \bw \|_p~. \label{finalsmo}
\end{align}
Inserting Eqs. (\ref{nocalasmo}-\ref{finalsmo}) into Eq. (\ref{lastlismo}) one obtains
 \beq  \label{diseqnp}
 \frac{d^+ \| \bw \|_p}{dt} \leqs -\mu \normP{\bw} + (\Gp_p \cald_p + \Kp_p \cald_{p+1})\normP{\bw}
 + \Gp_{pn} \normn{\bw} \normP{\bw} + \epsilon_p~;
 \eeq
we repeat that this holds in a neighborhood of any instant $t_0$ such that $\bw(t_0) \neq 0$.
\parn
Now, let us consider
a instant $t_0$ such that $\bw(t_0) = 0$. In this case
we use a general result on the Dini derivative
(see e.g. \cite{Petr}), ensuring that
\beq {d^{+} \| \bw \|_p \over d t}(t_0) \leqs \| {d \bw \over d t}(t_0) \|_p~;
\eeq
on the other hand, Eq. \rref{eqw} for ${d \bw /d t}$ and the assumption
$\bw(t_0)=0$ give $(d \bw /d t)(t_0) = - \bee(t_0)$ so that (recalling again
Eq. \rref{estiepsm} for $\ep_p$)
\beq {d^{+} \| \bw \|_p \over d t}(t_0) \leqs \| \bee(t_0) \|_p \leqs \ep_p(t_0)~.
\eeq
But $\ep_p(t_0)$ equals the right hand side of Eq. \rref{diseqnp} at
$t = t_0$, again by the assumption $\bw(t_0)=0$. \parn
In conclusion,
Eq. \rref{diseqnp} is proved at each instant in $[0,\min(T,T_a))$;
recalling that $\bw = \bu - \bu_a$, we see that Eq. \rref{diseqnp} coincides with
the thesis \rref{diseqwsmo}. \hfill $\square$
\begin{prop}
\label{mainsmo}
\textbf{Proposition.}
Consider a real $n > d/2 + 1$, and
assume there is a function $\RR_n \in C([0,T_c), \reali)$,
with $T_c \in (0,T_a]$, fulfilling the following
\textsl{control inequalities}:
\beq {d^{+} \RR_n \over d t} \geqs - \mu \RR_n
+ (\Gp_n \cald_n + \Kp_n \cald_{n+1}) \RR_n + \Gp_n \RR^2_n + \ep_n
~\mbox{everywhere on $[0,T_c)$}, \label{cont1smo} \feq
\beq \RR_n(0) \geqs \delta_n \label{cont2smo} \feq
(with $d^{+}/d t$ as in Eq. \rref{dpiuf}; note that \rref{cont1smo} \rref{cont2smo} are fulfilled as equalities
by a unique function in $C^1([0,\Tc),\reali)$ for a suitable, maximal $\Tc$).
Then, (i) and (ii) hold.\\
(i)  The maximal solution $\bu$ of the MHD Cauchy problem \rref{causm} \rref{cmhd} and
its time of existence $T$ are such that
\beq T \geqs T_c~, \label{ttasmo} \feq
\beq \| \bu(t) - \bu_a(t) \|_n \leqs \RR_n(t) \qquad \mbox{for $t \in [0,T_c)$} \label{furthsmo} \feq
(and Eq. \rref{furthsmo} of course implies $\RR_n(t) \geqs 0$).
In particular, if $\RR_n$ is global ($T_c = +\infty)$ then $\bu$ is global as well
($T=+\infty$). \\
(ii) Consider any real $p > n$, and let $\RR_p \in C([0,T_c), \reali)$ be a solution of the linear control inequalities
\beq {d^{+} \RR_p \over d t} \geqs - \mu \RR_p
+ (\Gp_p \cald_p + \Kp_p \cald_{p+1}+ \Gp_{pn}\RR_n) \RR_p  + \ep_p
~\mbox{everywhere on $[0,T_c)$}, \label{cont1npsmo} \feq
\beq \RR_p(0) \geqs \delta_p\, .\label{cont2npsmo} \feq
Then
\beq
 \| \bu(t) - \bu_a(t) \|_p \leqs \RR_p(t) \qquad \mbox{for $t \in [0,T_c)$}
\eeq
(which of course implies $\RR_p(t) \geqs 0$).
 Conditions (\ref{cont1npsmo}) (\ref{cont2npsmo}) are fulfilled as equalities by a unique function $\RR_p \in C^1 ([0,T_c) , \reali)$, given explicitly by
 \beq  \label{soleq}
 \RR_p(t)=\mathrm{e}^{-\mu t + \AA_p(t)}\left( \delta_p +\int_0^t ds \, \mathrm{e}^{\mu s -\AA_p (s)}
 \epsilon_p (s)\right) \qquad \mbox{for $t \in [0 , T_c)$},
 \eeq
 \beq  \label{defap}
 \AA_p(t) :=\int_0^t ds (\Gp_p \cald_p (s) + \Kp_p \cald_{p+1} (s) + \Gp_{pn} \RR_n(s)).
 \eeq
\parn
\end{prop}
\textbf{Proof.}
(i) The inequality \rref{diseqwsmo} of the
previous lemma and Eq. \rref{estdatumsm} for the estimator
$\delta_p$, with $p=n$, read
\beq\nonumber
{d^+ \over dt} \, \| \bu- \bu_a\|_{n} \leqs \eeq
$$ \leqs - \mu \| \bu - \bu_a \|_{n} +
(\Gp_n \cald_n + \Kp_n \cald_{n+1}) \| \bu - \bu_a \|_{n} +
\Gp_n \| \bu - \bu_a \|_{n}^2 + \epsilon_n $$
$$ \text{everywhere on} \, [0,\min(T,T_a))~; $$
$$ \lVert \bu_a (0) - \mathbf{u}_0 \rVert_{{n}} \leqs
\delta_n. $$
These inequalities for $\| \bu - \bu_a \|_n$
are like the control inequalities
\rref{cont1smo} \rref{cont2smo} for $\RR_n$, with the
reverse order relation; so,
the comparison theorem of \v{C}aplygin-Lakshmikhantam \cite{Las} \cite{Mitr} ensures that
\beq \| \bu(t) - \bu_a(t) \|_n \leqs \RR_n(t) \qquad \mbox{for $t \in [0, \min(T, T_a, T_c)) =
[0, \min(T, T_c))$}~. \feq
Finally, let us prove that
\beq T \geqs T_c \quad \mbox{(i.e., $\min(T, T_c) = T_c$)}~; \feq
indeed, if it were $T <  T_c$, for
all $t \in [0,T)$ we would have
$\| \bu(t) \|_n$ $\leqs \| \bu(t) - \bu_a(t) \|_n$ $+ \| \bu(t) \|_n$ $\leqs \RR_n(t) + \cald_n(t)$
and this would imply $\limsup_{t \to T^{-}} \|\bu(t) \|_n
\leqs \RR_n(T) + \cald_n(T) < +\infty$, contradicting the blow-up criterion (ii) of Proposition \ref{procau}. \parn
(ii) From (i) we know that $T \geqs \Tc$ and $\| \bu-\bu_a \|_n \leqs \RR_n$ on $[0, T_c)$.
Making use of this result and of the inequality (\ref{diseqwsmo}), which is valid on the interval $[0, \min (T, T_a))$,
we obtain that, on the shorter interval $[0, T_c)$, there holds
\beq
\begin{split}
&{d^+ \over dt} \, \| \bu- \bu_a\|_{p} \leqs\\
&\leqs -\mu \|\bu - \bu_a \|_p + (\Gp_p \cald_p +\Kp_p \cald_{p+1}))\| \bu - \bu_a \|_p
+ \Gp_{pn}\RR_n \| \bu -\bu_a \|_p +\epsilon_p~.  \label{diseqwnp}
\end{split}
\eeq
The inequality (\ref{diseqwnp}) and the relation $\| \mathbf{u}_a (0)- \mathbf{u}_0 \|_p \leqs \delta_p$
have the same structure as the relations (\ref{cont1npsmo}) and (\ref{cont2npsmo}), with the reversed order.
Therefore, as in the proof of (i) one can apply  a comparison argument \textsl{\`a la} \v{C}aplygin-Lakshmikhantam
ensuring that
\beq
 \| \bu(t) - \bu_a(t) \|_p \leqs \RR_p(t) \qquad \mbox{for $t \in [0,T_c)$}.
\eeq
Finally, using elementary facts on linear ODEs, one checks
that the function $\RR_p$ defined by (\ref{soleq}) and (\ref{defap})
is the unique $C^1$ function on $[0, T_c)$ satisfying  (\ref{cont1npsmo}) and (\ref{cont2npsmo}) as
equalities.  \fine
\begin{rema}
\label{remcont}
\textbf{Remark.} In the sequel we often refer to the case
mentioned just after Eqs.  \rref{cont1smo} \rref{cont2smo},
in which these control inequalities
are fulfilled as equalities by a $C^1$ function $\RR_n \in C^1([0,T_c),\reali)$, to
be determined; this gives rise
to the \textsl{control Cauchy problem}
\ben
\frac{d \calr_n}{dt} & = &- \mu \calr_n + (\Gp_n \cald_n + \Kp_n   \cald_{n+1})\calr_n + \Gp_n \calr_n^2 + \epsilon_n, \label{contreq}\\
\calr_n (0) & = &  \delta_n  \label{contreq0}
\een
for the unknown $\RR_n$.
\end{rema}

\section{Simple analytical estimates arising from
Proposition \ref{mainsmo}}
\label{secappan}
Let us consider again the Cauchy problem \rref{causm} \rref{cmhd};
throughout this section $\bu = (u,b) \in C^\infty([0,T),\HD{\infty})$ is
its maximal solution. In the sequel we present some elementary,
but useful consequences of Proposition \ref{mainsmo} based
on a very simple choice of the approximate solution $\bu_a$ mentioned
therein: the latter is assumed to be the zero function. \parn
\begin{prop}  \label{tczero}
\textbf{Lemma.}
Let us introduce the function
\beq
\bu_a : [0 , + \infty) \rightarrow \HD{\infty} , \qquad \bu_a(t):=0 \quad \mbox{for all $t$}
\eeq
and regard it as an approximate solution of problem \rref{causm}\rref{cmhd}.
The differential and datum errors of this approximate solution are
\beq
 \bee (\bu_a)(t)=\mathbf{0}, \quad \mbox{for all $t \in [0, + \infty)$}, \qquad \bu_a (0) - \mathbf{u}_0=-\mathbf{u}_0.
\eeq
Consequently, the zero approximate solution has the following
differential error, datum error and growth estimators of any order $p$:
\ben
\epsilon_p :=0~ , \qquad \delta_p=\| \mathbf{u}_0 \|_p~, \\
\cald_p(t) := 0~ .
\een
For any fixed $n > d/2 + 1$, the following holds: \parn
(i) The control Cauchy problem (\ref{contreq})-(\ref{contreq0}) with these estimators takes the form
\ben
\frac{d \calr_n}{dt}=- \mu \calr_n +  \Gp_n \calr_n^2 ,  \label{caucr}\\
\calr_n (0) =  \normn{\mathbf{u}_0}  \label{caucr0}
\een
and admits a solution $\calr_n \in C^1 ([0, T_c)),[0, + \infty))$, given by
\beq \RR_n(t) := {\| \mathbf{u}_0 \|_n e^{-\mu t} \over 1 - \Gp_n \| \mathbf{u}_0 \|_n e_{\mu}(t)} \qquad \mbox{for
$t \in [0,\Tc)$}~. \label{ern} \feq
Here
\beq \Tc := \left\{ \barray{ll} + \infty & \mbox{if $\mu >0$, $\| \mathbf{u}_0 \|_n \leqs {\mu/\Gp_n}$}~, \\
- \dd{1 \over \mu} \log\left(1 - \dd{\mu \over \Gp_n \| \mathbf{u}_0 \|_n} \right) & \mbox{if $\mu > 0$,
$\| \mathbf{u}_0 \|_n > {\mu/\Gp_n}$,} \\
\dd{1 \over \Gp_n \| \mathbf{u}_0 \|_n} & \mbox{if $\mu=0$} \farray \right. \label{ta} \feq
(intending $1/(\Gp_n \| \mathbf{u}_0 \|_n) := + \infty$ if $\mathbf{u}_0 = \mathbf{0}$), and
\beq e_{\mu}(t) := \left\{ \barray{ll} \dd{1 - e^{-\mu t} \over \mu} & \mbox{if $\mu > 0$}, \\
t & \mbox{if $\mu = 0$} \farray \right. \label{enu} \feq
(note that $t = \lim_{\mu \vain 0^{+}} \dd{1 - e^{-\mu t} \over \mu}$). \parn
(ii) For each real $p > n$, with $\RR_n$ defined by \rref{ern} \rref{enu} and
with the above mentioned estimators, the function $\RR_p$ of Eq.
\rref{soleq} is as follows:
\beq \RR_p(t) =
{\| \mathbf{u}_0 \|_p e^{-\mu t} \over \Big[1 - \Gp_n \| \mathbf{u}_0 \|_n e_{\mu}(t)\Big]^{\Gp_{p n}/\Gp_n}}
\qquad \mbox{for
$t \in [0,\Tc)$}~.
\feq
\end{prop}
\textbf{Proof.} It is obtained by elementary computations
(similar to those presented in \cite{Mor15}, page 305 for
the zero approximate solution of the NS Cauchy problem). \fine
The previous lemma allows to infer the following statement, similar to
Proposition 5.1 of \cite{Mor15} on the NS Cauchy problem.
\begin{prop}  \label{propzero}
\textbf{Proposition.}
Let $\bu = (u,b) \in C^\infty([0,T),\HD{\infty})$ be the maximal solution of the Cauchy
problem \rref{causm}\rref{cmhd}. Fix any real $n > d/2 + 1$ and define $T_c$ and $e_\mu$ as in Eqs. (\ref{ta}) and (\ref{enu}). Then
\beq
 T \geqs \Tc~, \qquad \| \bu(t) \|_n \leqs
{ \| \mathbf{u}_0 \|_n e^{-\mu t} \over 1 - \Gp_n \| \mathbf{u}_0 \|_n e_{\mu}(t)} \qquad \mbox{for $t \in [0,\Tc)$}~,
\label{cobg} \eeq
\beq \| \bu(t) \|_p \leqs
{\| \mathbf{u}_0 \|_p e^{-\mu t} \over \Big[1 - \Gp_n \| \mathbf{u}_0 \|_n e_{\mu}(t)\Big]^{\Gp_{p n}/\Gp_n}}
\qquad \mbox{for real $p > n$ and $t \in [0,\Tc)$}~.\eeq
In particular
\beq \mbox{$T=\Tc=+\infty$ \quad if ~$\| \mathbf{u}_0 \|_n \leqs \dd{\mu \over \Gp_n}$}~; \label{coglob}
\feq
in this case $\bu$ is global.
\end{prop}
\textbf{Proof.}
Use Proposition \ref{mainsmo} with the
approximate solution $\bu_a(t) :=0$, along with
the previous Lemma \ref{tczero}; in particular,
statement \rref{coglob} follows from Eq. \rref{ta}
of this Lemma. \fine
Hereafter we present two consequences of Proposition \ref{propzero};
these have close analogies with Corollaries 5.3 and 5.4
of \cite{appeul} on NS equations.
\begin{prop}
\label{coro}
\textbf{Corollary.} Consider again the maximal solution
$\bu$ of the Cauchy
problem \rref{causm}\rref{cmhd}. Assume that
\beq \| \bu(t_1) \|_n \leqs {\mu \over \Gp_n}
\qquad \mbox{for some $n > d/2 + 1$ and $t_1 \in [0,T)$}~. \label{assume} \feq
Then:
\beq T = + \infty, \qquad \| \bu(t) \|_n \leqs
{\| \bu(t_1) \|_n e^{-\mu (t-t_1)} \over 1 - \Gp_n \| \bu(t_1) \|_n \, e_{\mu}(t-t_1)}
 \quad \mbox{for $t \in [t_1, +\infty)$}~,
\label{tes0} \feq
with $e_{\mu}$ as in Eq. \rref{enu}.
\end{prop}
\textbf{Proof.}
The function $\bu \restriction [t_1, T)$
is the maximal solution of the Cauchy problem with initial datum $\bu(t_1)$ \textsl{specified
at time $t_1$}, rather than at time $0$; therefore, after a shift in the time
variable we can apply to this function  Eqs. \rref{cobg}\rref{coglob}, which yield the
thesis \rref{tes0}. \fine
\begin{prop}
\label{coro2}
\textbf{Corollary.}
Let $\bu_a = (u_a , b_a) \in  C^1([0, T_a), \HD{\infty} )$ be any approximate
solution of the Cauchy problem \rref{causm}\rref{cmhd}
with estimators $\ep_n, \delta_n, \cald_n, \cald_{n+1}$
for some $n > d/2 + 1$;
assume the control inequalities
\rref{cont1smo} \rref{cont2smo} to possess a solution $\Rr_n \in C([0,\Tc),\reali)$,
with $\Tc \in (0,T_a]$ (this is nonnegative, see Proposition
\ref{mainsmo}). Finally, assume
\beq (\cald_n + \Rr_n)(t_1) \leqs {\mu \over \Gp_n}
\qquad \mbox{for some $t_1 \in [0,\Tc)$}~. \label{iprn} \feq
Then, the maximal exact solution $\bu$ of the Cauchy
problem \rref{cmhd} has the following features:
\beq T = + \infty, \quad \| \bu(t) \|_n \leqs
{(\cald_n + \Rr_n)(t_1) e^{-\mu (t-t_1)} \over 1 - \Gp_n (\cald_n + \Rr_n)(t_1) \, e_{\mu}(t-t_1)}
~\mbox{for $t \in [t_1, +\infty)$}~.
\label{tes1} \feq
\end{prop}
\textbf{Proof.}
Writing $\| \bu(t_1) \|_n \leqs \| \bu_a(t_1) \|_n +
\| \bu(t_1) - \bu_a(t_1) \|_n$ and using
at time $t_1$ the bounds \rref{dednsm} (with $p=n$)
and \rref{furthsmo} we get
\beq \| \bu(t_1) \|_n \leqs (\cald_n + \Rr_n)(t_1)~. \label{furth1} \feq
Now the assumption \rref{iprn} gives the
inequality
$$ \| \bu(t_1) \|_n \leqs {\mu \over \Gp_n}~, $$
which has the form \rref{assume}. By Corollary \ref{coro} we have Eq. \rref{tes0}, and
inserting therein Eq. \rref{furth1} we obtain
the thesis \rref{tes1}. \fine
\salto
\noindent
\textbf{Applications to specific initial conditions.}
In this subsection the space dimension is
\beq d = 3, \feq
and we apply Proposition \ref{propzero}  with $n=3$.
Eq. \rref{coglob} from the cited proposition, together with Eq. \rref{k3g3} about
the constant $\Gp_3$, ensures the following: the MHD Cauchy problem
with a datum $\bu_0 \in \HD{\infty}$ has a global solution if
\beq \mu \geqs \Gp_3 \| \bu_0 \|_3~, \qquad \Gp_3 = 0.620. \label{coglob3} \feq
(in the above $\mu := \min(\nu,\eta)$, as in \rref{demu}). Hereafter we write explicitly
the condition \rref{coglob3} for two initial data often adopted in theoretical
studies on MHD turbulence: a three-dimensional Orszag-Tang vortex  and an Arnold-Beltrami-Childress (ABC) flow
with a perturbing magnetic field (see, e.g., \cite{Min06} \cite{Car09}).
\vskip 0.1cm \noindent
i) {\it Orszag-Tang vortex}. This is the datum $\bu_0 = (u_0, b_0)$ where, as in
\cite{Min06},
\parn{\vbox{
\beq u_0(x) := (-2 \sin x_2, 2 \sin x_1, 0)~, \feq
$$ b_0(x) := \beta (\, - 2 \sin(2 x_2) + \sin x_3, 2 \sin x_1 + \sin x_3, \sin x_1 + \sin x_2\,)
\quad (\beta \in \reali). $$
}}
We have the Fourier representations
\ben  \label{orst}
u_0=\sum_{k=\pm a_1 , \pm a_2 } u_{0k} e_k, \qquad b_0=\sum_{k=\pm a_1 , \pm a_2 , \pm a_3 , \pm a_4} b_{0k} e_k,
\een
\beq
a_1:=(1,0,0), \qquad a_2:=(0,1,0),   \label{a1a2}
\eeq
$$ a_3:=(0,0,1), \qquad a_4:=(0,2,0), $$
\begin{align}
& u_{0, \pm a_1} &:=& \mp (2\pi)^{3/2} i (0,1,0), \qquad & u_{0, \pm a_2} &:=& \pm (2\pi)^{3/2} & i (1,0,0),  \nonumber \\
& b_{0, \pm a_1} &:=& \mp (2 \pi)^{3/2} i \beta \left(0,1,\frac{1}{2}\right), \qquad
& b_{0, \pm a_2} &:= & \mp (2 \pi)^{3/2} & i \beta \left( 0,0, \frac{1}{2}\right),  \label{ot}\\
& b_{0, \pm a_3} &:=& \mp (2 \pi)^{3/2} i \beta \left(\frac{1}{2} , \frac{1}{2} ,0\right), \qquad
& b_{0, \pm a_4} &:=& \pm (2 \pi )^{3/2} & i \beta (1,0,0)  \nonumber
\end{align}
({\footnote{To avoid misunderstandings, let us explain the notations $\pm, \mp$
in Eq. \rref{ot} and in the subsequent Eq. \rref{abcdet}. As an example, the first
line in Eq. \rref{ot} means that
$u_{0, a_1} := - (2\pi)^{3/2} i (0,1,0)$, $u_{0, - a_1} := (2\pi)^{3/2} i (0,1,0)$ and
$u_{0, a_2}:= (2\pi)^{3/2} i (1,0,0)$, $u_{0,- a_2}:=- (2\pi)^{3/2} i (1,0,0)$.
}}).
We find
\beq
\lVert \bu_0 \rVert_3= (2 \pi)^{3/2} \sqrt{4  + 132 \beta^2}~;
\label{eqnorma} \eeq
from here, we infer that the condition \rref{coglob3} of global existence
for the MHD Cauchy problem with the Orszag-Tang datum (\ref{orst}) holds if
\beq
\mu \geqs 9.77 \,\sqrt{4 + 132 \beta^2}.
\eeq
\\
ii) {\it ABC flow with perturbing magnetic field}. This is the datum
$\bu_0 = (u_0,b_0)$ where, as in \cite{Car09},
\beq u_0(x) := ( B \cos x_2 + C \sin x_3, A \sin x_1 + C \cos x_3, A \cos x_1 + B \sin x_2)~, \label{abcxyz} \feq
$$ b_0(x) := D (\sin x_1 \cos x_2, - \cos x_1 \sin x_2, 0) \quad (A,B,C,D \in \reali)~. $$
We have the Fourier representations
\ben  \label{abc}
u_0=\sum_{k=\pm a_1 , \pm a_2 , \pm a_3} u_{0k} e_k, \qquad b_0=\sum_{k=\pm a_5 , \pm a_6} b_{0k} e_k,
\een
\beq a_1,a_2, a_3~\mbox{as in \rref{a1a2}}, \quad a_5 := (1,1,0), \quad a_6  :=  (1, -1,0), \feq
\begin{align}
& u_{0,  \pm a_1} &:= & \phantom{xx} (2 \pi)^{3/2} \frac{A}{2} (0,  \mp i , 1), \qquad
& u_{0,  \pm a_2} &:=  \phantom{x} (2 \pi)^{3/2} \frac{B}{2} (1, 0 , \mp i),  \nonumber \\
& u_{0,  \pm a_3} &:= & \phantom{xx} (2 \pi)^{3/2} \frac{C}{2} (\mp i, 1 , 0),
\qquad &  &  &  & ~ \label{abcdet}  \\
& b_{0,  \pm a_5} &:=& \pm (2 \pi)^{3/2} i \frac{D}{4} (-1, 1 , 0),
& b_{0,  \pm a_6} &:= \pm (2 \pi)^{3/2} i \frac{D}{4} (-1, -1 , 0), \nonumber
\end{align}
In this case
\beq
\lVert \bu_0 \rVert_3= (2 \pi)^{3/2} \, \sqrt{A^2 + B^2 +C^2 + 4 D^2}~;
\label{eqnorma2} \eeq
from here, we see that the condition \rref{coglob3} of global existence
for the MHD Cauchy problem with the datum (\ref{abc}) holds if
\beq
\mu \geqs 9.77 \sqrt{A^2 + B^2 +C^2 + 4 D^2}.
\label{muglob} \eeq
\section{The Galerkin approximate solutions for
the MHD equations, and their errors}
\label{galegen}
As well known, a Galerkin approximate solution
for the NS equations, the MHD equations or many other PDEs
is supported by finite sets of Fourier modes.
Hereafter we adapt to the MHD case the presentation
of the Galerkin approach already given
in \cite{appeul} for the NS case (on this
construction, see again the notice
at the end of the Introduction); in particular,
Definition \ref{defgal} and Propositions \ref{galeelledue}, \ref{galecompon}, \ref{lem67}
in the present section correspond, respectively, to Definition 6.3,
Lemma 6.4, Proposition 6.7 and Lemma 6.8 in \cite{appeul}.
\parn
Throughout the section we consider a set $G$ such that
\beq \G \subset \Zd_0~, \qquad \G~\mbox{finite}~, \qquad
k \in \G \Leftrightarrow - k \in \G~. \label{modesg} \feq
\salto
\textbf{Galerkin subspaces and projections.}
By definition, the Galerkin subspace and the projection corresponding to $\G$ are, respectively:
\beq \HM{G} := \{ v \in \Dsz~|~v_k = 0~\mbox{for $k \in \Zd_0 \setminus G$} \} \label{hg} \feq
$$ = \{ \sum_{k \in \G} v_k e_k~|~v_k \in \complessi^d, \overline{v_k} = v_{-k}, k
\sc \, v_k = 0~\mbox{for all $k$} \}~; $$
\beq \EG : \Dsz~ \vain \HM{G}~, \qquad v =
\sum_{k \in \Zd_0} v_k e_k \mapsto \EG v := \sum_{k \in \G} v_k e_k~. \label{pg} \feq
It is clear that
\beq \HM{G} \subset \HM{\infty};~\EG(\HM{p}) = \HM{G}~\mbox{for $p \in \reali \cup \{\infty\}$};~
\Delta(\HM{G}) = \HM{G}~. \feq
Moreover
\beq \la \EG v | w \ra_{p} = \la v | \EG w \ra_{p} \quad \mbox{for $p \in \reali$, $v,w \in \HM{p}$}~. \feq
Let us also mention that
\beq \| (1 - \EG) v \|_p \leqs {\| v \|_{q} \over
|\G|^{q-p}}~\mbox{for $p,q \in \reali$, $p \leqs q$, $v \in \HM{q}$},
\qquad |\G| := \min_{k \in \Zd_0 \setminus \G} |k| \label{egmp} \feq
(see e.g. \cite{appeul}, Lemma 6.2).
We can introduce a ``two-component'' Galerkin subspace and projection associated to $G$
which are, respectively,
\beq \HD{G} := \HM{G} \times \HM{G}~, \label{debeg} \feq
\beq \bEG : \bDd'_{\Sigma 0}~ \vain \HD{G}~, \qquad \bv = (v,b) \mapsto \bEG \bv := (\EG v,\EG b)~. \feq
The previous statements about $\HM{G}$ and $\EG$ have obvious implications for their two-component
analogues. In particular:
\beq \bA (\HD{G}) \subset \HD{G}~, \feq
\beq \la \bEG \bv | \bw \ra_{p} = \la \bv | \bEG \bw \ra_{p} \quad
\mbox{for $p \in \reali$, $\bv,\bw \in \HD{p}$},
\label{egsimm} \feq
\beq \| (1 - \bEG) \bv \|_p \leqs {\| \bv \|_{q} \over
|\G|^{q-p}}~~\mbox{for $p, q \in \reali$, $p \leqs q$, $\bv \in \HD{q}$}, \label{begmp} \feq
with $|G|$ as in Eq. \rref{egmp}.
\salto
\textbf{Galerkin approximate solutions.}
Let us be given
$\nu,\eta \in [0,+\infty)$ and $\bu_0 = (u_0,b_0) \in \HD{\infty}$.
\begin{prop}
\label{defgal}
\textbf{Definition.} The \textsl{Galerkin
approximate solution} of the MHD equations corresponding to
$\nu,\eta,\bu_0$ and to the set of modes $G$ is the maximal
(i.e., unextendable) solution
$\ug$
of the following Cauchy problem:
\beq \mbox{Find $\ug = (u_G,b_G) \in C^\infty([0,\Tg), \HD{G})$ such that} \label{galer} \feq
$$ {d \ug \over dt} = \bA \ug + \bEG \bP(\ug, \ug) ~, \qquad \ug(0) = \bEG \bu_0~. $$
\ffine
\end{prop}
Let us note that the Cauchy problem \rref{galer} rests on
the finite dimensional vector space $\HD{G}$ and on the $C^\infty$ function
$\HD{G} \vain \HD{G}$, $\bv \mapsto \bA \bv + \EG \bP(\bv, \bv)$.
Recalling the definitions \rref{deba} \rref{debp} \rref{debeg} of  $\bA$,$\bP$,$\bEG$ we can rephrase as follows
problem \rref{galer} in terms of the components $u_G,b_G$ of $\ug$: \parn
\vbox{
\beq \mbox{Find $u_G,b_G \in C^\infty([0,\Tg), \HM{G})$ such that} \label{galerco} \feq
$$ {d u_G \over dt} = \nu \Delta u_G + \EG \P(u_G, u_G) - \EG \P(b_G,b_G)~, \quad
{d b_G \over dt} = \eta \Delta b_G + \EG \P(u_G, b_G) - \EG \P(b_G,u_G)~, $$
$$ u_G(0) = \EG u_0~, ~~b_G(0) = \EG b_0~. $$
}
The standard theory of ODEs in finite dimension grants local existence
and uniqueness for the (maximal) solution of \rref{galer} (or \rref{galerco}); according to the same
theory, the finiteness of $\Tg$ would imply $\lim \sup_{t \to \Tg^{-}}
\| \ug(t) \| = + \infty$ for any norm $\|~\|$ on $\HD{G}$ (recall
that all norms on a finite dimensional vector space are equivalent). \parn
Hereafter
we consider, in particular, the evolution of the $L^2$ norm
$t \mapsto \| \ug(t) \|_{L^2}$ (giving twice the ``energy'' of the Galerkin
solution), and point out its
implications for $\Tg$:
\begin{prop}
\textbf{Proposition.}
\label{galeelledue}
Let us consider the maximal solution $\ug$ of problem \rref{galer}, of domain
$[0,\Tg)$. With $\mu$ as in \rref{demu}, one has
\beq {d \over d t} \| \ug \|^2_{L^2} = 2 \la \bA \ug | \ug \ra_{L^2} \leqs - 2 \mu \| \ug \|^2_{L^2}~.
\label{gathelast} \feq
This implies the following, for $t \in [0,\Tg)$:
\beq \| \ug(t) \|_{L^2} \left\{ \barray{ll} = \| \bEG \bu_0 \|_{L^2} & \mbox{if $\nu=\eta=0$,} \\
\leqs \| \bEG \bu_0 \|_{L^2} e^{- \mu t} & \mbox{for all $\nu,\eta \geqs 0$}. \farray \right.
\label{gaimp} \feq
A consequence of these estimates is that
\beq \Tg = + \infty~. \feq
\end{prop}
\textbf{Proof.}
Writing $\| \ug \|^2_{L^2} = \la \ug | \ug \ra_{L^2}$, taking the $t$ derivative and
using Eq. \rref{galer} we get
\beq {d \over d t} \| \ug \|^2_{L^2} = 2 \la {d \ug \over d t} | \ug \ra_{L^2}
=  2 \la \bA \ug | \ug \ra_{L^2} + 2 \la \bEG \bP(\ug,\ug) | \ug \ra_{L^2} ~. \label{uno} \feq
On the other hand, using Eq. \rref{egsimm} with $p=0$ and Eq. \rref{brecall0},
\beq \la \bEG \bP(\ug,\ug) | \ug \ra_{L^2} = \la \bP(\ug,\ug) | \bEG \ug \ra_{L^2} =
\la \bP(\ug,\ug) | \ug \ra_{L^2} = 0~; \label{is1} \feq
moreover, due to Eq. \rref{ba2} with $p=0$,
\beq \la \bA \ug | \ug \ra_{L^2} \leqs -\mu \| \ug \|^2_{L^2}~. \label{is2} \feq
Inserting Eqs. \rref{is1} \rref{is2} into Eq. \rref{uno} we obtain Eq. \rref{gathelast}.
Eq. \rref{gaimp} is a straightforward consequence of Eq. \rref{gathelast}
(in connection with this statement, let us recall that $\nu=\eta=0$ implies $\bA=0$). \parn
Finally, if $\Tg$ were finite we would have $\lim \sup_{t \to \Tg^{-}}
\| \ug(t) \|_{L^2} = + \infty$ (see the remark a few lines before the present proposition);
this would contradict Eq. \rref{gaimp}, so $\Tg = + \infty$. \fine
\salto
\textbf{Fourier representation of the Galerkin approximants.}
Let us fix $\nu,\eta \geqs 0$ and an initial datum $\bu_0 = (u_0,b_0) \in \HD{\infty}$;
we consider the Fourier expansions
\beq u_0 = \sum_{k \in \Zd_0} u_{0 k} e_k~, \qquad b_0 = \sum_{k \in \Zd_0} b_{0 k} e_k \label{fur2} \feq
with coefficients $u_{0 k},b_{0 k} \in \complessi^d$; these fulfill the conditions
\beq u_{0\,-k} = \overline{u}_{0 k}, \quad b_{0 \,-k} = \overline{b}_{0 k},
\quad k \sc u_{0 k} = k \sc b_{0 k} = 0~. \label{codato} \feq
Denoting again with $G$ a finite set of modes as in \rref{modesg}, we provisionally
write $\ug = (u_G,b_G)$ to indicate an unspecified function in $C^\infty([0,+\infty), \HD{G})$
and associate to it two families of Fourier coefficients
$\gamma_{k}, \beta_k \in C^\infty([0,+\infty),\complessi^d)$ ($k \in G$), defined by
\beq u_G(t) = \sum_{k \in G} \ga_{k}(t) e_k ~, \qquad b_G(t) = \sum_{k \in G} \be_{k}(t) e_k
\quad \mbox{for $t \in [0,+ \infty)$}~; \label{fur1} \feq
we note that
\beq \ga_{-k} = \overline{\ga}_k, \quad \be_{-k} = \overline{\be}_k, \quad k \sc \ga_k = k \sc \be_k =0~.
\label{cond} \feq
\begin{prop}
\label{galecompon}
\textbf{Proposition.}
$\ug$ fulfills the Cauchy
problem \rref{galer} (or \rref{galerco}) if and only if its coefficients $\ga_k, \be_k$ fulfill the
following for all $k \in G$:
\parn
\vbox{
\beq {d \ga_k \over d t} = - \nu |k|^2 \ga_k - {i \over (2 \pi)^{d/2}}
\sum_{h \in G} \Big(  [ \ga_h \sc (k-h) ] \LP_k \ga_{k-h} - [ \be_h \sc (k-h) ] \LP_k \be_{k-h}
\Big)  ~, \label{galecomp} \feq
$$ {d \be_k \over d t} =  - \eta |k|^2 \be_k - {i \over (2 \pi)^{d/2}}
\sum_{h \in G} \Big(  [ \ga_h \sc (k-h) ] \LP_k \be_{k-h} - [ \be_h \sc (k-h) ] \LP_k \ga_{k-h}
\Big) ~, $$
$$ \ga_k(0) = u_{0 k}~, \qquad \be_k(0) = b_{0 k} $$}
(intending $\ga_{k-h}, \be_{k-h} := 0$ if $k - h \not \in G$; as for $\LP_k$,
recall the explanations after Eq. \rref{defler}). \parn
\end{prop}
\textbf{Proof.} Clearly, $\ug$ fulfills problem \rref{galerco} if and only if, for all $k \in G$,
\parn
\vbox{
$$ {d \ga_k \over dt} = - \nu |k|^2 \ga_k + \P_k(u_G, u_G) - \P_k(b_G,b_G) ~, \quad
{d \be_k \over dt} = - \eta |k|^2 \be_k + \P_k(u_G, b_G) - \P_k(b_G,u_G) ~, $$
\beq \ga_k(0) = u_{0 k}~, ~~\be_k(0) = b_{0 k}~. \label{thus} \feq
}
Using the representation \rref{eq2} for the Fourier component $\P_k(~,~)$
with the fact that $(u_G)_h = \ga_h$, $(b_G)_h = \be_h$ for $h \in G$ and $(u_G)_h = (b_G)_h = 0$
for $h \in \Zd \setminus G$, we obtain
$$ \P_k(u_G, u_G) = - \frac{i}{(2 \pi)^{d/2}} \sum_{h \in G} [\ga_h
\sc (k-h)] \LL_k \ga_{k-h}~, $$
$$ \P_k(b_G, b_G) = - \frac{i}{(2 \pi)^{d/2}} \sum_{h \in G} [\be_h
\sc (k-h)] \LL_k \be_{k-h} $$
and so on, thus Eq. \rref{thus} coincides with Eq.\rref{galecomp}. \fine
\begin{rema}
\textbf{Remark.}
We can regard the system \rref{galecomp} as a Cauchy problem for finitely many unknown
functions $\ga_k, \be_k \in C^\infty([0,+\infty),\complessi^d)$ ($k\in G$).
An elementary argument based on Eqs. \rref{galecomp} and \rref{codato}
shows that the unique solution $(\ga_k,\be_k)_{k \in G}$
of this Cauchy problem automatically fulfills the conditions
\rref{cond} (a similar statement on the Galerkin approximants
for the NS equations is proved in \cite{appeul}, Proposition 6.7).
\end{rema}
\salto
\textbf{The Galerkin solutions in the framework of Section
\ref{secapp}.} From now on we consider,
for given $\nu,\eta \geqs 0$ and $\bu_0 = (u_0, b_0) \in \HD{\infty}$: \parn
i) the MHD Cauchy problem \rref{causm} \rref{cmhd} and its
maximal solution $\bu \in C^\infty([0,T), \HD{\infty})$; \parn
ii) the Galerkin approximant $\ug = (u_G,b_G) \in C^\infty([0,+\infty), \HD{G})$
defined by Eq. \rref{galer}, for a finite set $G$ of modes
as in Eq. \rref{modesg}. We also refer to the Fourier representations
\rref{fur2}\rref{fur1}\rref{galecomp} of $\bu_0, \ug$ and of the Galerkin
Cauchy problem.
\parn
We regard $\ug$ as an approximate solution of the MHD Cauchy problem
\rref{causm} \rref{cmhd}, to be treated using the methods of
Section \ref{secapp} (and \ref{secappan}); to this purpose,
we need growth and error estimators for $\ug$. \parn
Concerning the growth of $G$, we have the tautological growth estimators
\beq \cald_p(t) := \| \ug(t) \|_p = \sqrt{ \sum_{k \in G} |k|^{2 p} \big(|\ga_k(t)|^2 + |\beta_k(t)|^2\big)}
\quad (p \in \reali, t \in [0,+\infty))~; \feq
the errors of $\ug$ and their estimators are discussed heferafter.
\begin{prop}
\label{lem67}
\textbf{Proposition.}
(i) The Galerkin solution $\ug$ has the datum error
\beq \ug(0) - \bu_0 = -(1 - \bEG) \bu_0 = - \big( \sum_{k \in \Zd_0 \setminus G} u_{0 k} e_k,
\sum_{k \in \Zd_0 \setminus G} b_{0 k} e_k \big)~.
\label{ovv1} \feq
For each $p \in \reali$, the datum error has the tautological estimator
\beq \delta_p := \| \ug(0) - \bu_0 \|_{p} =
\sqrt{\sum_{k \in \Zd_0 \setminus G} |k|^{2 p} \big(|u_{0 k}|^2 + | b_{0 k}|^2\big)}~. \label{ovv2} \feq
and a rougher estimator, depending on another real number $q \geqs p$,
\beq \| \ug(0) - \bu_0 \|_{p} \leqs \delta'_{p q}~, \qquad \delta'_{p q} :=
{\| \bu_0 \|_{q} \over |G|^{q-p}}\,.
\label{ovv3} \feq
(ii) The differential error of $\ug$ is
\beq \bee(\ug) = - (1 - \bEG) \bP(\ug, \ug) =
- \left(\sum_{k \in dG} \rho_k e_k~, \sum_{k \in dG} \sigma_k e_k \right) \label{prov1} \feq
where:
\beq dG := (G + G) \setminus (G \cup \{0\})~, \label{630} \feq
$$ \rho_k := - {i \over (2 \pi)^{d/2}}\sum_{h \in G}
\Big(  [ \ga_h \sc (k-h) ] \LP_k \ga_{k-h} - [ \be_h \sc (k-h) ] \LP_k \be_{k-h}
\Big)~, $$
$$ \sigma_k := - {i \over (2 \pi)^{d/2}}\sum_{h \in G}
\Big(  [ \ga_h \sc (k-h) ] \LP_k \be_{k-h} - [ \be_h \sc (k-h) ] \LP_k \ga_{k-h}
\Big)~.
$$
(In the above: $G + G := \{p+q | ~p, q \in G\}$; $\setminus$ is the
set-theoretical difference; again, $\ga_{k-h} := 0$ and $\be_{k-h} := 0$
if $k - h \not\in G$.)\par \noindent
For each $p \in \reali$, the differential error has the tautological estimator
\beq \ep_p :=
\| \bee(\ug) \|_p = \sqrt{\sum_{k \in dG} |k|^{2 p} \big(|\rho_k|^2 + |\sigma_k|^2 \big)}~;
\label{nora} \feq
there is a rougher estimator, depending on a second real number $q \geqs p$, of the form
\beq \ep'_{p q} := {\Kp_{q} \over |G|^{q-p}} \| \ug \|_{q} \| \ug \|_{q+1}~\label{norb} \feq
where $\Kp_q \in (0,+\infty)$ is constant fulfilling \rref{basic2} with $p$ replaced
by $q$ (i.e., $\| \bP(\bv, \bw) \|_q \leqs \Kp_q \| \bv \|_q \| \bw \|_{q+1}$
for all $\bv \in \HD{q}$, $\bw \in \HD{q+1}$).
\end{prop}
\textbf{Proof.}
(i) Eqs. \rref{ovv1} \rref{ovv2} are self-evident. To derive Eq.
\rref{ovv3}, write $\| \ug(0) - \bu_0 \|_{p}  =
\| (1 - \bEG) \bu_0 \|_{p}$ and use the inequality \rref{begmp}. \par \noindent
(ii) Definition \ref{deferrsm} for the differential error
and Eq. \rref{galer} for $\ug$ give
$$ \bee(\ug) = {d \ug \over d t} - \bA \ug - \bP(\ug, \ug) $$
\beq = \bEG \bP(\ug, \ug) - \bP(\ug, \ug) = - (1 - \bEG) \bP(\ug, \ug)~; \label{sum1} \feq
this proves the first equality in \rref{prov1}.
In order to derive the second equality in \rref{prov1} we must compute
the Fourier representation of $(1 - \bEG) \bP(\ug, \ug)$. Let us
start from the equation
\beq \bP(\ug,\ug) = \big(\P(u_G, u_G) - \P(b_G,b_G),\P(u_G, b_G) - \P(b_G,u_G) \big) \feq
and use the Fourier representation \rref{eq2} of $\P$, recalling again that
$(u_G)_h = \ga_h$, $(b_G)_h = \be_h$ for $h \in G$ and $(u_G)_h = (b_G)_h = 0$ for
$h \not\in G$; this readily
gives
\beq \bP(\ug,\ug) = \left(\sum_{k \in \Zd_0} \rho_k e_k~, \sum_{k \in \Zd_0} \sigma_k e_k \right)~,
\label{refor} \feq
where $\rho_k, \sigma_k$ are defined following Eq. \rref{630} for all $k \in \Zd_0$.
Let us consider, for example, the coefficient $\rho_k$, which is a sum over $h \in G$
containing terms of the form $\gamma_{k -h}$ and $\beta_{k-h}$. If $k \not \in G+G$,
for all $h \in G$ we have $k - h \not \in G$ (since $k - h \in G$
would imply $k = (k - h) + h \in G+G$); but
$k-h \not \in G$ implies $\gamma_{k-h} =0$ and $\beta_{k-h}=0$. In conclusion
we have $\rho_k = 0$ for $k \not \in G+G$; for similar reasons we have
$\sigma_k=0$ for $k \not \in G+G$. Summing up, we can reformulate
Eq. \rref{refor} as
\beq \bP(\ug,\ug) = \left(\sum_{k \in (G + G) \setminus \{0 \}}
\rho_k e_k~, \sum_{k \in (G+G) \setminus \{0 \}} \sigma_k e_k \right)~.
\label{refor2} \feq
Application of $1-\bEG$ to the above sums deletes all terms with $k \in G$;
since $(G + G) \setminus (G \cup \{0\})= dG$ (see Eq. \rref{630}), we obtain
\beq (1 - \bEG) \bP(\ug,\ug) = \left(\sum_{k \in dG}
\rho_k e_k~, \sum_{k \in dG} \sigma_k e_k \right)~.
\label{refor3} \feq
Eqs. \rref{sum1} \rref{refor3} fully justify Eq. \rref{prov1}. \parn
Once one has Eq. \rref{prov1}, statement \rref{nora} is obvious.
The subsequent statement \rref{norb} is proved using Eq. \rref{begmp} and
the inequality involving $\Kp_q$, which imply
$$ \| (1 - \bEG) \bP(\ug, \ug) \|_p \leqs {1 \over |G|^{q-p}} \| \bP(\ug, \ug) \|_q
\leqs {\Kp_{q} \over |G|^{q-p}} \| \ug \|_{q} \| \ug \|_{q+1}~. $$
\fine
\begin{rema}
\textbf{Remarks.} We think it is conceptually important to
propose here the analogues of two remarks made in \cite{appeul}
about the Galerkin approach to NS equations. \parn
(i) The ``rough'' error estimator $\ep'_{p q}$ of Eq.
\rref{norb} is determined by the norm
$\| \ug \|_{q} = \big(\sum_{k \in G} |k|^{2 q} (|\ga_k|^2 + |\be_k|^2)\big)^{1/2}$
and by the analogous norm of order $q+1$, whose computation
involves sums over $G$. The tautological estimator
$\ep_p$ of Eq. \rref{nora} is obviously more precise,
but involves a sum over the set $dG$ which is significantly
bigger than $G$. In applications with a large $G$, the sum
over $dG$ becomes too expensive from a computational viewpoint
and one is led to use the rough estimator \rref{norb}. \parn
(ii) The Galerkin equations \rref{galecomp}
are usually solved numerically; of course, this procedure
does not give the exact solution
$(\ga_k, \be_k)$ ($k \in G$) but, rather, some approximant whose
distance from $(\ga_k, \be_k)$ should be estimated. In the application presented
in the next section, relying on a relatively small set $G$ of modes,
we have assumed this distance to be negligible; this viewpoint should be revised
if $G$ were much larger. ({\footnote{To get reliable results for computations in many modes,
one could perhaps use an ODE solver implementing a standard numerical method
and its theoretical error estimates via
a software for certified numerical computations, like \texttt{INTLAB} \cite{intlab} or \texttt{arb} \cite{arb}.
For general considerations on certified computations, including applications
to ODEs, see \cite{Rump}.}})
\ffine
\end{rema}
\section{An application of the Galerkin method}
\label{galespec}
\begin{figure}
\centering
\begin{subfigure}{.49\textwidth}
\includegraphics[width=6cm]{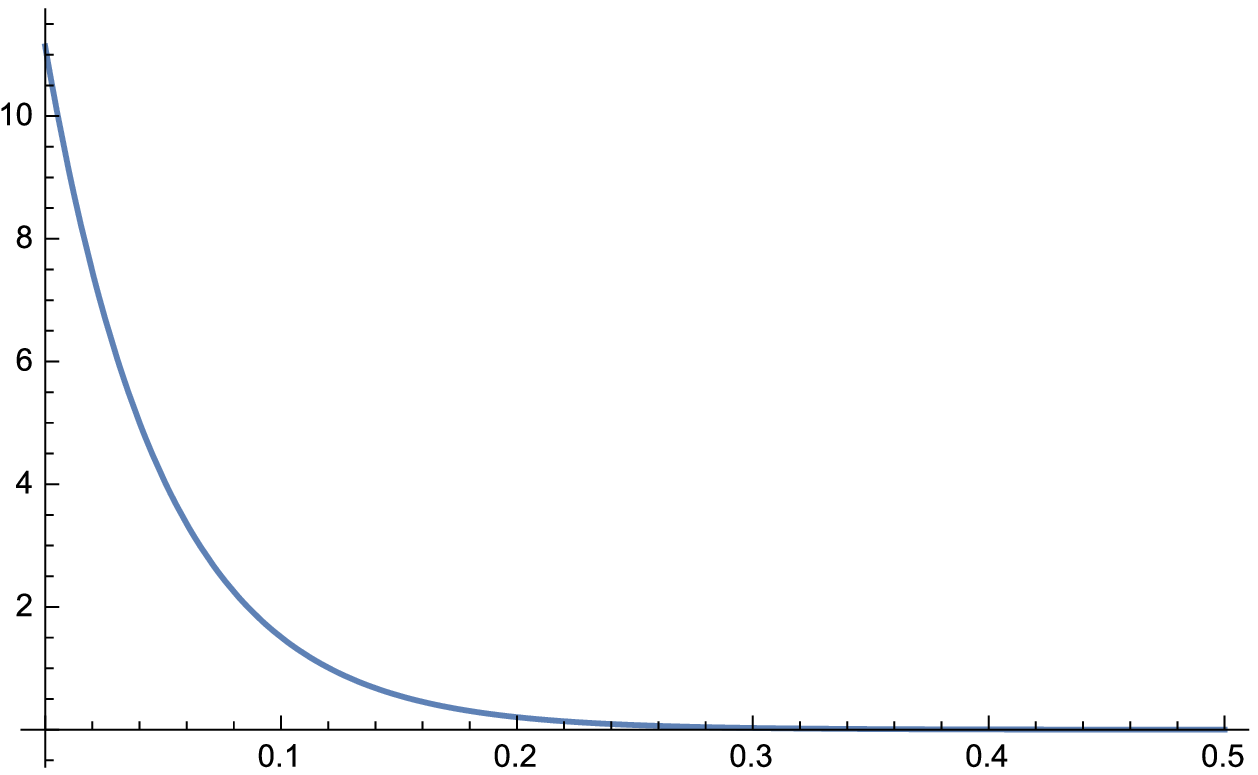}
\caption{ $\mu=20$. Graph of $| \gamma_k (t) |$ for $k=(0,1,0)$ and for $t \in [0,0.5)$.}  \label{fig1a}
\end{subfigure}
\begin{subfigure}{.49\textwidth}
\includegraphics[width=6cm]{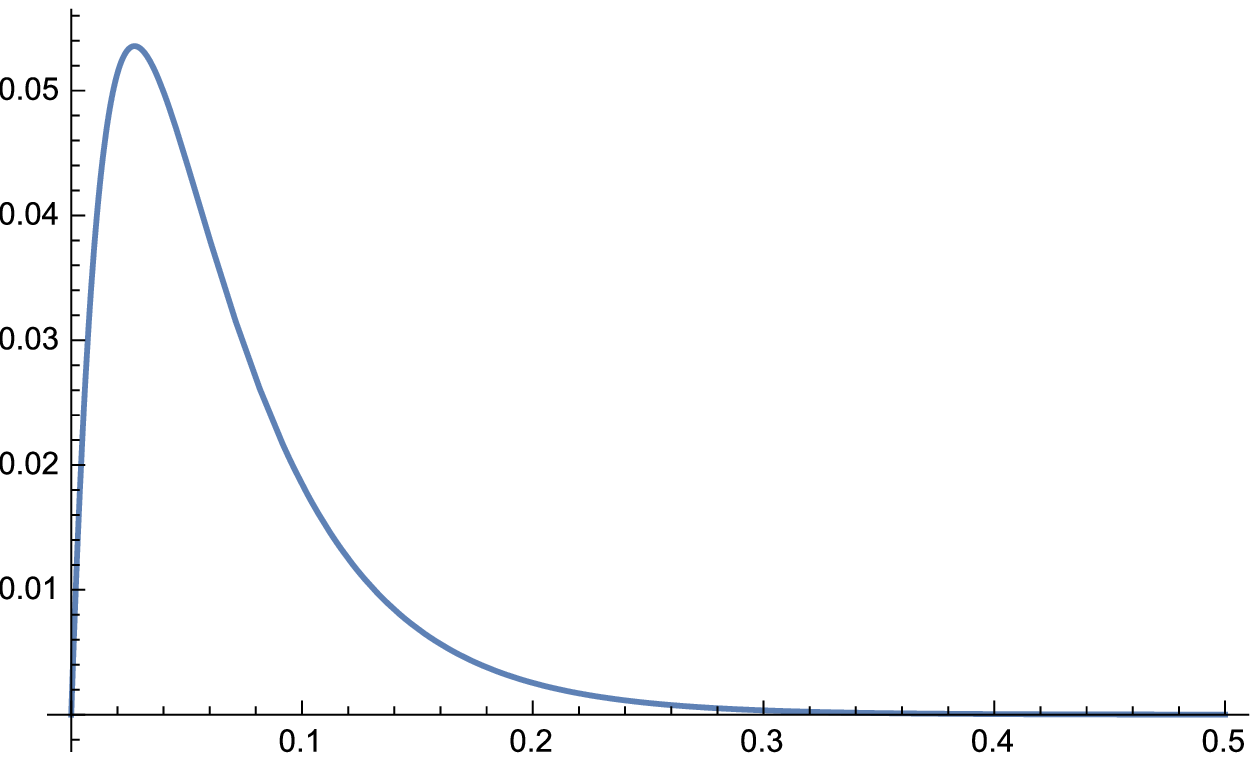}
\caption{$\mu=20$. Graph of $| \beta_k (t) |$ for $k=(0,1,0)$ and for $t \in [0,0.5)$.}  \label{fig1b}
\end{subfigure}\\
\begin{subfigure}{.49\textwidth}
\includegraphics[width=6cm]{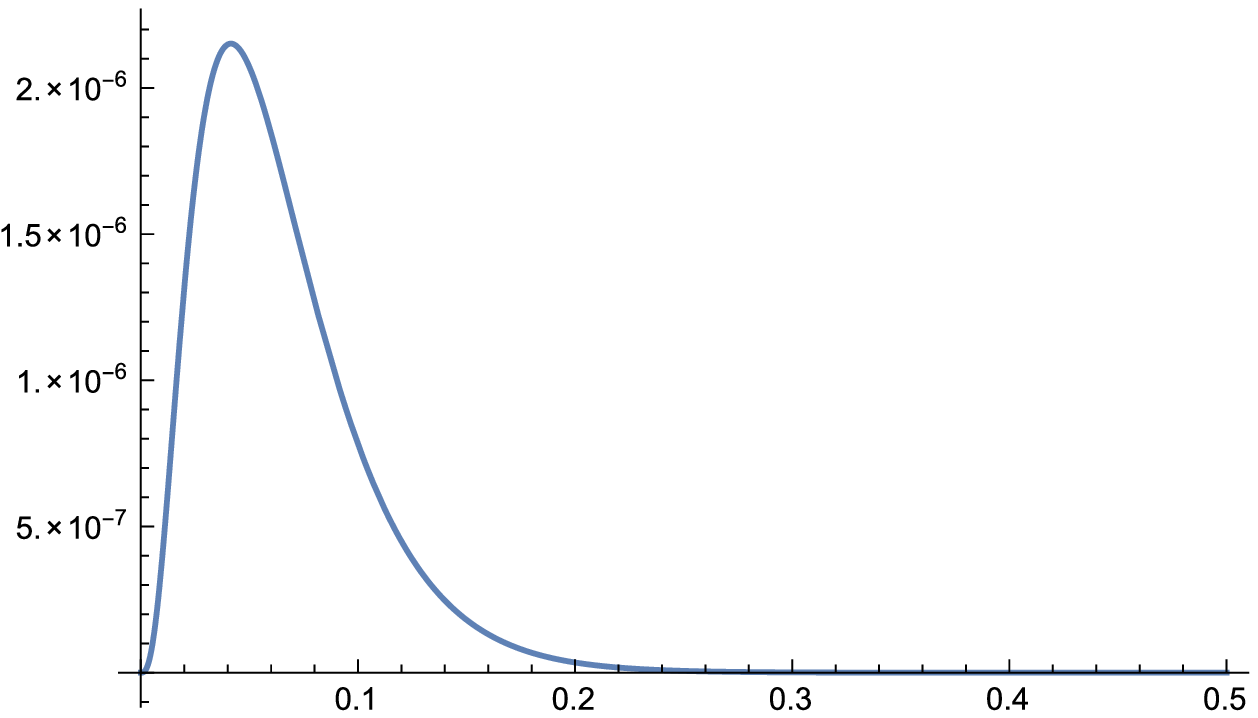}
\caption{$\mu=20$. Graph of $| \gamma_k (t) |$ for $k=(1,1,0)$ and for $t \in [0,0.5)$.}  \label{fig1c}
\end{subfigure}
\begin{subfigure}{.49\textwidth}
\includegraphics[width=6cm]{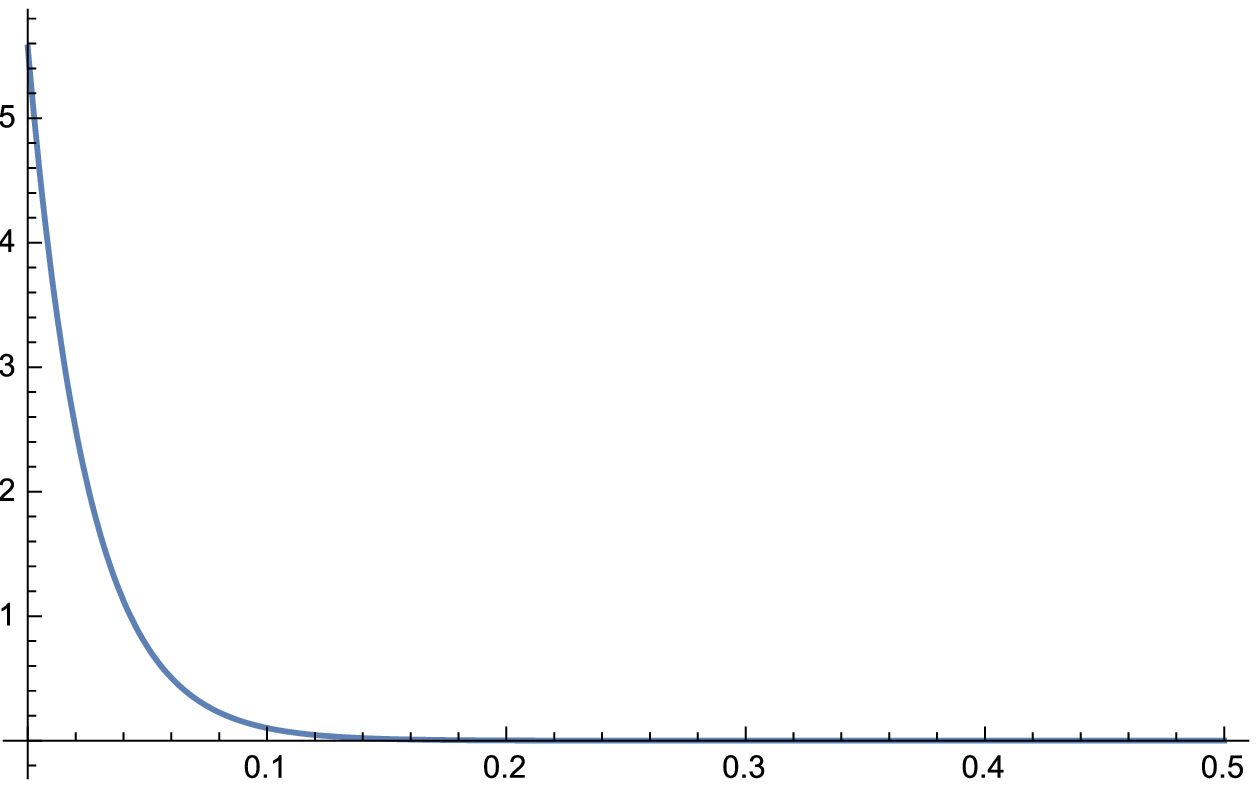}
\caption{$\mu=20$. Graph of $| \beta_k (t) |$ for $k=(1,1,0)$ and for $t \in [0,0.5)$.}  \label{fig1d}
\end{subfigure}\\
\begin{subfigure}{.49\textwidth}
\includegraphics[width=6cm]{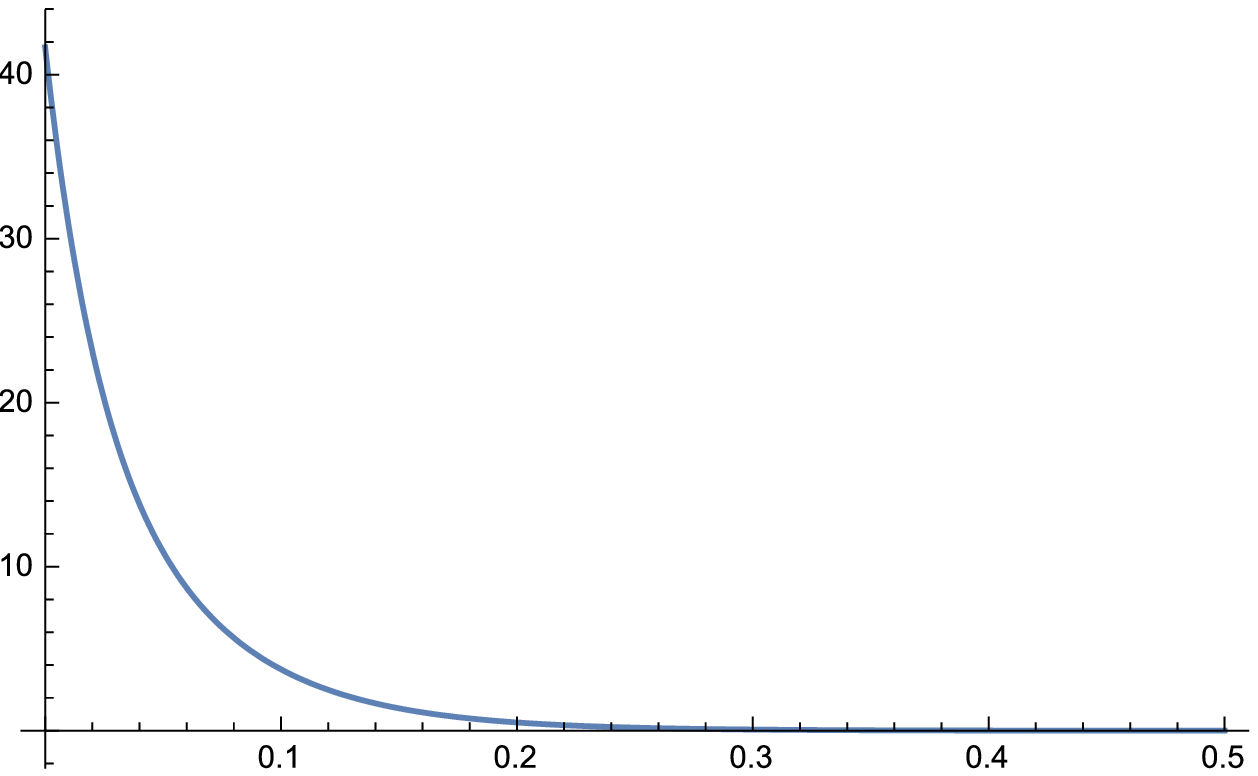}
\caption{$\mu=20$. Graph of $\cald_3 (t)$ for $t \in [0,0.5)$.}  \label{fig1e}
\end{subfigure}
\begin{subfigure}{.49\textwidth}
\includegraphics[width=6cm]{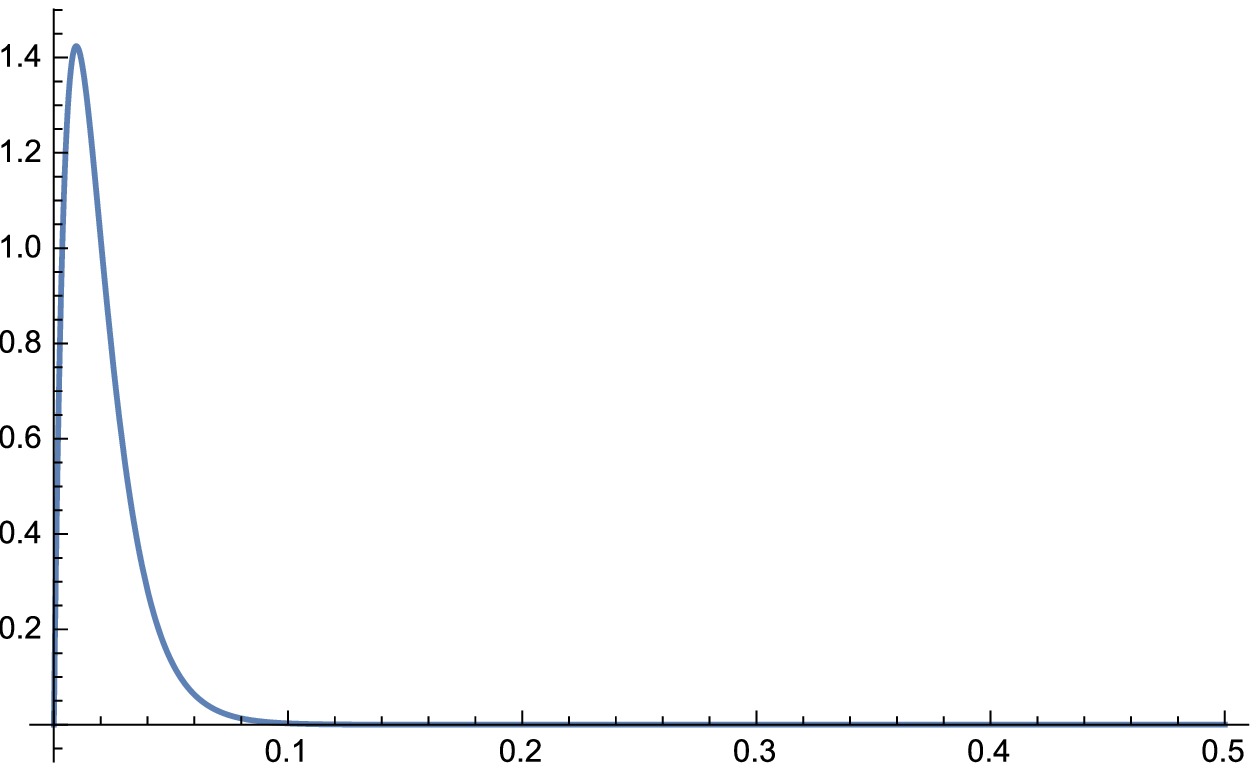}
\caption{$\mu=20$. Graph of $\epsilon_3 (t)$ for $t \in [0,0.5)$.}  \label{fig1f}
\end{subfigure}\\
\begin{subfigure}{.49\textwidth}
\includegraphics[width=6cm]{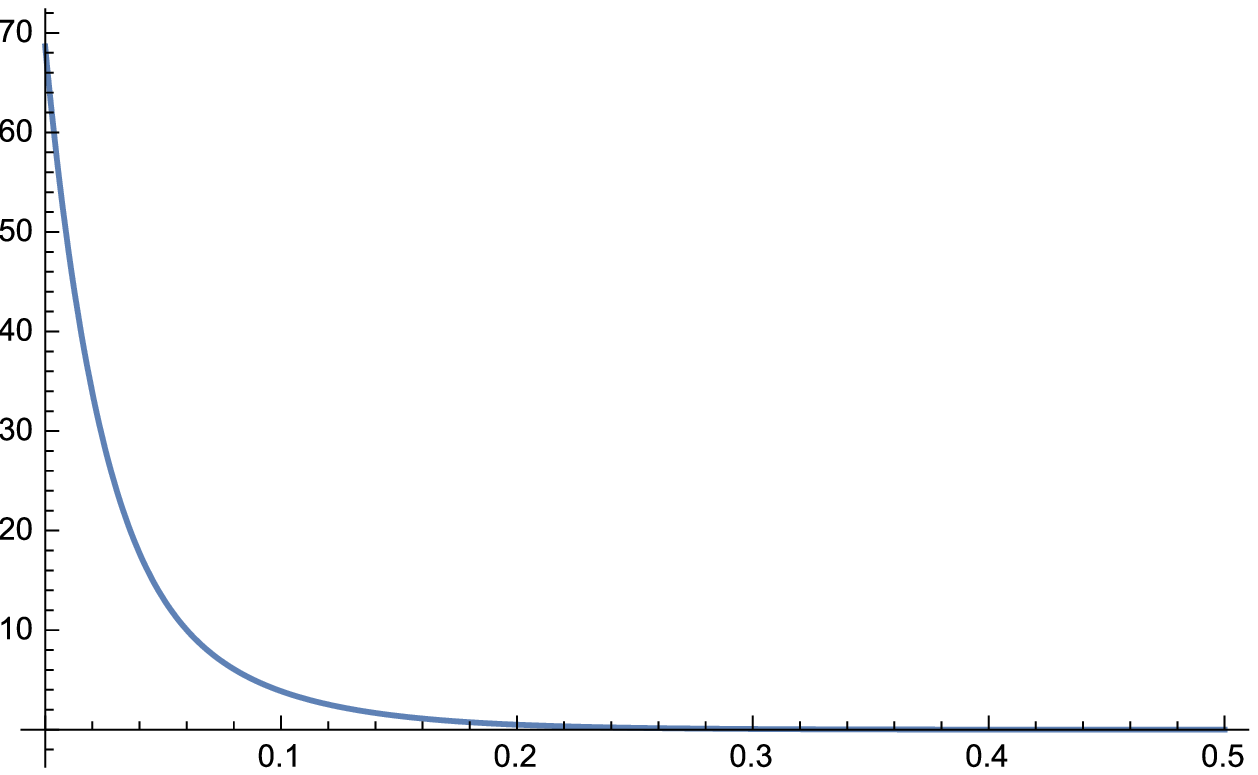}
\caption{$\mu=20$. Graph of $\cald_5 (t)$ for $t \in [0,0.5)$.}  \label{fig1g}
\end{subfigure}
\begin{subfigure}{.49\textwidth}
\includegraphics[width=6cm]{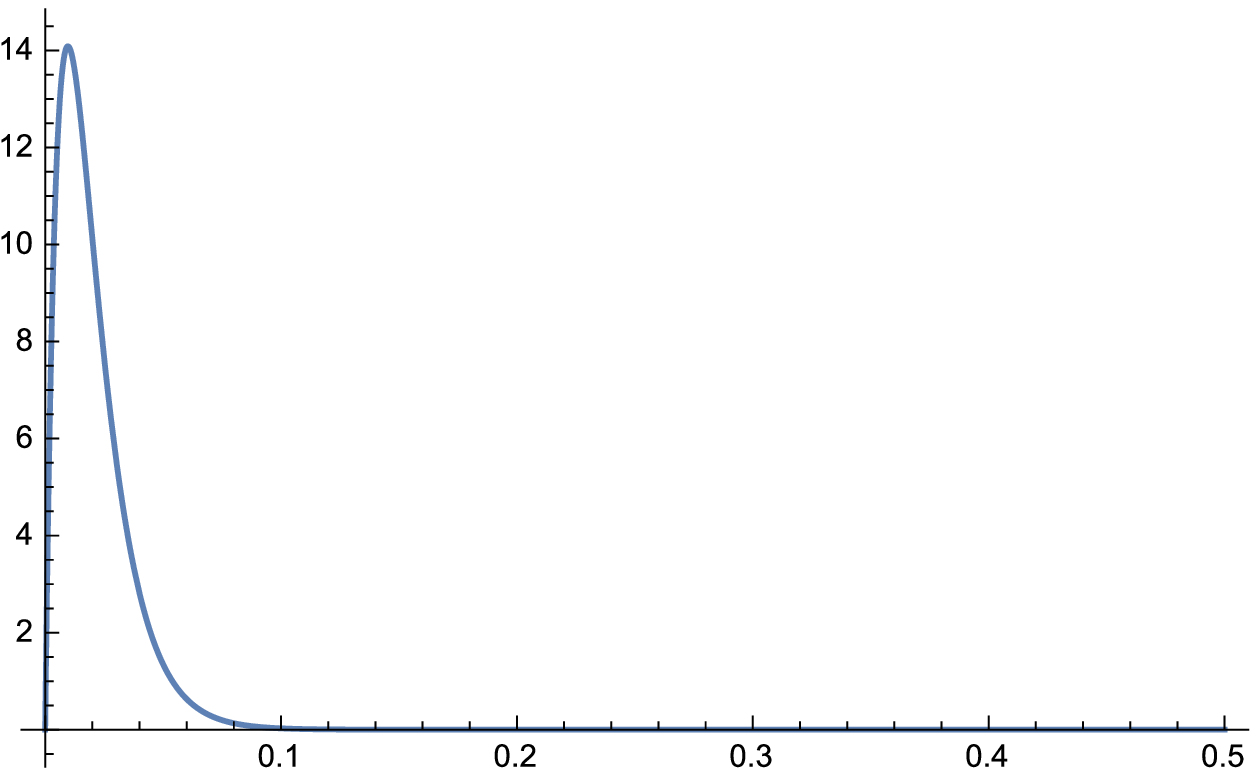}
\caption{$\mu=20$. Graph of $\epsilon_5 (t)$ for $t \in [0,0.5)$.}  \label{fig1h}
\end{subfigure}\\
\begin{subfigure}{.49\textwidth}
\includegraphics[width=6cm]{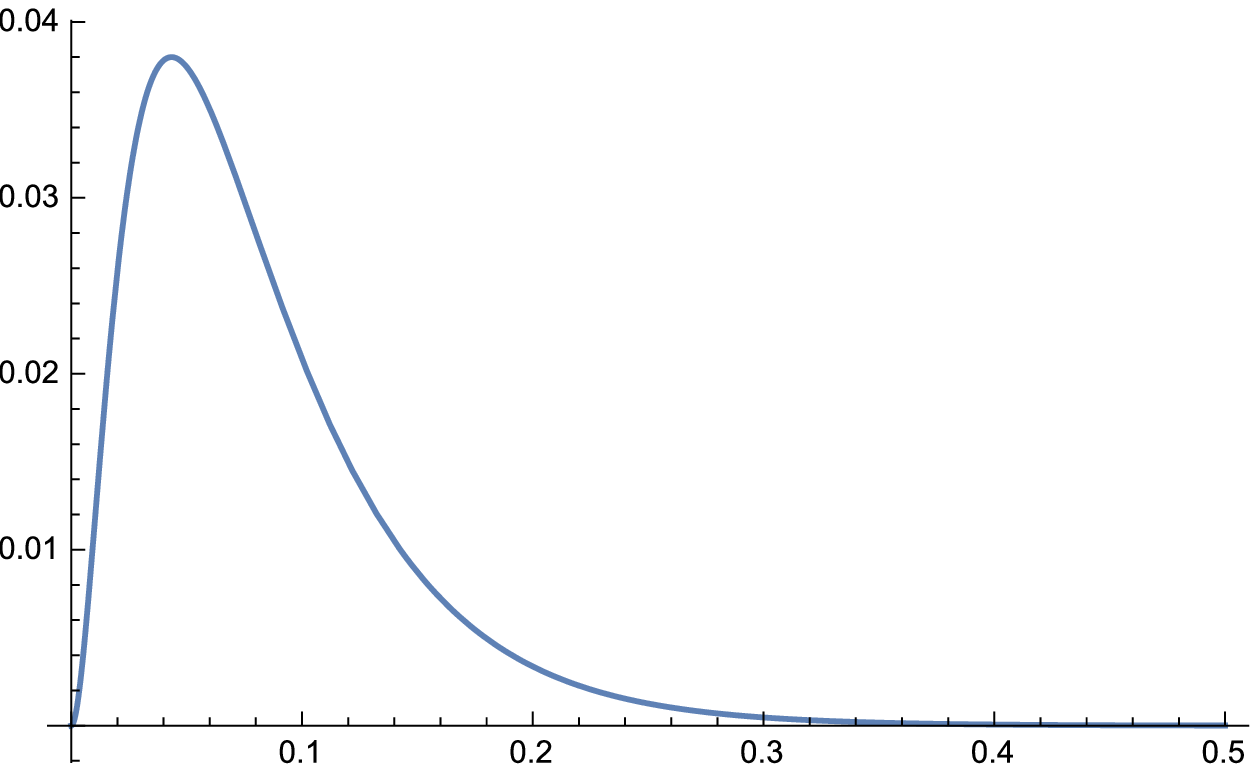}
\caption{$\mu=20$. Graph of $\RR_3 (t)$ for $t \in [0,0.5)$.}  \label{fig1i}
\end{subfigure}
\begin{subfigure}{.49\textwidth}
\includegraphics[width=6cm]{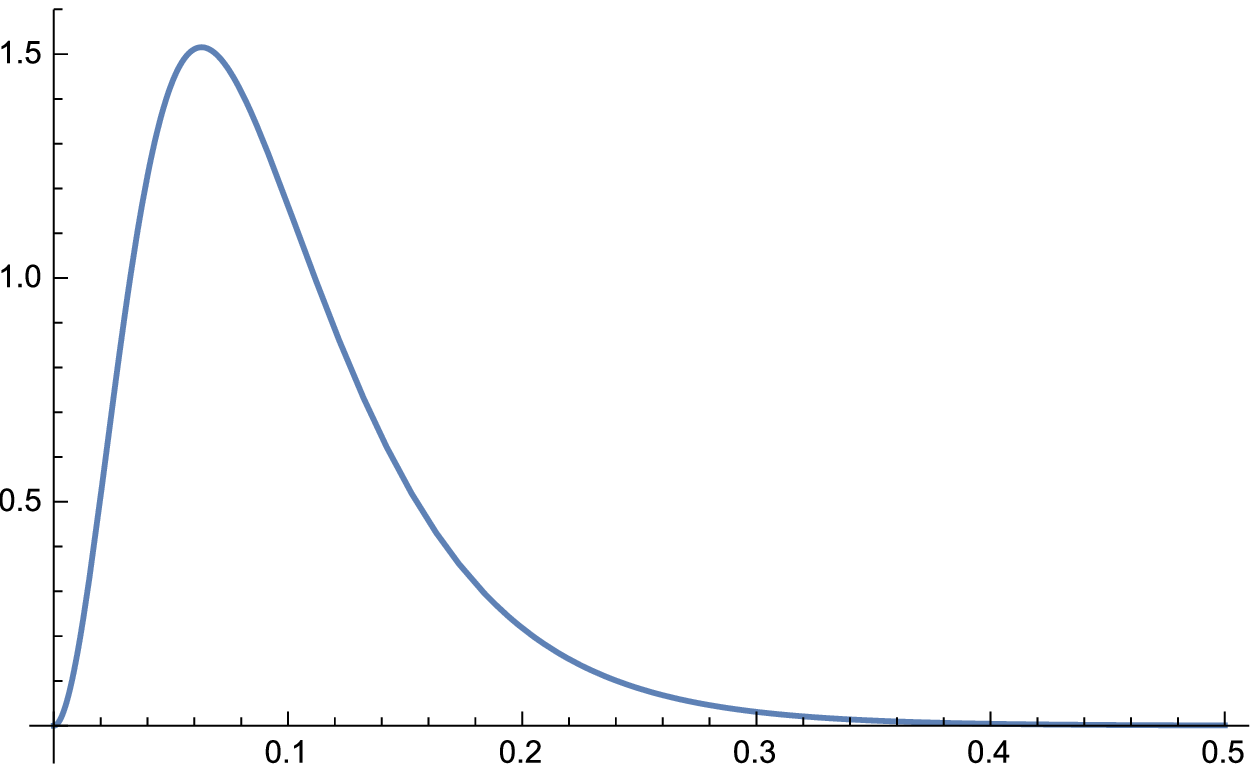}
\caption{$\mu=20$. Graph of $\RR_5 (t)$ for $t \in [0,0.5)$.}  \label{fig1j}
\end{subfigure}
\caption{Plots related to the case $\mu=20$.}
\end{figure}

\begin{figure}
\centering
\begin{subfigure}{.49\textwidth}
\includegraphics[width=6cm]{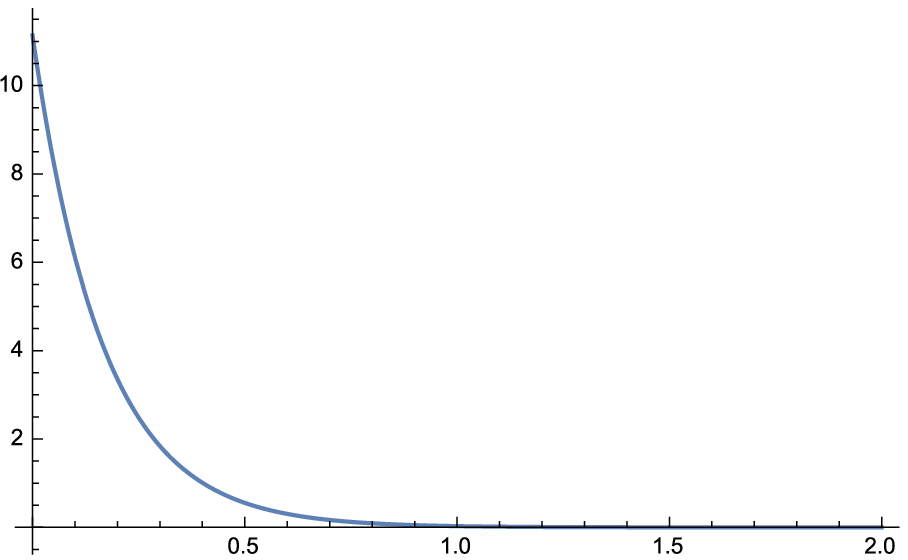}
\caption{ $\mu=6$. Graph of $| \gamma_k (t) |$ for $k=(0,1,0)$ and for $t \in [0,2)$.}  \label{fig2a}
\end{subfigure}
\begin{subfigure}{.49\textwidth}
\includegraphics[width=6cm]{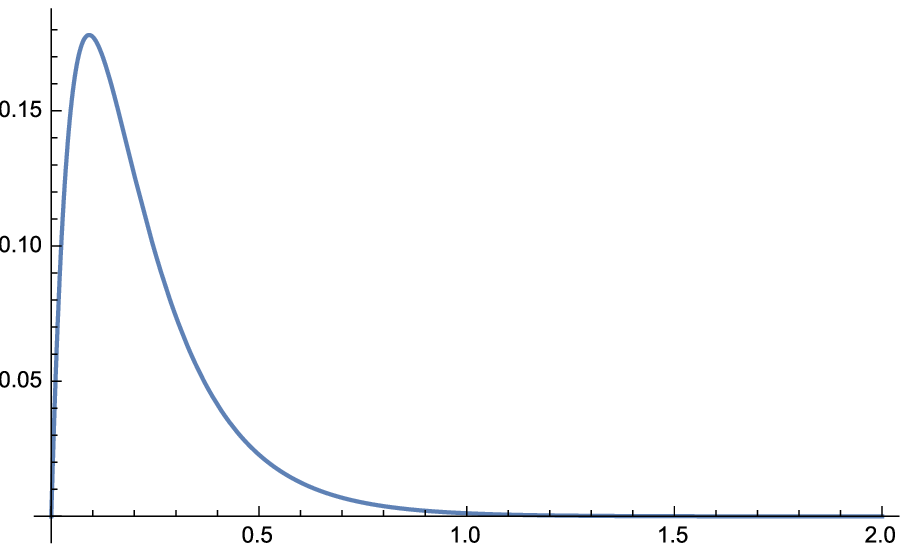}
\caption{$\mu=6$. Graph of $| \beta_k (t) |$ for $k=(0,1,0)$ and for $t \in [0,2)$.}  \label{fig2b}
\end{subfigure}\\
\begin{subfigure}{.49\textwidth}
\includegraphics[width=6cm]{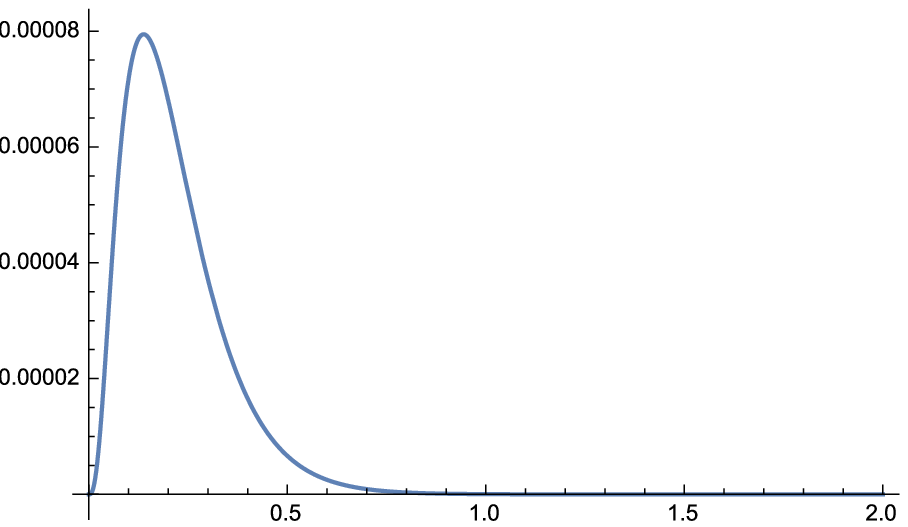}
\caption{$\mu=6$. Graph of $| \gamma_k (t) |$ for $k=(1,1,0)$ and for $t \in [0,2)$.}  \label{fig2c}
\end{subfigure}
\begin{subfigure}{.49\textwidth}
\includegraphics[width=6cm]{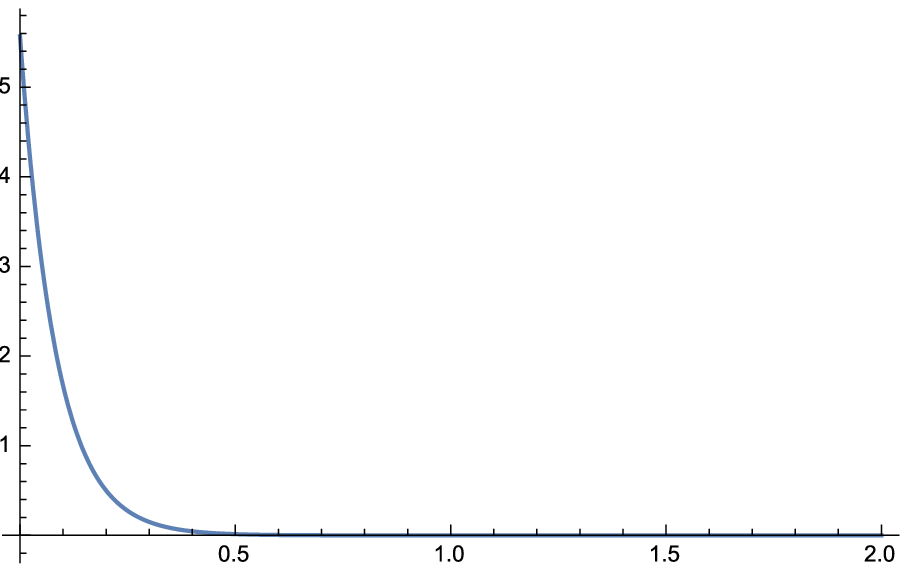}
\caption{$\mu=6$. Graph of $| \beta_k (t) |$ for $k=(1,1,0)$ and for $t \in [0,2)$.}  \label{fig2d}
\end{subfigure}\\
\begin{subfigure}{.49\textwidth}
\includegraphics[width=6cm]{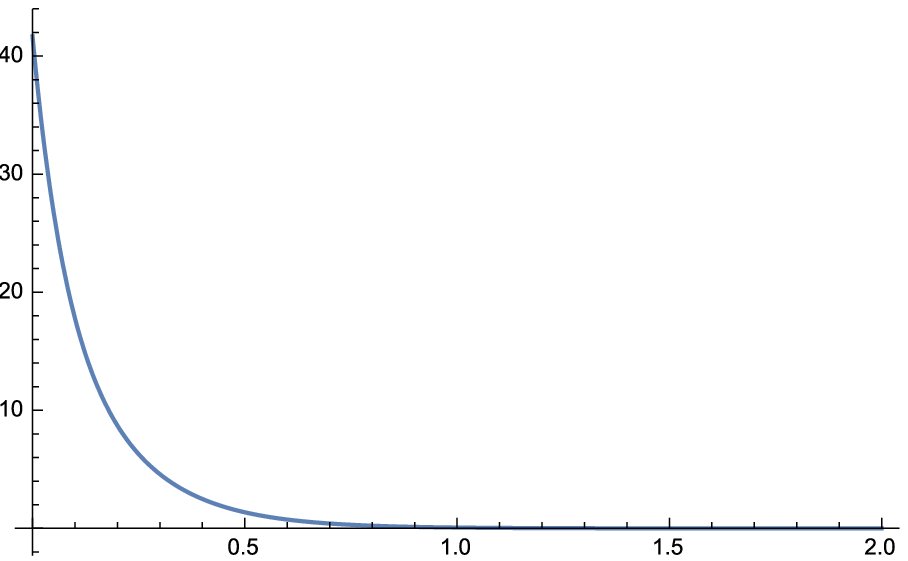}
\caption{$\mu=6$. Graph of $\cald_3 (t)$ for $t \in [0,2)$.}  \label{fig2e}
\end{subfigure}
\begin{subfigure}{.49\textwidth}
\includegraphics[width=6cm]{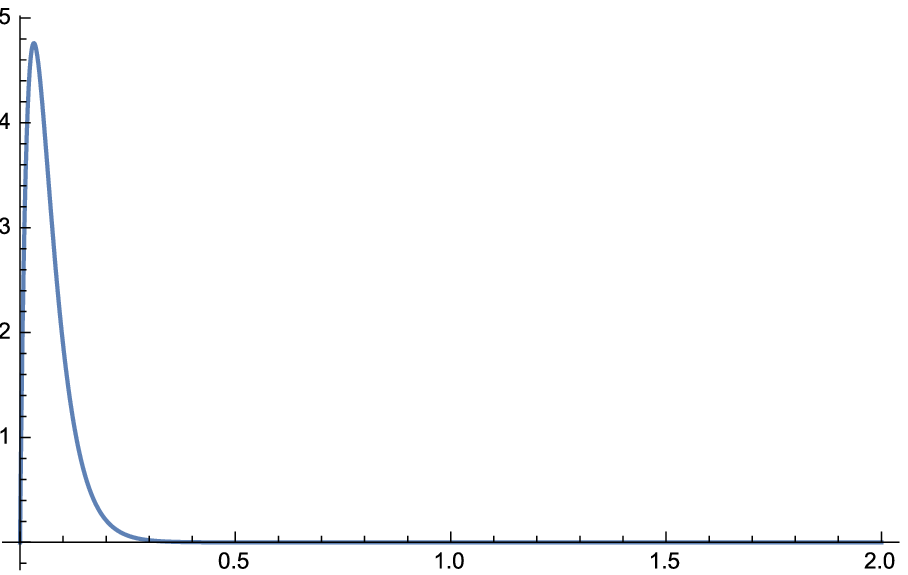}
\caption{$\mu=6$. Graph of $\epsilon_3 (t)$ for $t \in [0,2)$.}  \label{fig2f}
\end{subfigure}\\
\begin{subfigure}{.49\textwidth}
\includegraphics[width=6cm]{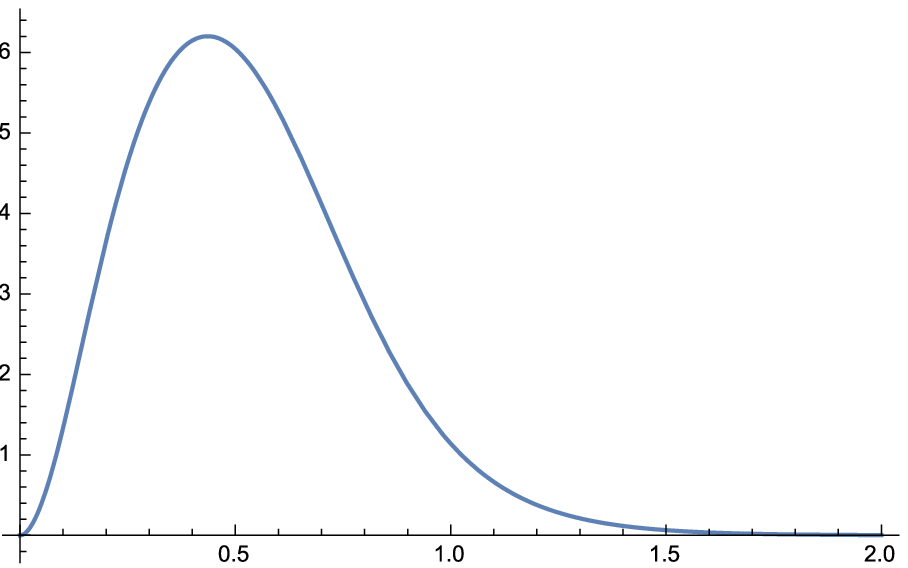}
\caption{$\mu=6$. Graph of $\RR_3 (t)$.}  \label{fig2g}
\end{subfigure}
\caption{Plots related to the case $\mu=6$.}
\end{figure}

\begin{figure}
\centering
\begin{subfigure}{.49\textwidth}
\includegraphics[width=6cm]{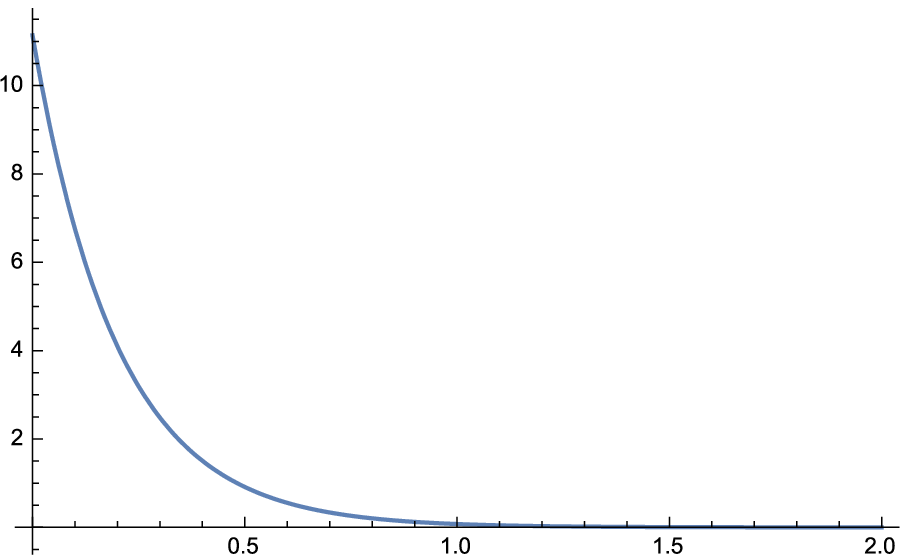}
\caption{ $\mu=5$. Graph of $| \gamma_k (t) |$ for $k=(0,1,0)$ and for $t \in [0,2)$.}  \label{fig3a}
\end{subfigure}
\begin{subfigure}{.49\textwidth}
\includegraphics[width=6cm]{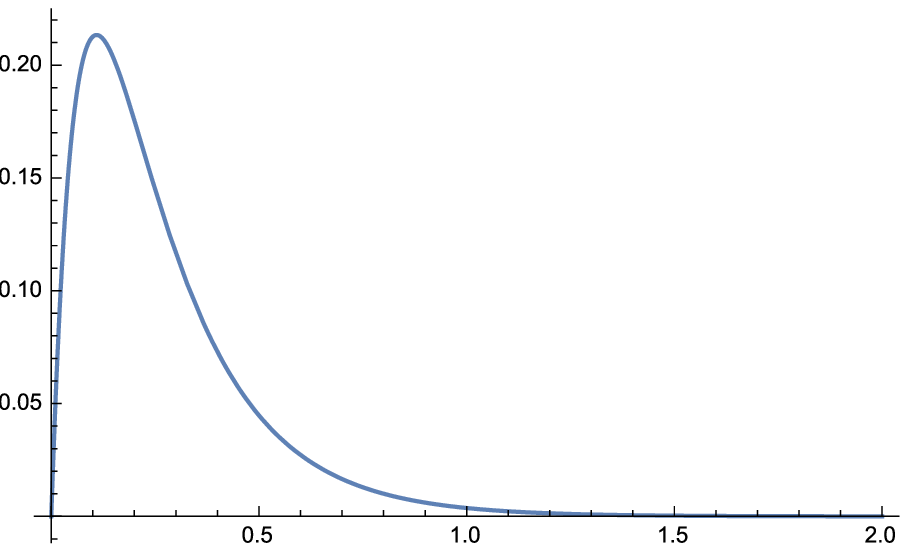}
\caption{$\mu=5$. Graph of $| \beta_k (t) |$ for $k=(0,1,0)$ and for $t \in [0,2)$.}  \label{fig3b}
\end{subfigure}\\
\begin{subfigure}{.49\textwidth}
\includegraphics[width=6cm]{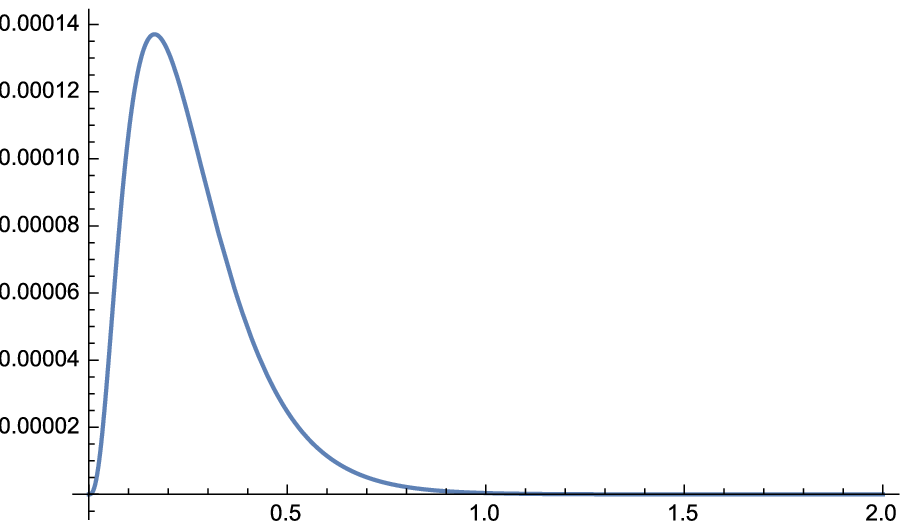}
\caption{$\mu=5$. Graph of $| \gamma_k (t) |$ for $k=(1,1,0)$ and for $t \in [0,2)$.}  \label{fig3c}
\end{subfigure}
\begin{subfigure}{.49\textwidth}
\includegraphics[width=6cm]{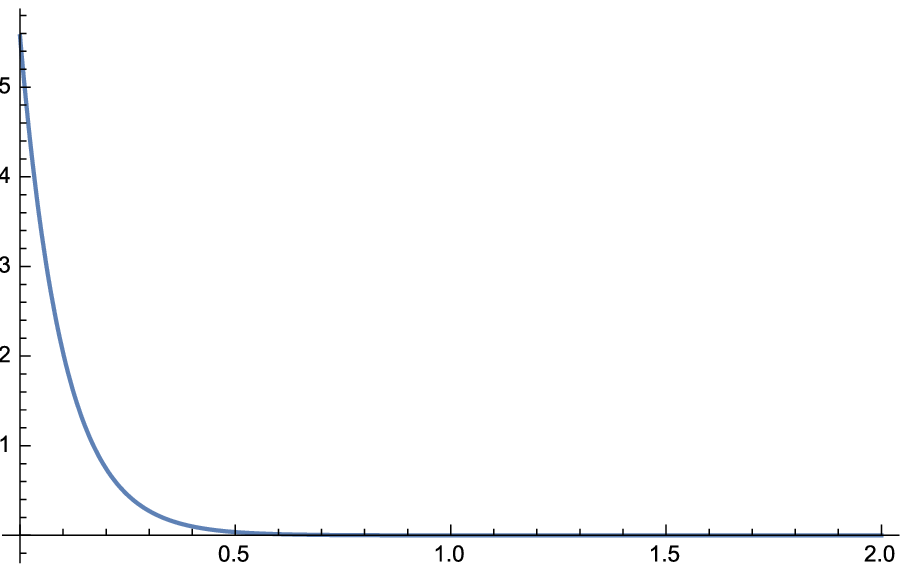}
\caption{$\mu=5$. Graph of $| \beta_k (t) |$ for $k=(1,1,0)$ and for $t \in [0,2)$.}  \label{fig3d}
\end{subfigure}\\
\begin{subfigure}{.49\textwidth}
\includegraphics[width=6cm]{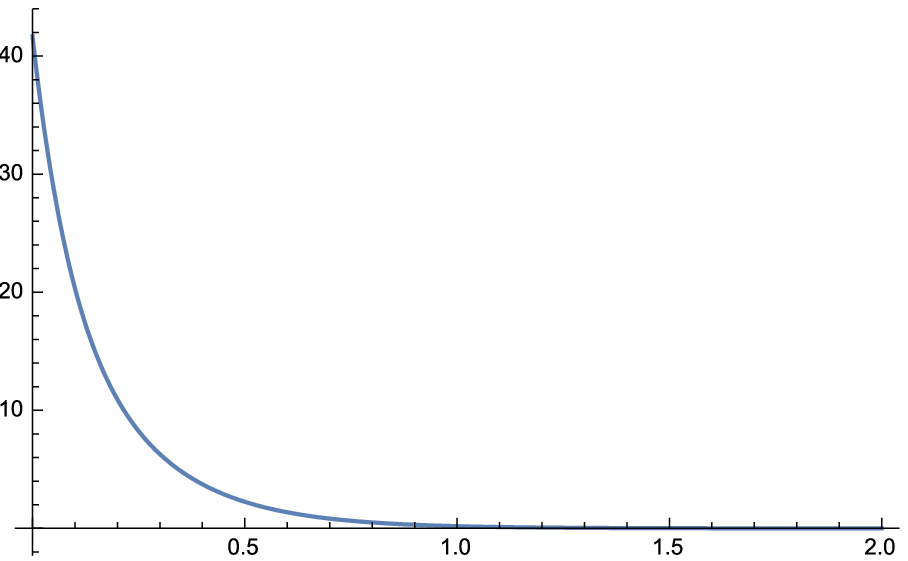}
\caption{$\mu=5$. Graph of $\cald_3 (t)$ for $t \in [0,2)$.}  \label{fig3e}
\end{subfigure}
\begin{subfigure}{.49\textwidth}
\includegraphics[width=6cm]{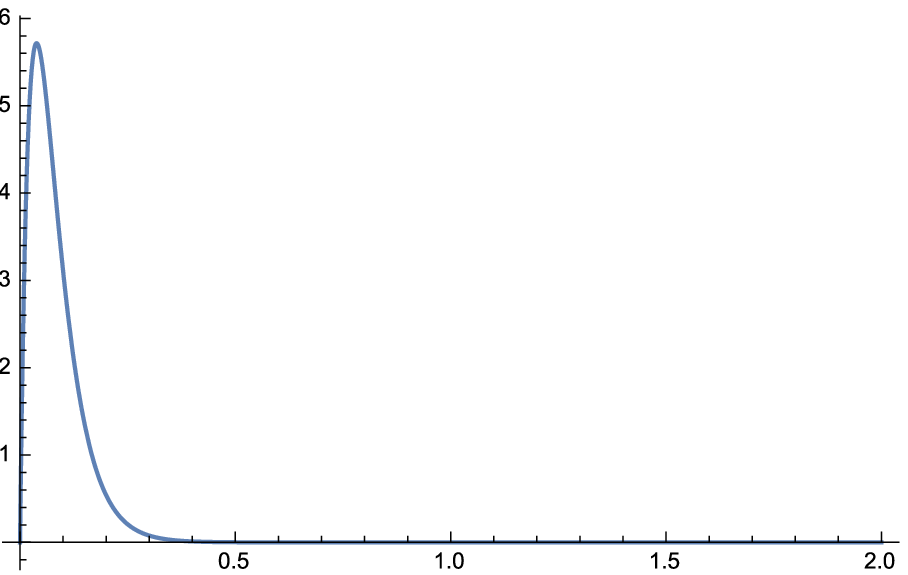}
\caption{$\mu=5$. Graph of $\epsilon_3 (t)$ for $t \in [0,2)$.}  \label{fig3f}
\end{subfigure}\\
\begin{subfigure}{.49\textwidth}
\includegraphics[width=6cm]{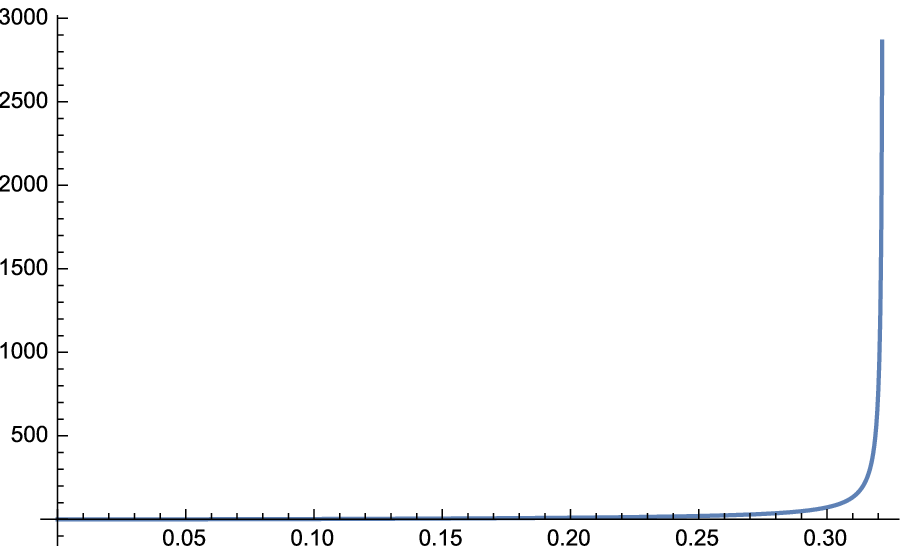}
\caption{$\mu=5$. Graph of $\RR_3 (t)$. This function diverges as $t\rightarrow T_c$ with $T_c= 0.3238...$.}  \label{fig3g}
\end{subfigure}
\caption{Plots related to the case $\mu=5$.}
\end{figure}

\begin{figure}
\centering
\begin{subfigure}{.49\textwidth}
\includegraphics[width=6cm]{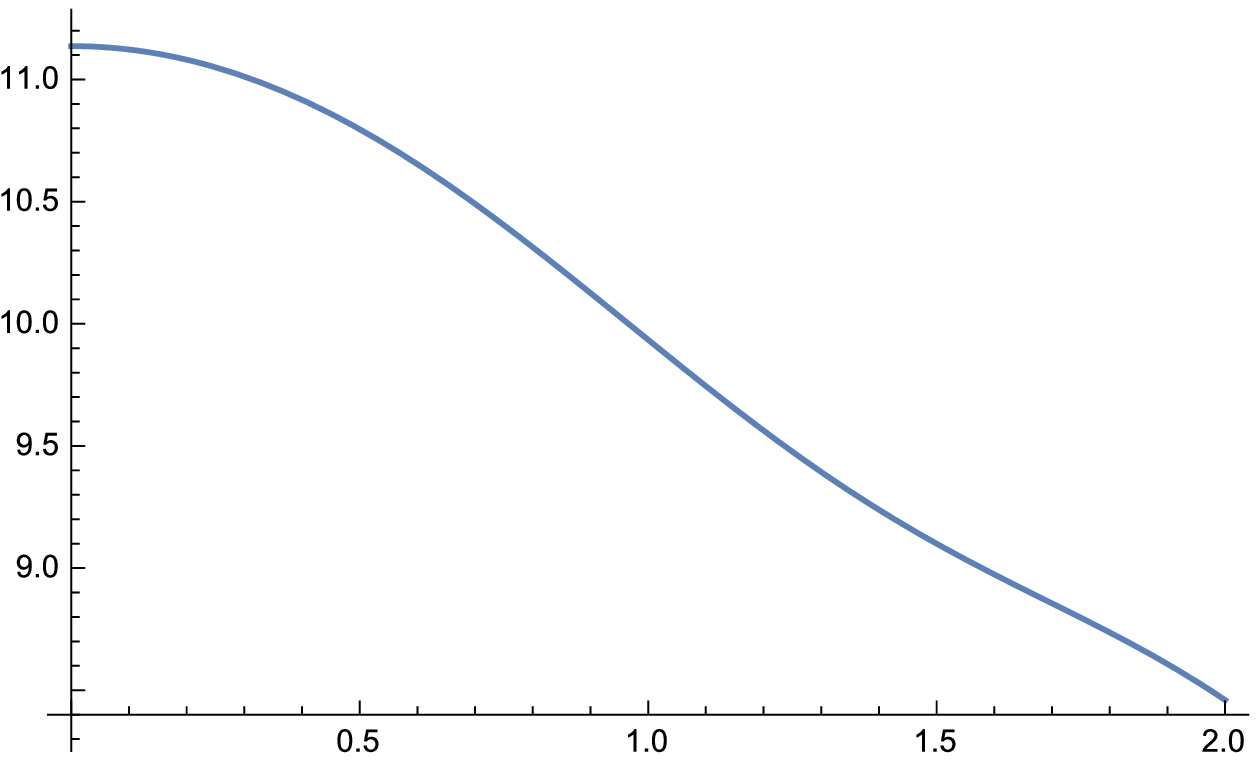}
\caption{ $\mu=0$. Graph of $| \gamma_k (t) |$ for $k=(0,1,0)$ and for $t \in [0,2)$.}  \label{fig4a}
\end{subfigure}
\begin{subfigure}{.49\textwidth}
\includegraphics[width=6cm]{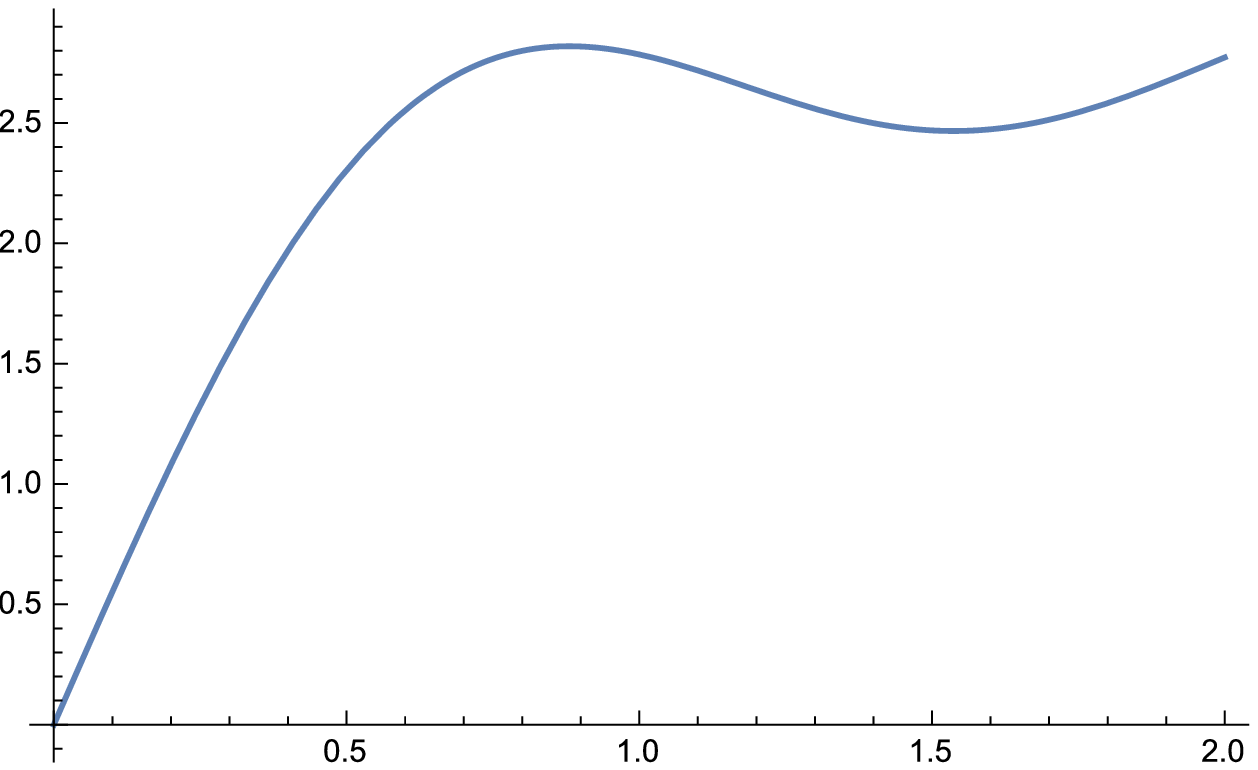}
\caption{$\mu=0$. Graph of $| \beta_k (t) |$ for $k=(0,1,0)$ and for $t \in [0,2)$.}  \label{fig4b}
\end{subfigure}\\
\begin{subfigure}{.49\textwidth}
\includegraphics[width=6cm]{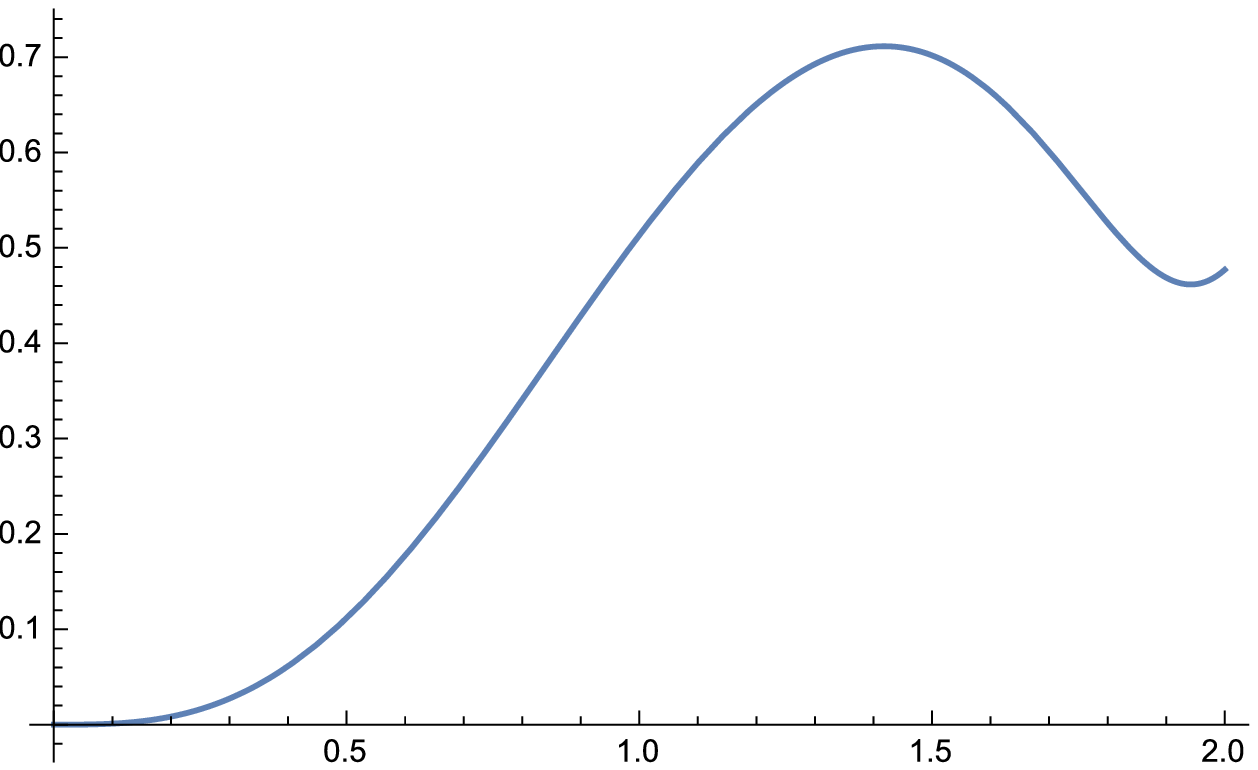}
\caption{$\mu=0$. Graph of $| \gamma_k (t) |$ for $k=(1,1,0)$ and for $t \in [0,2)$.}  \label{fig4c}
\end{subfigure}
\begin{subfigure}{.49\textwidth}
\includegraphics[width=6cm]{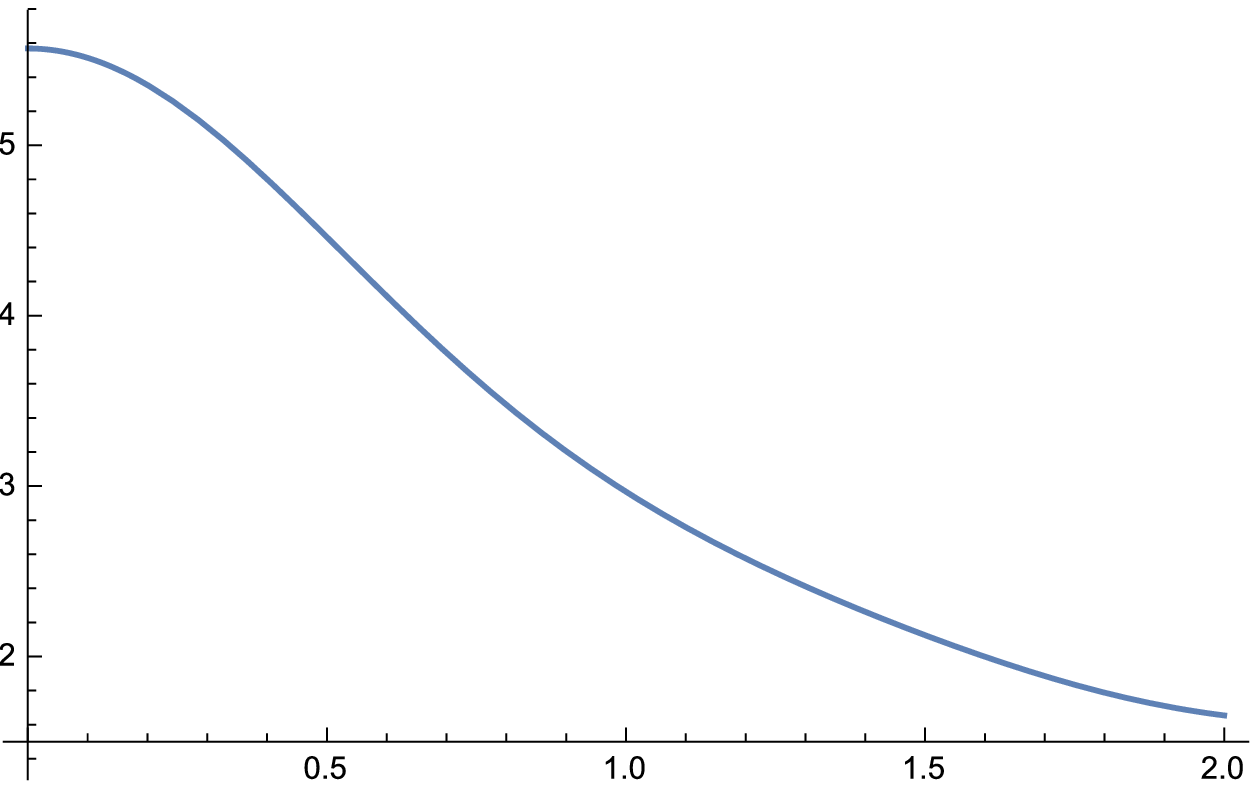}
\caption{$\mu=0$. Graph of $| \beta_k (t) |$ for $k=(1,1,0)$ and for $t \in [0,2)$.}  \label{fig4d}
\end{subfigure}\\
\begin{subfigure}{.49\textwidth}
\includegraphics[width=6cm]{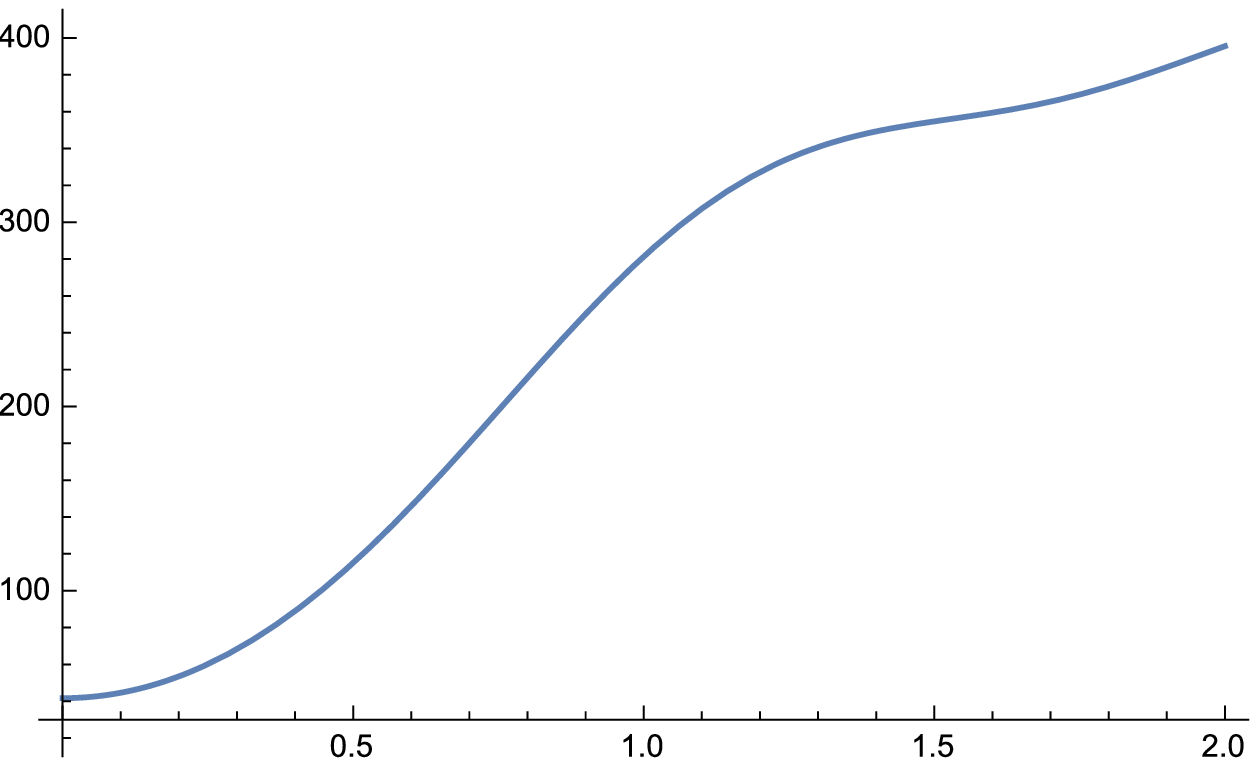}
\caption{$\mu=0$. Graph of $\cald_3 (t)$ for $t \in [0,2)$.}  \label{fig4e}
\end{subfigure}
\begin{subfigure}{.49\textwidth}
\includegraphics[width=6cm]{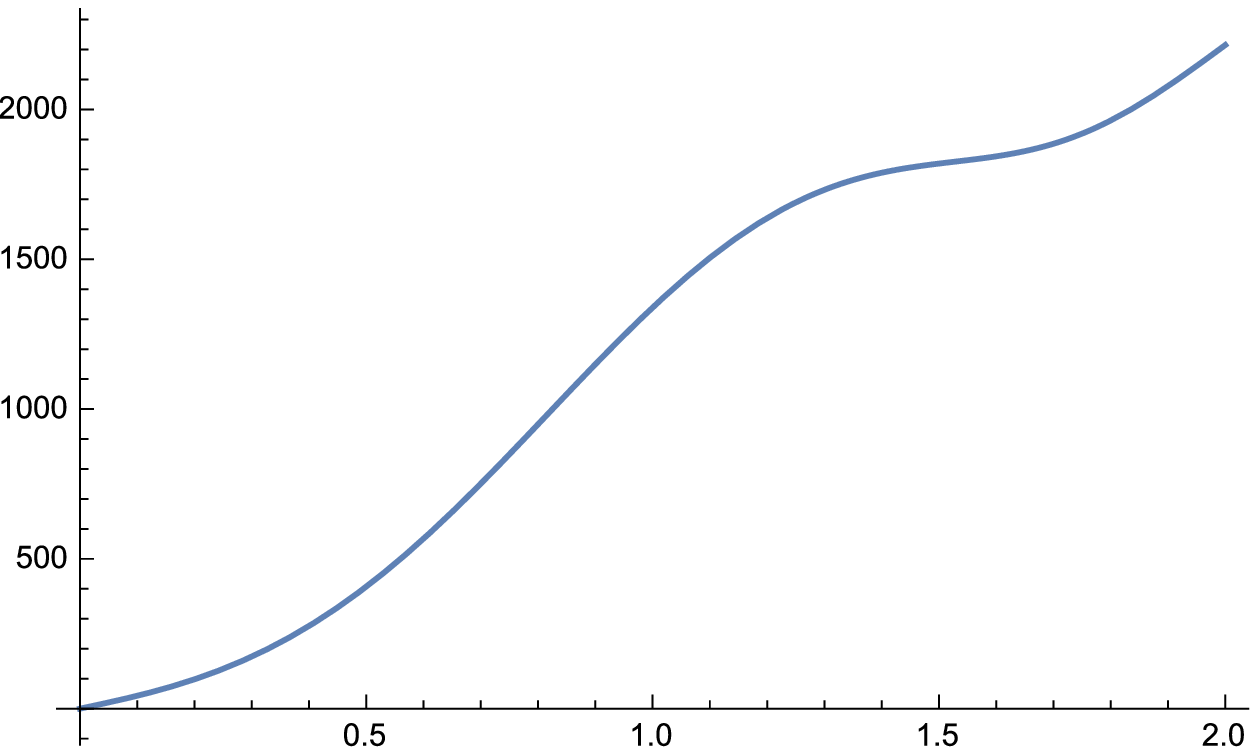}
\caption{$\mu=0$. Graph of $\epsilon_3 (t)$ for $t \in [0,2)$.}  \label{fig4f}
\end{subfigure}\\
\begin{subfigure}{.49\textwidth}
\includegraphics[width=6cm]{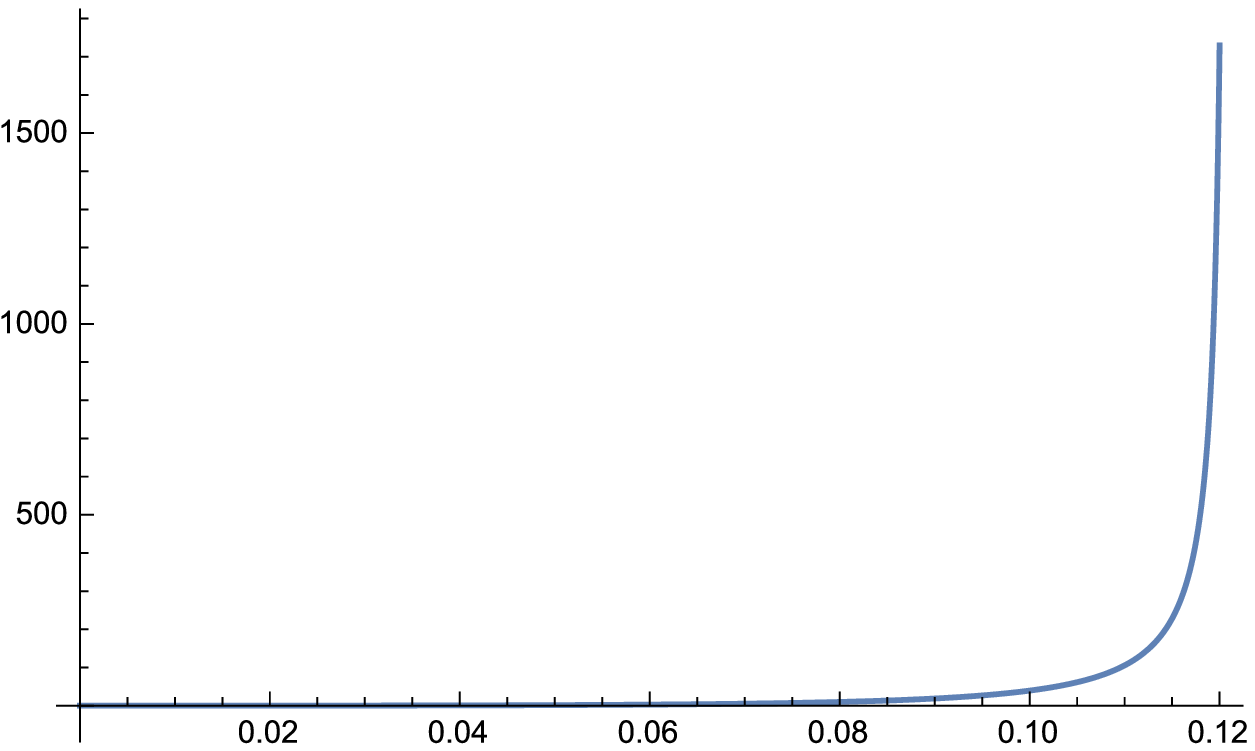}
\caption{$\mu=0$. Graph of $\RR_3 (t)$. This function diverges as $t\rightarrow T_c$ with $T_c= 0.1211...$.}  \label{fig4g}
\end{subfigure}
\caption{Plots related to the case $\mu=0$.}
\end{figure}

In this section we apply the general framework developed in Secs. \ref{prelim}-\ref{approxmhd} adopting, as approximate solutions,
the Galerkin solutions described in Sec. \ref{galegen}.
We will work in space dimension
\beq
d=3.
\eeq
Moreover, we will specialize our analysis
to the case where the dimensionless viscosity and resistivity are equal:
\beq
\nu=\eta \equiv \mu \in [0,+\infty)~.
\eeq
We will choose as initial datum the ABC flow with perturbing magnetic field, given by
Eqs. (\ref{abcxyz})-(\ref{abcdet}), with the following
values for the parameters appearing therein:
\beq
A=B=C=D=1.
\eeq
The previous sections frequently refer to a basic Sobolev order
$n > d/2 + 1$; here we will take
\beq
n=3.
\eeq
We remark that, with our choice for the initial datum, one has
\beq
\lVert \bu_0  \rVert_3=41.6695...
\eeq
and the criterion \rref{muglob} grants global existence
for the solution of the Cauchy problem \rref{cmhd} if
\beq \mu \geqs 25.9~. \label{resmu} \feq
As shown hereafter, the use of a Galerkin approximant for
this Cauchy problem allows, amongst else, to improve
significantly the bound \rref{resmu}. \parn
So, let us consider the Galerkin approximate
solution $\ug (t)=(u_G(t), b_G(t))$ for a suitable,
finite set $G$ of Fourier modes. Following Eq. (\ref{fur1}), we write
\beq
u_G(t) = \sum_{k \in G} \ga_{k}(t) e_k ~, \qquad b_G(t) = \sum_{k \in G} \be_{k}(t) e_k
\quad \mbox{for $0 \leqs t < + \infty$}. \label{fur1bis}
\eeq
with $\gamma_{k}, \beta_k \in C^\infty([0,+\infty),\complessi^3)$, to be determined.
We choose $G := \{ k = (k_1,k_2, k_3) \in \Zt_0~|~-2 \leqs k_1,k_2,k_3 \leqs 2\}$; this set
consists of $124$ modes and admits the representation
\beq
G:=S \cup -S, \qquad -S:= \{ -k : k \in S\},
\label{modesgs} \eeq
where $S$ is the following set of $62$ modes:
\beq
\begin{split}
& S:=\{ (0, 0, 1), (0, 0, 2), (0, 1, -2), (0, 1, -1), (0, 1, 0), (0, 1, 1),
 (0, 1, 2), (0, 2, -2), \\
 & (0, 2, -1), (0, 2, 0), (0, 2, 1), (0, 2, 2),
 (1, -2, -2), (1, -2, -1), (1, -2, 0), (1, -2, 1), \\
 & (1, -2, 2), (1, -1, -2),
 (1, -1, -1), (1, -1, 0), (1, -1, 1), (1, -1, 2), (1, 0, -2), (1, 0, -1), \\
 & (1, 0, 0), (1, 0, 1), (1, 0, 2), (1, 1, -2), (1, 1, -1), (1, 1, 0),
 (1, 1, 1), (1, 1, 2), \\
 & (1, 2, -2), (1, 2, -1), (1, 2, 0), (1, 2, 1),
 (1, 2, 2), (2, -2, -2), (2, -2, -1), (2, -2, 0), \\
 & (2, -2, 1), (2, -2, 2),
 (2, -1, -2), (2, -1, -1), (2, -1, 0), (2, -1, 1), (2, -1, 2), (2, 0, -2),\\
& (2, 0, -1), (2, 0, 0), (2, 0, 1), (2, 0, 2), (2, 1, -2), (2, 1, -1),
 (2, 1, 0), (2, 1, 1), \\
 & (2, 1, 2), (2, 2, -2), (2, 2, -1), (2, 2, 0),
 (2, 2, 1), (2, 2, 2) \}
 \end{split}
\label{modess} \eeq
The Galerkin approximation and its
implications about the exact solution
$\bu$ of the MHD Cauchy problem \rref{cmhd}
have been considered for several
values of $\mu$ between $0$ and $20$,
following for each $\mu$ the scheme
(i)(ii)(iii) described hereafter; the related
numerical computations have been implemented
using Mathematica on a PC. Here is
the scheme, for a given value of $\mu$. \parn
(i) First of all, the Galerkin approximate
solution $\bu_G$ is computed numerically on a finite time interval $[0,T_F)$,
for the set of modes $G$ in Eqs. \rref{modesgs}
\rref{modess}; this amounts to solve numerically on $[0,T_F)$ the
system of equations \rref{galecomp}
for the unknowns $\ga_k, \be_k$.
Due to the relations $\gamma_{-k}(t) = \overline{\gamma}_{k}(t)$
and $\be_{-k}(t) = \overline{\be}_{k}(t)$, known from Section \ref{galegen}, the computation
is reduced to modes $k \in S$.
\parn
In our numerical computations, $T_F$ is  between
$0.5$ and $2$ (more details on this are given
in the sequel); the CPU time required to solve
the system \rref{galecomp} on $[0,T_F)$ is
of the order of $1$ minute in all cases considered. The rather small
number of modes in $G$ and the precision of
the Mathematica routines for ODEs presumably make negligible
the numerical errors in the treatment of
\rref{galecomp}. Our analysis assumes this and confuses the numerical solution
of \rref{galecomp} via Mathematica with the exact solution $(\ga_k,\be_k)$ ($k \in G$).
\parn
(ii) The next step is to determine the growth and error estimators
for $\bu_G$. Our attention is focused on the tautological growth
estimators
\beq \cald_p(t) := \| \ug(t) \|_p = \sqrt{ \sum_{k \in G} |k|^{2 p} \big(|\ga_k(t)|^2 + |\beta_k(t)|^2\big)}
\quad (p=3,4,5,6)~ \label{cald34} \feq
and on the tautological, differential error estimators
\beq \ep_p(t) :=
\| \bee(\ug)(t) \|_p = \sqrt{\sum_{k \in dG} |k|^{2 p} \big(|\rho_k(t)|^2 + |\sigma_k(t)|^2 \big)}~
\quad (p=3,5)
\label{nora3} \feq
with $dG$ determined by $G$ and $\rho_k, \sigma_k$ determined by the components $\gamma_k$ and $\beta_k$
$(k\in G)$ according to Eq. \rref{630}. The choice of the orders $p$
in Eqs. \rref{cald34} \rref{nora3} will be clarified by the subsequent item (iii).
\parn
In the case under analysis, the initial datum $\bu_0$ belongs to
the Galerkin subspace $\HD{G}$, so the datum error
$-(1 - \bEG) \bu_0$ vanishes, and the corresponding estimators can be set
to zero:
\beq \delta_p = 0 \quad (p \in \reali)~. \feq
For the computation of
$\cald_p$ ($p=3,...,6$) and $\ep_p$ ($p=3,5$) via
Mathematica, we have
used a two-step approach. First of all,
we have used Eqs. \rref{cald34} \rref{nora3}
at a grid of about $40$ points in the interval $[0,T_F)$;
then we have asked Mathematica to interpolate
the results.
The calculation of $\ep_3$ and $\ep_5$ at the above mentioned grid of
points in $[0,T_F)$ is the most expensive part of the
present scheme in terms of time, since it requires
$15$ minutes approximately for each one of the two
estimators; all the other computations
for the present item (ii) are performed within few seconds.
\parn
In the sequel, it is assumed that the interpolating
functions produced in this way can be
confused with the actual functions $\cald_p,\ep_p$.
\parn
(iii) Having the necessary estimators, we can pass to
the control inequalities. In particular, the control Cauchy
problem of order $3$ reads: find
$\calr_3 \in C^1([0,T_c), \reali)$ ($0 < T_c \leqs T_F$) such that
\ben
\frac{d \calr_3}{dt} & = &- \mu \calr_3 + (\Gp_3 \cald_3 + \Kp_3   \cald_{4})\calr_3
+ \Gp_3 \calr_3^2 + \epsilon_3, \label{contreq3}\\
\calr_3 (0) & = &  0  \label{contreq03}
\een
with $\Kp_3, \Gp_3$ as in Eq.\rref{k3g3}
(see Remark \ref{remcont}, here used
with $n=3$ and $\delta_n = 0$).
The numerical solution of problem
\rref{contreq3} \rref{contreq03}
is performed almost instantaneously
by Mathematica, which easily detects
the possible blow-up of $\calr_3$ at a time $T_c < T_F$.
In the sequel we assume that the numerical
solution provided by Mathematica can be confused
with the actual solution
$\calr_3$ of \rref{contreq3} \rref{contreq03},
even for what concerns its domain. \parn
On the grounds of Proposition
\ref{mainsmo}, the solution $\bu$ of the MHD
Cauchy problem \rref{cmhd} is granted to exist at least up to time
$\Tc$, and to fulfill
\beq \| \bu(t) - \bu_G(t) \|_3 \leqs \Rr_3(t) \quad \mbox{for $t \in [0,\Tc)$}.
\feq
Let us also recall Corollary \ref{coro2} which grants the following: if
\beq (\cald_3 + \Rr_3)(t_1) \leqs {\mu \over \Gp_3}
\qquad \mbox{for some $t_1 \in [0,\Tc)$}~, \label{ipr3} \feq
the solution $\bu$ of \rref{cmhd} exists up to
$T = + \infty$, and
\beq \| \bu(t) \|_3 \leqs
{(\cald_3 + \Rr_3)(t_1) e^{-\mu (t-t_1)} \over 1 - \Gp_3 (\cald_3 + \Rr_3)(t_1) \, e_{\mu}(t-t_1)}
~\mbox{for $t \in [t_1, +\infty)$}~.
\label{tes13} \feq
($e_\mu$ as in Eq. \rref{enu}).
Once $\Rr_3$ is known, using item (ii) of Proposition
\ref{mainsmo} one could construct for each $p > 3$
a function $\Rr_p$ with the same domain $[0,T_c)$,
giving a bound on $\| \bu(t) - \bu_G(t) \|_p$.
For example, if $p=5$ we have
a function $\Rr_5 \in C^1([0,T_c),\reali)$
fulfilling as equalities the relations \rref{cont1npsmo} \rref{cont2npsmo}
with $p=5$ and $n=3$, i.e.:
\beq {d \RR_5 \over d t} = - \mu \RR_5
+ (\Gp_5 \cald_5 + \Kp_5 \cald_{6}+ \Gp_{5 3}\RR_3) \RR_5  + \ep_5
~\mbox{everywhere on $[0,T_c)$}, \label{cont153smo} \feq
\beq \RR_5(0) = 0~, \label{cont253smo} \feq
with $\Kp_5, \Gp_5, \Gp_{5 3}$ as in Eq.\rref{k3g3}.
For the function $\RR_5$ we have an integral representation,
provided by Eqs. \rref{soleq}\rref{defap}; however, the direct
numerical solution of the Cauchy problem
\rref{cont153smo} \rref{cont253smo} via Mathematica
is almost instantaneous, and it has been preferred to
the computation of the integrals in \rref{soleq}\rref{defap}.
\parn
Given the numerical solution $\RR_5$ of Eqs. \rref{cont153smo} \rref{cont253smo} (that we confuse
with the exact solution),
we have the bound
\beq \| \bu(t) - \bu_G(t) \|_5 \leqs \Rr_5(t)~\mbox{for $t \in [0,\Tc)$.} \label{bou5} \feq
In the subsequent applications of the scheme (i)(ii)(iii) for
several values of $\mu$, the graph of $\Rr_5$ is reported only
in a case with rather large $\mu$, in which $\Rr_5(t)$ is
small with respect to $\cald_5(t) = \| \bu_G(t) \|_5$; this fact
makes the bound \rref{bou5} interesting. In the other cases
considered, $\Rr_5(t)$ is sensibly larger than $\cald_5(t)$;
this makes the bound \rref{bou5} less interesting, so
the graph of $\Rr_5$ is not so useful.
\vskip 0.2cm \noindent
Let us pass to exemplify the procedure (i)(ii)(iii) for
some values of $\mu$.
\vskip 0.1cm \noindent
\textbf{Case $\boma{\mu=20}$.} System
\rref{galecomp} for the unknowns $\ga_k, \be_k$
($k \in G$) has been integrated
on a time interval of length $T_F = 0.5$.
Figures \ref{fig1a}-\ref{fig1d} report, as
examples, the graphs of $|\ga_k(t)|$, $|\be_k(t)|$
for $k=(0,1,0)$ and $k = (1,1,0)$ (where $|z| :=
\sqrt{\sum_{i=1}^3 |z_i|^2}$ is the standard
$\complessi^3$ norm).
Figures \ref{fig1e}-\ref{fig1h} report the graphs
of the estimators $\cald_p, \epsilon_p$ for $p=3,5$
(the graphs of $\cald_p$ for $p=4,6$ are omitted
just for brevity). \parn
The solution $\calr_3$ of the control Cauchy problem
\rref{contreq3} \rref{contreq03} is found
to exist on the whole interval $[0,T_F) = [0,0.5)$;
its graph is given by Figure \ref{fig1i}. \parn
It turns out that condition \rref{ipr3}
$(\cald_3 + \Rr_3)(t_1) \leqs {\mu/\Gp_3}$
is fulfilled for any $t_1 \in [0.01,0.5)$;
this ensures that the solution of
the MHD Cauchy problem \rref{cmhd}
is global ($T=+\infty$) and decays
exponentially as indicated by
\rref{tes13}. For example, let us
write down the estimate \rref{tes13}
choosing $t_1 = 0.25$; with appropriate roundings we have
$(\cald_3 + \Rr_3)(0.25) = 0.186$ and
$\Gp_3 (\cald_3 + \Rr_3)(0.25) = 0.115$, so
the cited equation gives
\beq \| \bu(t) \|_3 \leqs
{0.186 e^{-20 (t- 0.25)} \over 1 - 0.115 \, e_{20}(t-0.25)}
~\mbox{for $t \in [0.25, +\infty)$}~.
\label{tes1venti} \feq
($e_{20}$ as in \rref{enu}). Figure \ref{fig1j} gives the graph of $\RR_5(t)$
for $t \in [0,0.5)$. Of course, we have
\beq \| \bu(t) - \bu_G(t) \|_3 \leqs \RR_3(t),~
\| \bu(t) - \bu_G(t) \|_5 \leqs \RR_5(t)
\quad \mbox{for $t \in [0,0.5)$}. \feq
The first of these estimates is
certainly interesting on the whole
interval $[0,0.5)$, where
$\RR_3(t)$ is always much smaller than
$\cald_3 (t) := \| \bu_G(t)\|_3$
($\RR_3(t) < \cald_3(t)/100$ for $t \in [0,0.5)$).
Concerning the second estimate, it
should be pointed out that $\RR_5(t)$
is smaller than $\cald_5(t) := \| \bu_G(t)\|_5$
on the whole interval $[0,0.5)$, and sensibly
smaller on a shorter interval
($\RR_5(t) < \cald_5(t)/2$ for $t \in [0,0.5)$
and $\RR_5(t) < \cald_5(t)/10$ for $t \in [0,0.047)$).
\vskip 0.1cm \noindent
\textbf{Case $\boma{\mu=6}$.} System
\rref{galecomp} for $\ga_k, \be_k$
($k \in G$) has been integrated
on a time interval of length $T_F = 2$.
Figures \ref{fig2a}-\ref{fig2d} give
the graphs of $|\ga_k(t)|$, $|\be_k(t)|$
for $k=(0,1,0)$ and $k = (1,1,0)$.
Figures \ref{fig2e}-\ref{fig2f} report the graphs
of the estimators $\cald_3$, $\epsilon_3$. \parn
The solution $\calr_3(t)$ of the control Cauchy problem
\rref{contreq3} \rref{contreq03} is found
to exist on the whole interval $[0,T_F) = [0,2)$;
its graph is given by Figure \ref{fig2g}.
Condition \rref{ipr3}
$(\cald_3 + \Rr_3)(t_1) \leqs {\mu/\Gp_3}$
is fulfilled for all $t_1 \in [0.32,2)$;
this ensures that the solution $\bu$ of
the MHD Cauchy problem \rref{cmhd}
is global ($T=+\infty$) and decays
exponentially as indicated by
\rref{tes13}. For example, let us
write down the estimate \rref{tes13}
choosing $t_1 = 1$; with appropriate roundings we have
$(\cald_3 + \Rr_3)(1) = 1.09$ and
$\Gp_3 (\cald_3 + \Rr_3)(1) = 0.68$, so
the cited equation gives
\beq \| \bu(t) \|_3 \leqs
{1.09 e^{-6 (t- 1)} \over 1 - 0.68 \, e_{6}(t-1)}
~\mbox{for $t \in [1, +\infty)$}~.
\label{tes1ventibis} \feq
($e_{6}$ as in \rref{enu}). We have
\beq \| \bu(t) - \bu_G(t) \|_3 \leqs \RR_3(t)
\quad \mbox{for $t \in [0,2)$}. \feq
This estimate is especially interesting when
$\RR_3(t)$ is sensibly smaller than
$\cald_3 (t) := \| \bu_G(t)\|_3$. This does not
hold on the whole interval $[0,2)$, but it
is true on shorter intervals: for example, $\RR_3(t) < \cald_3(t)/10$
for $t \in [0,0.11]$.
\vskip 0.1cm \noindent
\textbf{Case $\boma{\mu=5}$.}
System
\rref{galecomp} for $\ga_k, \be_k$
($k \in G$) has been integrated
on a time interval of length $T_F = 2$.
Figures \ref{fig3a}-\ref{fig3d} give
the graphs of $|\ga_k(t)|$, $|\be_k(t)|$
for $k=(0,1,0)$ and $k = (1,1,0)$.
Figures \ref{fig3e}-\ref{fig3f} report the graphs
of the estimators $\cald_3$, $\epsilon_3$. \parn
The solution $\calr_3$ of the control Cauchy problem
\rref{contreq3} \rref{contreq03} has
domain $[0,\Tc)$ where $\Tc = 0.3238...$; $\calr_3(t)$
is found to diverge for $t \to \Tc$. Figure \ref{fig3g}
contains the graph of $\calr_3$. \parn
Due to the features of our general scheme,
the solution $\bu$ of the MHD Cauchy problem \rref{cmhd}
is granted to exist (at least) up to time $\Tc$.
We have
\beq \| \bu(t) - \bu_G(t) \|_3 \leqs \RR_3(t)
\quad \mbox{for $t \in [0,\Tc) = [0,0.3238...)$}; \feq
this inequality is especially interesting
in the smaller interval $[0,0.1]$, where
$\RR_3(t)$ is sensibly smaller than $\cald_3 (t) := \| \bu_G(t)\|_3$
($\RR_3(t) < \cald_3(t)/10$ for $t \in [0,0.1]$).
\vskip 0.1cm \noindent
\textbf{Case $\boma{\mu=0}$.}
Again, system
\rref{galecomp} for $\ga_k, \be_k$
($k \in G$) has been integrated
on a time interval of length $T_F = 2$.
Figures \ref{fig4a}-\ref{fig4d} give
the graphs of $|\ga_k(t)|$, $|\be_k(t)|$
for $k=(0,1,0)$ and $k = (1,1,0)$.
Figures \ref{fig4e}-\ref{fig4f} report the graphs
of the estimators $\cald_3$, $\epsilon_3$. \parn
The solution $\calr_3$ of the control Cauchy problem
\rref{contreq3} \rref{contreq03} has
domain $[0,\Tc)$ where $\Tc = 0.1211...$, and diverges
for $t \to \Tc$. The graph of $\calr_3$ is presented in
Figure \ref{fig4g}. \parn
The solution $\bu$ of the MHD Cauchy problem \rref{cmhd}
is granted to exist up to time $\Tc$.
We have
\beq \| \bu(t) - \bu_G(t) \|_3 \leqs \RR_3(t)
\quad \mbox{for $t \in [0,\Tc) = [0,0.1211...)$}; \feq
this inequality is especially interesting
in the smaller interval $t \in [0,0.067]$, where
$\RR_3(t)$ is sensibly smaller than $\cald_3 (t) := \| \bu_G(t)\|_3$
($\RR_3(t) < \cald_3(t)/10$ for $t \in [0,0.067]$)
\vskip 0.1cm \noindent
\textbf{Other cases.} We have performed computations similar
to those described before even for $\mu=10$ and $\mu=3$.
As in the cases $\mu=20$ and $\mu=6$, for $\mu=10$ we can
grant global existence for the exact solution $\bu$
of the MHD Cauchy problem \rref{cmhd}. As in the cases
$\mu=5$ and $\mu=0$, for $\mu=3$ we can grant existence
of the solution $\bu$ only on a finite interval; in fact
the solution $\calr_3(t)$ of the control Cauchy problem
\rref{contreq3} \rref{contreq03} diverges for $t \to
\Tc = 0.1853...$, so we can ensure existence of $\bu$
only up to $\Tc$.
\vskip 0.6cm
\noindent
\textbf{\Large{Acknowledgments}}
\vskip 0.2cm \noindent
L.P. acknowledges support from: INdAM, Gruppo Nazio\-nale per la Fisica Matematica;
Istituto Nazionale di Fisica Nucleare; MIUR, PRIN 2010 Research Project
``Geometric and analytic theory of Hamiltonian systems in finite and infinite dimensions'';
Universit\`{a} degli Studi di Milano.


\vfill \eject \noindent
\appendix
\section{Appendix. Proof of Eqs. \rref{basic2np}- \rref{weknowg}}
\label{appehel}
In this appendix we frequently use the following inequalities, holding for
all $\alpha,\beta,\gamma,\delta$ $\in \reali$:
\beq \alpha \gamma + \beta \delta \leqs \sqrt{\alpha^2 + \beta^2} \sqrt{\gamma^2 + \delta^2}~, \label{elem1} \eeq
\beq \gamma + \delta \leqs \sqrt{2} \sqrt{\gamma^2 + \delta^2}~. \label{elem2} \eeq
Eq. \rref{elem1} is just the Schwartz inequality for the standard inner product of $\reali^2$;
Eq. \rref{elem2} is the specialization of \rref{elem1} to the case $\alpha=\beta=1$. We will also
use the parallelogram law
\beq \| x + y \|^2 + \| x - y \|^2 = 2 \| x \|^2 + 2 \| y \|^2~, \label{parall} \feq
holding for all elements $x,y$ of any Hilbert space with norm $\|~\|$.
\begin{prop}   \label{propbasic2np}
\textbf{Proposition}.
Consider two reals $p \geqs n>d/2$. For $\bv \in \HD{p}$, $\bw \in \HD{p+1}$
one has the inequality
\beq
\normP{\bP (\bv,\bw)} \, \, \leqs
\frac{1}{2} (\sqrt{2} K_{pn}) (\normP{\bv} \normpnb{\bw}+ \normnb{\bv} \normPp{\bw} ),
\label{onehask} \eeq
where $K_{pn}$ is a constant fulfilling Eq. \rref{basineqa}. Thus, Eqs.
\rref{basic2np} \rref{weknowk} hold.
\end{prop}
\textbf{Proof.} Let us fix $p \geqs n > d/2$ and
$\bv=(v,b) \in \HD{p}$, $\bw=(w,c) \in \HD{p+1}$.
The definition \rref{debp} of $\bP$ gives
\beq
\normP{ \bP (\bv , \bw)}= \sqrt{\normP{\P(v,w) - \P(b,c)}^2 + \normP{\P(v,c)-\P(b,w)}^2}~.
\eeq
Let us note that
\begin{align}
&\normP{\P(v,w) - \P(b,c)} \leqs \normP{\P(v,w)}+\normP{\P(b,c)}  \nonumber \\
&\leqs \frac{1}{2}K_{pn}(\normP{v}\normp{w}+\normn{v}\normPp{w}
+ \normP{b}\normp{c}+\normn{b}\normPp{c})  \nonumber \\
&\leqs \frac{1}{2}K_{pn}(\sqrt{\normP{v}^2+\normP{b}^2}\sqrt{\normp{w}^2 + \normp{c}^2}
+\sqrt{\normn{v}^2+\normn{b}^2}\sqrt{\normPp{w}^2 + \normPp{c}^2}\,) \nonumber \\
&= \frac{1}{2}K_{pn}(\normP{\bv}\normp{\bw}+\normn{\bv} \normPp{\bw}).  \label{dis1}
\end{align}
In the above chain of relations, to go from the first to the second line we have used the inequality \rref{basineqa};
to go from the second to the third line, after exchanging the order of summands we have used twice
Eq. \rref{elem1}. \par \noindent
Similarly, we obtain
\beq
\normP{\P(v,c) - \P(b,w)} \leqs \frac{1}{2}K_{pn}(\normP{\bv}\normp{\bw}+\normn{\bv} \normPp{\bw}).  \label{dis2}
\eeq
Consequently, from (\ref{dis1}) and (\ref{dis2}), we obtain
\begin{align}
&\normP{ \bP ( \bv , \bw )} =\sqrt{\normP{\P(v,w) - \P(b,c)}^2 + \normP{\P(v,c)-\P(b,w)}^2} \nonumber \\
& \leqs \sqrt{\frac{K_{pn}^2}{4}(\normP{\bv}\normp{\bw}+\normn{\bv} \normPp{\bw})^2
+\frac{K_{pn}^2}{4}(\normP{\bv}\normp{\bw}+\normn{\bv} \normPp{\bw})^2} \nonumber \\
&=\frac{K_{pn}}{\sqrt{2}}(\normP{\bv}\normp{\bw}+\normn{\bv} \normPp{\bw}), \nonumber
\end{align}
which yields the inequality (\ref{onehask}).  \fine
\begin{prop} \label{propkatop2np}
\textbf{Proposition.}
Consider two reals $p \geqs n>d/2 +1$. For $\bv \in \HD{p}$, $\bw \in \HD{p+1}$
one has the inequality
\beq
\vert \la \bP ( \bv,\bw) \vert \bw \ra _{{p}}\vert \leqs
\frac{1}{2} (\sqrt{2} G_{pn}) (\normP{\bv} \normnb{\bw}+\normnb{\bv}\normP{\bw})\normP{\bw},
\label{tesi0} \eeq
where $G_{p n}$ is a constant fulfilling \rref{katineqa}.
Thus, Eqs. \rref{katop2np} \rref{weknowg} hold.
\end{prop}
\noindent
\textbf{Proof.}
In the sequel $p \geqs n > d/2+1$ and $\bv=(v,b) \in \HD{p}$, $\bw=(w,c) \in \HD{p+1}$
are fixed; we proceed in several steps.
\parn
\textsl{Step 1. One has}
\beq  \label{tesi1}
\la \bP (\bv,\bw) \vert \bw \ra _{p}=
\la \P(v,w)| w \ra_p + \la \P(v,c)|c \ra_p + \eeq
$$ - \frac{1}{2} \la \P(b,w + c)| w + c \ra_p + \frac{1}{2} \la \P(b,w - c)| w - c \ra_p~. $$
To prove this, we note that Eq. \rref{debp} for $\bP$ implies
\beq \label{1pp}
\begin{split}
\la \bP (\bv,\bw) \vert \bw\ra _{p}=\la( \P(v,w)-\P(b,c) , \P(v,c)-\P(b,w)) \vert (w,c) \ra _{p} \\
=\la \P(v,w)\vert w\ra _{p} - \la \P(b,c)\vert w\ra _{p}+\la \P(v,c)\vert c\ra _{p}
-\la  \P(b,w)\vert c\ra _{p}~.
\end{split}
\eeq
On the other hand, by elementary manipulations relying on the bilinearity of
$\P$ and $\la~|~\ra_{p}$ we get
\beq \la \P(b,c) | w \ra_p + \la \P(b,w) | c \ra_p =
\frac{1}{2} \la \P(b,w + c)| w + c \ra_p - \frac{1}{2} \la \P(b,w - c)| w - c \ra_p
\eeq
and inserting this result into \rref{1pp} we get the thesis
\rref{tesi1}. \parn
\textsl{Step 2. One has}
\beq  \label{tesi2}
|\la \bP (\bv,\bw) \vert \bw \ra _{p} | \leqs
|\la \P(v,w)| w \ra_p| + |\la \P(v,c)|c \ra_p| + \eeq
$$ + \frac{1}{2} |\la \P(b,w + c)| w+c \ra_p| + \frac{1}{2} |\la \P(b,w-c)| w-c \ra_p |~. $$
This is an obvious consequence of \rref{tesi1}. \parn
\textsl{Step 3. One has}
\beq |\la \P(v,w)| w \ra_p| + |\la \P(v,c)| c \ra_p| \leqs
\frac{1}{2} G_{p n} \| v \|_p \| \bw \|_n \| \bw \|_p + \frac{1}{2} G_{p n} \| v \|_n \| \bw \|^2_p~.
\label{tesi3} \eeq
In fact, due to \rref{katineqa},
\beq |\la \P(v,w)| w \ra_p| \leqs \frac{1}{2} G_{p n} (\| v \|_p \| w \|_n + \| v \|_n \| w \|_p) \| w \|_p
\label{dains} \eeq
$$ = \frac{1}{2} G_{p n} (\| v \|_p \| w \|_n \| w \|_p + \| v \|_n \| w \|^2_p)~; $$
one treats similarly the term $|\la \P(v,c)| c \ra_p|$, so
\beq |\la \P(v,w)| w \ra_p| + |\la \P(v,c)| c \ra_p| \eeq
$$ \leqs \frac{1}{2} G_{p n} (\| v \|_p \| w \|_n \| w \|_p + \| v \|_n \| w \|^2_p) +
\frac{1}{2} G_{p n} (\| v \|_p \| c \|_n \| c \|_p + \| v \|_n \| c \|^2_p) $$
$$ = \frac{1}{2} G_{p n} \| v \|_p ( \| w \|_n \| w \|_p + \| c \|_n \| c \|_p )
+ \frac{1}{2} G_{p n} \| v \|_n (\| w \|^2_p + \| c \|^2_p)~. $$
On the other hand, due to \rref{elem1}
\beq \| w \|_n \| w \|_p + \| c \|_n \| c \|_p \leqs
\sqrt{ \| w \|^2_n + \| c \|^2_n} \sqrt{ \| w \|^2_p + \| c \|^2_p}
= \| \bw \|_n \| \bw \|_p~, \label{ins1}\feq
while
\beq \| w \|^2_p + \| c \|^2_p = \| \bw \|^2_p~; \label{ins2} \feq
inserting Eqs. \rref{ins1} \rref{ins2} into \rref{dains} we
get the thesis \rref{tesi3}.
\parn
\textsl{Step 4. One has}
\beq \frac{1}{2} |\la \P(b,w+c)| w+c \ra_p| + \frac{1}{2} |\la \P(b,w-c)| w - c\ra_p |
\label{tesi4} \eeq
$$ \leqs \frac{1}{2} G_{p n} \| b \|_p \| \bw \|_n \| \bw \|_p + \frac{1}{2} G_{p n} \| b \|_n \| \bw \|^2_p~. $$
In fact, using the inequality \rref{katineqa} for each one of the above two
terms we get
\beq \frac{1}{2} |\la \P(b,w+c)|w+c \ra_p| + \frac{1}{2} |\la \P(b,w-c)| w-c \ra_p | \eeq
$$ \leqs \frac{1}{4} G_{p n} (\| b \|_p \| w + c \|_n + \| b \|_n \| w + c \|_p) \| w + c \|_p +
\frac{1}{4} G_{p n} (\| b \|_p \| w - c \|_n + \| b \|_n \| w - c \|_p) \| w - c \|_p $$
$$ = \frac{1}{4} G_{p n} \| b \|_p ( \| w + c \|_n \| w + c \|_p + \| w - c \|_n \| w - c \|_p )
+ \frac{1}{4} G_{p n} \| b \|_n ( \| w + c \|^2_p + \| w - c \|^2_p )~; $$
from here and from Eq. \rref{elem1} we infer
\beq \frac{1}{2} |\la \P(b,w + c)| w + c \ra_p| + \frac{1}{2} |\la \P(b,w - c)| w - c \ra_p | \label{dains2} \eeq
$$ \leqs \frac{1}{4} G_{p n} \| b \|_p \sqrt{ \| w + c \|^2_n + \| w - c \|^2_n}
\sqrt{ \| w + c \|^2_p + \| w - c \|^2_p} + \frac{1}{4} G_{p n} \| b \|_n ( \| w + c \|^2_p + \| w - c \|^2_p )~. $$
On the other hand, the parallelogram law \rref{parall} for the Hilbert spaces $\HM{n},\HM{p}$ gives
\beq \| w + c \|^2_n + \| w - c \|^2_n = 2 \| w \|^2_n  + 2 \| c \|^2_n = 2 \| \bw \|^2_n~, \label{parahn} \feq
$$ \| w + c \|^2_p + \| w - c \|^2_p = 2 \| w \|^2_p  + 2 \| c \|^2_p = 2 \| \bw \|^2_p $$
and inserting Eq. \rref{parahn} into \rref{dains2} we get the thesis \rref{tesi4}.
\parn
\textsl{Step 5. The inequality \rref{tesi0} holds (so the proof is concluded).}
In fact, from Eqs. \rref{tesi2} \rref{tesi3} \rref{tesi4} we get:
\beq |\la \bP (\bv,\bw) \vert \bw \ra _{p} | \label{weget} \feq
$$ \leqs \frac{1}{2} G_{p n} (\| v \|_p + \| b \|_p)
\| \bw \|_n \| \bw \|_p + \frac{1}{2} G_{p n} (\| v \|_n + \| b \|_n) \| \bw \|^2_p~.
$$
On the other hand, Eq. \rref{elem2} gives
\beq \| v \|_p + \| b \|_p \leqs
\sqrt{2} \sqrt{ \| v \|^2_p + \| b \|^2_p} = \sqrt{2}\,\| \bv \|_p~, \feq
$$ \| v \|_n + \| b \|_n \leqs
\sqrt{2} \sqrt{ \| v \|^2_n + \| b \|^2_n} = \sqrt{2}\,\| \bv \|_n~,
$$
and inserting these inequalities into \rref{weget} we obtain the thesis \rref{tesi0}. \fine
\vfill \eject \noindent

\vfill \eject 


\noindent
\begin{thebibliography}{90}

\bibitem{Che} S.I. Chernyshenko, P. Constantin, J.C. Robinson, E.S. Titi,
\textsl{A posteriori regularity of the three-dimensional Navier-Stokes
equations from numerical computations}, J. Math. Phys. \textbf{48} (2007), 065204/1-10.

\bibitem{DR} M. Dashti, J.C. Robinson, \textsl{An a posteriori condition on the numerical
approximation of the Navier-Stokes equations for the existence of a strong solution},
SIAM J. Numer. Anal. \textbf{46} (2008), 3136-3150.

\bibitem{Rob} J.C. Robinson, W. Sadowski, \textsl{Numerical verification of
regularity in the three-dimensional Navier-Stokes equations for bounded sets of
initial data}, Asymptot. Anal. \textbf{59} (2008), 39-50.

\bibitem{appeul} C. Morosi and L. Pizzocchero,
\textsl{On approximate solutions of the incompressible Euler and Navier-Stokes equations}, Nonlinear Analysis,
\textbf{75} (2012), 2209-2235.

\bibitem{reylarge} C. Morosi, M. Pernici, L. Pizzocchero,
\textsl{Large order Reynolds expansions for the Navier-Stokes equations},
Appl. Math. Letters \textbf{49} (2015), 58-66.

\bibitem{hyp2012} C. Morosi, M. Pernici, L. Pizzocchero,
\textsl{A posteriori estimates for Euler and Navier-Stokes equations},
in ``Hyperbolic Problems:
Theory, Numerics and Applications. Proceedings of the XIV International Conference
held in Padova (June 25-29, 2012)'', edited by F. Ancona, A. Bressan, P. Marcati, A. Marson,
AIMS Series on Applied Mathematics \textbf{8} (2014), 847-855.

\bibitem{Mor15} C. Morosi, L. Pizzocchero,
\textsl{Smooth solutions of the Euler and Navier-Stokes equations
from the a posteriori analysis of approximate solutions}, Nonlinear Analysis, \textbf{113} (2015), 298-308.

\bibitem{BlomNol} D. Bl\"omker, C. Nolde, J. Robinson,
\textsl{Rigorous Numerical Verification of Uniqueness and Smoothness in a Surface Growth Model},
Journal of Mathematical Analysis and Applications \textbf{429}(1) (2015), 311-325.

\bibitem{BlomRom}
D. Bl\"omker, M. Romito,
\textsl{Stochastic PDEs and lack of regularity.
(A surface growth equation with noise: existence, uniqueness, and blow-up)},
 Jahresbericht der Deutschen Mathematiker-Vereinigung, \textbf{117}(4) (2015), 233-286.

\bibitem{Nol}
C. Nolde, ``Global Regularity and Uniqueness of Solutions
in a Surface Growth Model Using Rigorous A-Posteriori Methods'',
ISBN: 978-3-8325-4453-9,
Logos Verlag Berlin (2017).

\bibitem{BlomKam}
D. Bl\"omker, M. Kamrani,
\textsl{Numerically Computable A Posteriori-Bounds for SPDEs},
arXiv:1702.01347v3 [math.NA] 13 Nov 2017, to appear in BIT Numerical Mathematics

\bibitem{Mor12b} C. Morosi and L. Pizzocchero, \textsl{On the constants in a Kato inequality
for the Euler and Navier-Stokes equations}, Comm. on Pure and Applied Analysis, \textbf{11} (2012), 557-586.

\bibitem{Mor13} C. Morosi and L. Pizzocchero, \textsl{On the constants in a basic inequality for
the Euler and Navier-Stokes equations}, Applied Mathematics Letters \textbf{26} (2013), 277-284.

\bibitem{Mor17} C. Morosi, M. Pernici, L. Pizzocchero,
\textsl{New results on the constants in some inequalities for the Navier-Stokes quadratic nonlinearity},
Applied Mathematics and Computation, \textbf{308} (2017), 54-72.

\bibitem{Min06} P. D. Mininni, A. G. Pouquet and D. C. Montgomery,
\textsl{Small-scale structures in three-dimensional magnetohydrodynamic turbulence},
Phys. Rev. Lett., \textbf{97} (2006), 244503/1-4.

\bibitem{Car09} C. Cartes, M.D. Bustamante, A. Pouquet and M.E. Brachet,
\textsl{Capturing reconnection phenomena using generalized Eulerian-Lagrangian
description in Navier-Stokes and resistive MHD}, Fluid Dyn. Res., \textbf{41} (2009), 011404/1-14.



\bibitem{Kat72} T. Kato, \textsl{Nonstationary flows of viscous and ideal fluids in $\mathbb{R}^3$},
J. Funct. Anal., \textbf{9} (1972), 296-305.

\bibitem{CoFo} P. Constantin, C. Foias, ``Navier Stokes equations'', Chicago University Press (1988).


\bibitem{Tem} R. Temam, \textsl{Local existence of $C^\infty$ solutions of the Euler equation
of incompressible perfect fluids}, in ``Turbulence and Navier Stokes equation'',
Proceedings of the Orsay Conference,
Lecture Notes in Mathematics \textbf{565} (1976), 184-193.

\bibitem{BKM} J.T. Beale, T. Kato, A.J. Majda, \textsl{Remarks on the breakdown of smooth solutions
for the 3D Euler equations}, Commun. Math. Phys. \textbf{94} (1984), 61-66.

\bibitem{RSS}
J. C. Robinson, W. Sadowski, R. P. Silva,
\textsl{Lower bounds on blow up solutions of the three-dimensional Navier-Stokes equations
in homogeneous Sobolev spaces}, J. Math. Phys. \textbf{53} (2012), 115618, 15pp.

\bibitem{Lions}
G. Duvaut, J. L. Lions,
\textsl{In\'equations en thermo\'elasticit\'e et ma\-gn\'eto\-hydro\-dynamique},
Arch. Rational Mech. Anal. \textbf{46} (1972), 241-279.

\bibitem{Ser83} M. Sermange and R. Temam, \textsl{Some mathematical questions related to the MHD equations},
Comm. on Pure and Applied Mathematics, \textbf{36} (1983), 635-664.

\bibitem{Sch88} P. G. Schmidt, \textsl{On a magnetohydrodynamic problem of Euler type},
Journal of Diff. Equations, \textbf{74} (1988), 318-335.

\bibitem{Cafl} R. E. Caflisch, I. Klapper, G. Steele, \textsl{Remarks on Singularities, Dimension and Energy
Dissipation for Ideal Hydrodynamics and MHD}, Commun. Math. Phys. \textbf{184} (1997), 443-455.

\bibitem{Fan}
J. Fan, T. Ozawa,
\textsl{Regularity criteria for the magnetohydrodynamic equations with partial viscous terms
and the Leray-$\alpha$- MHD model}, Kinetic $\&$ Related Models \textbf{2}(2) (2009), 293-305.

\bibitem{Fef14}
C. L. Fefferman, D. S. McCormick, J. C. Robinson, J.L. Rodrigo,
\textsl{Higher order commutator estimates and local existence for the non-resistive
MHD equations and related models}, Journal of Functional Analysis \textbf{267} (2014), 1035-1056.

\bibitem{Fef17}
C. L. Fefferman, D. S. McCormick, J. C. Robinson, J.L. Rodrigo,
\textsl{Local Existence for the Non-Resistive MHD
Equations in Nearly Optimal Sobolev Spaces},
Arch. Rational Mech. Anal. \textbf{223} (2017), 677-691.


\bibitem{Petr}
T. Petry, \textsl{On the stability
of the Abramov transfer for differential–algebraic equations of index 1}, SIAM J. Numer. Anal.
\textbf{35} (1998), 201-216.

\bibitem{Las} V. Lakshmikantham, S. Leela, ``Differential and
integral inequalities'', Volume I, Academic Press, New York (1969).

\bibitem{Mitr} D.S. Mitrinovic, J.E. Pecaric, A.M. Fink, ``Inequalities
involving functions and their integrals and derivatives'', Kluwer, Dordrecht (1991).

\bibitem{intlab} S. M. Rump, \textsl{INTLAB}, \texttt{http://www.ti3.tu-harburg.de/rump/intlab/}.

\bibitem{arb} Fredrik Johansson, \textsl{Arb - a C library for arbitrary-precision ball arithmetic},
\texttt{http://arblib.org/}.

\bibitem{Rump} S. M. Rump, \textsl{Verification methods:
rigorous results using floating point arithmetic}, Acta Numerica \textbf{19} (2010), 287-449.

\end{thebibliography}
\end{document}